\documentclass[12pt]{article}
\usepackage{latexsym,amsmath,amssymb}
\newenvironment{demo}[1]{\par\smallskip\par\begin{trivlist}
\item[]{\bf #1}\ }{\end{trivlist}\par\smallskip\par}
\newcommand{\Proof}{\begin{demo}{{\it Proof.\ }}}
\newcommand{\qed}{\end{demo}}
\newcommand{\toy}{\ \rule[0em]{0.5ex}{1.8ex}}
\newcommand{\QED}{\toy\end{demo}}
\newtheorem{thm}{Theorem}[section]
\newtheorem{prop}[thm]{Proposition} 
\newtheorem{cor}[thm]{Corollary}
\newtheorem{lem}[thm]{Lemma}
\newtheorem{remark}[thm]{Remark}

\makeatletter \@addtoreset{equation}{section} \makeatother
\makeatletter
\ifnum\@ptsize=0

\else

\fi
\makeatother

\newcommand{\ve}{\varepsilon}

\newcommand{\la}{\langle}
\newcommand{\ra}{\rangle}

\newcommand{\Lm}{\Lambda}

\newcommand{\del}{\partial}

\newcommand{\nn}{\nonumber}
\newcommand{\de}{\partial}

\begin{document}
\title{ 
A stochastic Taylor-like  expansion in the rough path theory
}
\author{Yuzuru INAHAMA 
\\
\\
Graduate School of Mathematics,
\\
Nagoya University, 
\\
Furocho,  Chikusa-ku, Nagoya, 464-8602,  Japan
\\
E-mail: inahama@math.nagoya-u.ac.jp
}

\date{}
\maketitle 
\raggedbottom
%
%
%
%
%
\begin{center}
{\bf  Abstract}
\end{center}
In this paper we establish a Taylor-like 
expansion in the context of the rough path theory
for a family of
It\^{o} maps indexed by a small parameter.
We treat not only the case that the roughness $p$ satisfies $[p]=2$,
but also the case that $[p] \ge 3$.
As an application, we discuss the Laplace asymptotics 
for It\^{o} functionals of Brownian rough paths.
\vspace{3mm}
%
\section{Introduction and the main result}
Let ${\cal V},{\cal W}$ be real Banach spaces
and let $X:[0,1] \to {\cal V}$ be a nice path in $\cal V$.
Let us consider
the following ${\cal W}$-valued ordinary differential equation (ODE);
\begin{equation}\label{int1}
dY_t = \sigma (Y_t) dX_t, 
\qquad
\text{ with $Y_0=0$.}
\end{equation}
Here, 
$\sigma$ is a nice function from
${\cal W}$ to the space $L({\cal V},{\cal W})$
of bounded linear maps.
The correspondence $X \mapsto Y$ is called the It\^{o} map
and will be denoted by $Y=\Phi (X)$.

In the rough path theory of T. Lyons,
the equation (\ref{int1}) is significantly generalized.
First, the space of geometric rough paths on ${\cal V}$
with roughness $p \ge 1$,
which contains all the nice paths, 
is introduced.
It is denoted by $G\Omega_p ({\cal V})$ and its precise definition
will be given in the next section.
Then, the It\^{o} map $\Phi$ extends to a continuous map
from $G\Omega_p ({\cal V})$
to $G\Omega_p ({\cal W})$.
In particular, when $2<p<3$ and $\dim ({\cal V}), \dim ({\cal W}) <\infty$, 
this  equation (\ref{int1}) corresponds to a stratonovich-type
stochastic differential equation (SDE).
(See Lyons and Qian \cite{lq} for the facts in this paragraph.)

In many fields of analysis it is quite important to 
investigate how the output of a map behaves asymptotically
when the input is given small perturbation.
The Taylor expansion in the calculus is a typical example.
In this paper we investigate the behaviour 
of $\Phi(\ve X +\Lm)$ as $\ve \searrow 0$ 
for a nice path $\Lm$ and $X \in G\Omega_p ({\cal V})$.
Slightly generalizing it, we will consider the asymptotic behaviour
of $Y^{(\ve)}$, which is defined by (\ref{int2}) below,
 as $\ve\searrow 0$;
\begin{equation}\label{int2}
dY^{(\ve)}_t = \sigma (\ve, Y^{(\ve)}_t) \ve dX_t
+b(\ve, Y^{(\ve)}_t) d\Lm_t,
\qquad
\text{ with $Y^{(\ve)}_0=0$.}
\end{equation}
Then, we will obtain an asymptotic expansion as follows;
there exist $Y^0, Y^1, Y^2,\ldots $ such that
\[
 Y^{\ve} \sim Y^0 +\ve Y^1+\cdots + \ve^n Y^n + \cdots
\qquad
\text{as $\ve \searrow 0$.}
\]
We call it a stochastic Taylor-like expansion  around a point $\Lm$.
Despite its name, this is purely real analysis and no probability
measure is involved in the argument.

This kind of expansion in the context of the rough path theory 
was first done by Aida \cite{aida2, aida}
(for the case where coefficients $\sigma, b$ are independent of $\ve$,
$[p]=2$, ${\cal V}, {\cal W}$ are finite dimensional).
Then, Inahama and Kawabi \cite{inah-KB-2} (see also \cite{inah})
extended it to the infinite dimensional case
in order to investigate the Laplace asymptotics for 
the Brownian motion over loop groups. 
(The methods in \cite{aida} and \cite{inah-KB-2}
are slightly different. 
In  \cite{aida}, unlike in \cite{inah-KB-2},
the derivative equation of the given equation
is explicitly used. 
see Introduction of \cite{inah-KB-2}.)

The main result (Theorems \ref{thm.main.estI} and \ref{thm.rem.est})
in this paper is to generalize 
the stochastic Taylor-like expansion in \cite{inah-KB-2}.
The following points are improved:
\begin{enumerate}
\item
The roughness $p$ satisfies $2 \le p <\infty$.
In other words, not only the case $[p]=2$,
but also the case $p \ge 3$ is discussed.
 
\item
The coefficients $\sigma$ and $b$ depend on the small parameter
$\ve >0$.
In other words, we treat not just one fixed It\^{o} map,
but a family of It\^{o} maps indexed by $\ve$.

\item
The base point $\Lm$ of the expansion is a
continuous $q$-variational path for any $1 \le q <2$ with $1/p+1/q >1$.
In \cite{inah-KB-2}, $\Lm$ is a continuous bounded variational path
(i.e., the case $q=1$).

\item
Not only estimates of the first level paths of
$Y^0, Y^1, Y^2,\ldots$,
but also estimates of the higher level paths are given.
\end{enumerate}

The organization of this paper is as follows:
In Section 2, we briefly recall the definition and basic facts
on geometric rough paths.
We also prove simple lemmas on continuous 
$q$-variational paths ($1 \le q<2$).
In the end of this section we prove a few lemmas,
including  an extension of Duhamel's principle,
for later use.

In Section 3, 
we first slightly generalize the local Lipschitz continuity of the
integration map as the integrand varies
(Proposition \ref{thm.int.cont}).
Put simply, the proposition states that the map 
$$
(f ,X) \in   
C_{b,loc}^{[p]+1} ({\cal V},L({\cal V},{\cal W}) ) \times G\Omega_p ({\cal V})
\mapsto \int f(X)dX   \in G\Omega_p ({\cal W})
$$
is continuous.
Here, $C_{b,loc}^{[p]+1}$ denotes the space of 
$[p]+1$-times Fr\'echet differentiable maps whose 
derivatives of order $0,1,\ldots, [p]+1$ are bounded on every bounded sets.
Note that in Lyons and Qian \cite{lq},
the integrand (or the coefficients of ODE) is always fixed.
In the path space analysis,
a path on a manifold is often regarded as a current-valued path.
This generalization is also necessary for such a viewpoint
in the rough path context.

In the latter half of the section,
using the above fact,
we slightly improve Lyons' continuity theorem (also known as 
the universal limit theorem)
when the coefficient of the ODE varies.
(Theorem \ref{thm.loc.Lip} and Corollary \ref{cor.loc.Lip}).
Put simply,
the correspondence
\[
(\sigma ,X, y_0) \in   
C_{M}^{[p]+2} ({\cal V},L({\cal V},{\cal W}) ) 
\times G\Omega_p ({\cal V}) \times {\cal W}
\mapsto Z=(X,Y)  \in G\Omega_p ({\cal V} \oplus {\cal W})
\]
is continuous.
Here, $Z=(X,Y)$ is the solution of (\ref{int1})
with the initial condition replaced with $y_0$
and $C_{M}^{[p]+1} ({\cal V},L({\cal V},{\cal W}) ) $ ($M>0$)
is a subset of $C_{b}^{[p]+1} ({\cal V},L({\cal V},{\cal W}) ) $
(a precise definition is given later).

In Section 4, as we stated above,
we prove the main theorems in this paper 
(Theorems \ref{thm.main.estI} and \ref{thm.rem.est}).

In  Section 5,
as an application of the expansion in Section 4,
we improve the Laplace asymptotics for the Brownian rough path
given in \cite{inah-KB-2}.
In this paper, we are now able to treat the case where the coefficients
of the ODE
are dependent on the small parameter $\ve >0$ 
(see Remark \ref{rem.Laplace}).

\begin{remark}\label{rem.intro}
In Coutin and Qian \cite{cq}
they showed that, when the Hurst parameter is larger than $1/4$,
the fractional Brownian rough paths exist
and the rough path theory is applicable to 
the study of SDEs 
driven by the fractional Brownian motion.
In particular, if the Hurst parameter is between $1/4$ and $1/3$,
the roughness satisfies $[p]=3$ and the third level path
plays a role.

(Recently, Friz and Victoir \cite{fv2} showed existence of a geometric rough path
over a multidimensional Gaussian process provided that its covariance function,
in the sense of two dimensional functions,
 is of finite $p$-variation
with $p<2$.)

Since Millet and Sanz-Sol\'e \cite{ms} proved the large deviation
principle for the fractional Brownian rough paths,
it is natural to guess that the Laplace asymptotics 
as in Theorem \ref{thm.expansion}
for the fractional Brownian rough paths is also true.
This was proved in the author's recent preprint \cite{inah2}.

Note that Baudin and Coutin \cite{bc} proved a similar asymptotic 
problem (the short time asymptotics for finite dimensional,
one fixed differential equation) 
for the fractional Brownian rough paths.
Friz and Victoir \cite{fv} also studied a problem similar to \cite{bc}
for finite dimensional (fractional) Brownian rough paths.
\end{remark}

%
\section{The space of geometric rough paths}
\subsection{Definition}
Let $p \ge 2$ and let ${\cal V}$ be a real Banach space.
In this section we recall the definition of $G\Omega_p({\cal V})$,
the space 
of geometric rough paths over ${\cal V}$.
For details, see Lyons and Qian \cite{lq}.

On the tensor product ${\cal V} \otimes \hat{\cal V}$ of two (or more)
Banach spaces ${\cal V}$ and $\hat{\cal V}$,
various Banach norms can be defined.
In this paper, however, we only consider the projective norm
on ${\cal V} \otimes \hat{\cal V}$.
The most important property of the projective norm is 
the following isometrical isomorphism;
$L({\cal V} \otimes \hat{\cal V}, {\cal W}) \cong 
L^{2} ({\cal V},\hat{\cal V}; {\cal W})$.
Here, the right hand side denotes the space of 
bounded bilinear functional from ${\cal V} \times \hat{\cal V}$
to another real Banach space ${\cal W}$.
(For definition and basic properties of the projective norm,
see Diestel and Uhl \cite{du}.)

For a real Banach space $\cal V$ and $n \in {\mathbb N}=\{1,2,\ldots\}$, 
we set 
$T^{(n)}({\cal V})= {\mathbb R} \oplus {\cal V} \oplus \cdots  
\oplus {\cal V}^{\otimes n}$.
For two elements 
$a=(a^0, a^1, \ldots, a^n), b=(b^0, b^1, \ldots, b^n) \in T^{(n)}({\cal V})$,
the multiplication and the scalar action 
are defined as follows;
\begin{eqnarray*}
a \otimes b
&=&
 \bigl(
a^0b^0, a^1b^0+a^0b^1, a^2b^0+a^1 \otimes b^1+a^0b^2,
\ldots,
\sum_{i=0}^n a^{n-i} \otimes b^i
\bigr),
\\
ra &=& (a^0, r a^1, r^2 a^2, \ldots, r^n a^n),
\qquad r \in {\mathbb R}.
\end{eqnarray*}
Note that, $a \otimes b \neq b\otimes a$ in general. 
The non-commutative algebra $T^{(n)}({\cal V})$ is called the 
truncated tensor algebra of degree $n$.
As usual $T^{(n)}({\cal V})$ is equipped with the direct sum norm.

Let $\triangle =\{ (s,t) ~|~ 0\le s \le t \le 1 \}.$
We say $X= (1, X^1,\ldots,X^{[p]}) :\triangle \to T^{([p])}({\cal V})$ 
is a rough path over ${\cal V}$ of roughness $p$
if it  is continuous and satisfies the following; 
\begin{eqnarray*}
X_{s,u} \otimes X_{u,t}=X_{s,t}
\qquad
\text{ for all $(s,u), (u,t) \in \triangle$,}
\nn\\
\|X^j\|_{p/j}: =  
\Bigl\{
\sup_{  D}
\sum_{i=1}^N |X_{t_{i-1}, t_i  }^j|^{p/j}
\Bigl\}^{j/p} < \infty
\qquad
\text{ for all $j=1,\ldots,[p]$}.
\end{eqnarray*}
Here, $D=\{0=t_0<t_1<\cdots <t_N=1 \}$ 
runs over all the finite partitions of $[0,1]$.
The first identity above is called Chen's identity.
The set of all the  rough paths over ${\cal V}$ of roughness $p$
is denoted by $\Omega_p ({\cal V})$.
The distance on $\Omega_p ({\cal V})$ is defined by
\[
d(X,Y)= \sum_{j=1}^{[p]} \|X^j - Y^j \|_{p/j},
\qquad
\text{  $X,Y \in \Omega_p ({\cal V})$.}
\]
With this distance, $\Omega_p ({\cal V})$ is a complete metric space.
For $X \in \Omega_p ({\cal V})$,
we set 
$\xi(X)=\sum_{j=1}^{[p]} \| X^j \|_{p/j}^{1/j}$.
It is obvious that
$\xi(r X)=|r|\xi(X)$ for $r \in {\mathbb R}$.

Let ${\rm BV} ({\cal V}) =\{ X \in C([0,1],{\cal V}) ~|~  
\text{$X_0=0$ and  $\|X\|_1 <\infty$}  \}$
be the space of continuous, bounded variational paths
starting at $0$.
By using the Stieltjes integral, we can define a rough path 
as follows ($p \ge 2$);
\[
X^{j}_{s,t} :=\int_{ s<t_1 <\cdots <t_j<t} 
dX_{t_1} \otimes \cdots\otimes dX_{t_j},
\qquad
(s,t) \in \triangle, \quad
j=1,\ldots,[p].
\]
This rough path is called be 
the smooth rough path lying above $X \in {\rm BV} ({\cal V})$
and is again denoted by $X$
(when there is no possibility of confusion).
The $d$-closure of the totality of all the smooth rough paths
is denoted by $G\Omega_p ({\cal V})$,
which is called the space of geometric rough paths.
This is a complete metric space.
(If ${\cal V}$ is separable, then  $G\Omega_p ({\cal V})$ is 
also separable,
which can easily be seen from Corollary \ref{lem.inj.q} below.)

\subsection{On basic properties of $q$-variational paths ($1\le q<2$).}
%
%
%

Let $1\le q <2$.
For a real Banach space ${\cal V}$,
set 
\[
C_{0,q}({\cal V})=\{ X \in C([0,1], {\cal V}) \,|\, 
\text{ $X_0=0$ and $\|X\|_q <\infty$} \},
\]
where $\| \,\cdot\, \|_q$ denotes the $q$-variation norm.
When $q=1$, ${\rm BV}({\cal V}) =C_{0,q}({\cal V})$.

Let ${\cal P}=\{0=t_0<t_1<\cdots<t_N=1 \}$ 
be a (finite) partition of $[0,1]$.
We denote by $\pi_{{\cal P}}X$ 
a piecewise linear path associated with  ${\cal P}$
(i.e., $\pi_{{\cal P}}X_{t_i}=X_{t_i}$ 
and $\pi_{{\cal P}}X$ is linear on $[t_{i-1},t_i]$ for all $i$).
The following lemma states that
$\pi_{{\cal P}} : C_{0,q}({\cal V}) \to C_{0,q}({\cal V})$
is uniformly bounded as the partition ${\cal P}$ varies.

\begin{lem}\label{lem.u.bdd}
Let $1 \le q <2$.
Then, there exists a positive constant $c=c_q$ depending only on $q$
such that
$\|\pi_{{\cal P}}  X\|_q \le c \|X\|_q $
for any partition ${\cal P}$ of $[0,1]$.
\end{lem}

\Proof
Fix ${\cal P}=\{0=t_0<t_1<\cdots<t_N=1 \}$.
For a partition 
${\cal Q}=\{0=s_0<s_1<\cdots<s_M=1 \}$, we set
\[
%
S_{{\cal Q}}=
\sum_{k =1}^M |\pi_{{\cal P}} X_{s_i} -\pi_{{\cal P}} X_{s_{i-1}}  |^q .
\]
Suppose that $[t_{i-1}, t_i] \cap {\cal Q}$ contains
three  points (namely, $s_{j-1}<s_{j}<s_{j+1}$).
Then, since $\pi_{{\cal P}}  X$ is linear on $[s_{j-1}, s_{j+1}]$,
we see that
\[
 |\pi_{{\cal P}} X_{s_{j-1}} -\pi_{{\cal P}} X_{s_{j}}  |^q
+
 |\pi_{{\cal P}} X_{s_j} -\pi_{{\cal P}} X_{s_{j+1}}  |^q
\le
 |\pi_{{\cal P}} X_{s_{j-1}} -\pi_{{\cal P}} X_{s_{j+1}}  |^q,
\]
which implies that $S_{{\cal Q}} \le S_{{\cal Q} \setminus \{s_j\}}$.
(The same argument holds if it contains more than three points.)
Therefore, we have only to consider ${\cal Q}$'s
such that $|[t_{i-1}, t_i] \cap {\cal Q}| \le 2$ for all $i=1,\ldots,N$.

Let ${\cal Q}$ be as such.
If 
$[t_{i-1}, t_i] \cap {\cal Q}=\{ s_j \}$,
then define $\hat{s}_j=t_{i-1}$ except if $s_j=1$.
(If so we set $\hat{s}_j=1$.)
If 
$[t_{i-1}, t_i] \cap {\cal Q}=\{ s_j <s_{j+1} \}$,
then define $\hat{s}_j=t_{i-1}$ and $\hat{s}_{j+1}=t_{i}$.
Note that 
$0=\hat{s}_0 \le \hat{s}_1 \le \cdots \le \hat{s}_M=1$.
Some of $\hat{s}_j$'s may be equal.
If so, we only collect distinct $\hat{s}_j$'s and 
call the collection $\hat{{\cal Q}}$.
Noting that $S_{ \hat{{\cal Q}} } \le \|X\|^q_q$
and that 
\[
 |\pi_{{\cal P}} X_{s_j} -\pi_{{\cal P}} X_{\hat{s}_{j}}  |
\le 
|\pi_{{\cal P}} X_{t_i} -\pi_{{\cal P}} X_{t_{i-1}}|
=
| X_{t_i} -  X_{t_{i-1}}|,
\]
if $s_j \in [t_{i-1}, t_i]$,
we see that
\begin{eqnarray*}
S_{{\cal Q}}
&=&
\sum_{j=1}^M |\pi_{{\cal P}} X_{s_j} -\pi_{{\cal P}} X_{s_{j-1}}  |^q 
\\
&\le&
c_q
\Bigl[
\sum_{k =1}^M |\pi_{{\cal P}} X_{s_j} -\pi_{{\cal P}} X_{\hat{s}_{j}}  |^q 
+
\sum_{k =1}^M 
|\pi_{{\cal P}} X_{s_{j-1}} -\pi_{{\cal P}} X_{\hat{s}_{j-1}}  |^q 
+
S_{ \hat{{\cal Q}} }
\Bigr]
\\
&\le&
c'_q  S_{ \hat{{\cal Q}} }
\le c'_q \|X\|^q_q.
\end{eqnarray*}
Taking supremum over such ${\cal Q}$'s, we complete the proof.
\QED

\begin{cor}\label{cor.appro.cq}
Let $1 \le q <q' <2$ and $X \in C_{0,q}({\cal V})$.
Then, 
\[
\lim_{|{\cal P}| \to 0}   
 \| X - \pi_{{\cal P}}X   \|_{q'} =0.
\]
Here, $|{\cal P}|$ denotes the mesh of the partition ${\cal P}$.
\end{cor}

\Proof
This is easy from Lemma \ref{lem.u.bdd} and 
the fact that $\pi_{{\cal P}}X  \to X$ as $|{\cal P}| \to 0$
in the uniform topology.
\QED


For $1\le q <2$ and $X \in C_{0,q}({\cal V})$,
we set 
\[
X^j_{s,t} = 
\int_{s< t_1<\cdots< t_j  <t}  
dX_{t_1} \otimes \cdots \otimes dX_{t_j},
\qquad
j=1,2,\ldots,[p].
\]
Here, the right hand side is the Young integral.
As before, $X=(1,X^1,\ldots,X^{[p]})$ satisfies Chen's identity
and, if $|X_t-X_s| \le \omega(s,t)^{1/q}$ for some control 
function $\omega$ (in the sense of p. 16, Lyons and Qian \cite{lq}),
then $X^j$ is of finite $q/j$-variation 
(i.e., there exists a positive constant $C$ such that
$|X^j_{s,t}| \le C\omega(s,t)^{j/q}$).
So, $C_{0,q}({\cal V}) \subset \Omega_{p} ({\cal V})$.

By Theorem 3.1.3 in \cite{lq}, 
if $X,Y \in C_{0,q}({\cal V})$ satisfy that
\begin{eqnarray*}
|X_t-X_s|,  \, |Y_t-Y_s| &\le& \omega(s,t)^{1/q},
\\
|(X_t-X_s)  -(Y_t-Y_s)| &\le& \ve \omega(s,t)^{1/q}.
\end{eqnarray*}
for some control function,
then there exists a positive constant $C$ depending only on 
$p, q, \omega(0,1)$  
such that
\begin{eqnarray*}
| X^j_{s,t} - Y^j_{s,t}| &\le& C \ve \omega(s,t)^{j/q},
\qquad
\text{ for all $(s,t)\in \triangle$ and $j=1,\ldots,[p]$.}
\end{eqnarray*}
In particular, 
the injection 
$X \in C_{0,q}({\cal V}) \mapsto 
X =(1,X^1,\ldots,X^{[p]}) \in \Omega_{p} ({\cal V})$
is continuous.
%
%
%
Combining this with Corollary \ref{cor.appro.cq},
we obtain the following corollary.
(By taking sufficiently small $q'(>q)$.)
Originally, $G\Omega_{p} ({\cal V})$ is defined as the closure 
of ${\rm BV}({\cal V})$ in $\Omega_{p} ({\cal V})$.
In the the following corollary, we prove that $G\Omega_{p} ({\cal V})$ is
also obtained as the closure 
of $C_{0,q}({\cal V})$.
(The ingredients of Corollary \ref{lem.shift.q} is partially in \cite{fv2}.)
\begin{cor}\label{lem.inj.q}
Let $1 \le q<2$ and $p \ge 2$.
For any $X \in C_{0,q}({\cal V})$, $\pi_{{\cal P}}X \in {\rm BV} ({\cal V})$ 
converges to $X$ in $\Omega_{p} ({\cal V})$ as $|{\cal P}| \to 0$.
In particular, 
we have the following continuous inclusion;
${\rm BV} ({\cal V}) \subset C_{0,q}({\cal V}) 
\subset  G\Omega_{p} ({\cal V})$.
\end{cor}

%

\begin{cor}\label{lem.shift.q}
Let $1 \le q<2$ and $p \ge 2$ with $1/p+1/q >1$.
Let ${\cal V}$ and ${\cal W}$ be real Banach spaces.
Then, the following (1) and (2) hold:
\\
\noindent
{\rm (1).}~
For $X \in G\Omega_p({\cal V})$ and $H \in C_{0,q}({\cal V})$,
the natural shift $X+H \in G\Omega_p({\cal V})$ is well-defined.
Moreover, it is continuous
as a map  from $G\Omega_p({\cal V}) \times C_{0,q}({\cal V})$ to $G\Omega_p({\cal V})$.
\\
\noindent
{\rm (2).}~
For $X \in G\Omega_p({\cal V})$ and $H \in C_{0,q}({\cal W})$,
$(X,H) \in G\Omega_p({\cal V} \oplus {\cal W})$ is well-defined.
Moreover, it is continuous
as a map  from $G\Omega_p({\cal V}) \times C_{0,q}({\cal W})$ 
to $G\Omega_p({\cal V} \oplus {\cal W})$.
\end{cor}

\Proof
The shift as a map 
 from $\Omega_p({\cal V}) \times C_{0,q}({\cal V})$ to $\Omega_p({\cal V})$
is continuous.
(See Section 3.3.2 in \cite{lq}.)
Therefore, we have only to prove that
$X+H \in G\Omega_p({\cal V})$ if  $X \in G\Omega_p({\cal V})$.
However, it is immediately shown from Corollary \ref{lem.inj.q}.

The second assertion can be verified in the same way.
\QED

%
%
Now we consider linear ODEs of the following form:
For a given $L({\cal V},{\cal V})$-valued path $\Omega$, we set
\begin{eqnarray}
dM_t &=& d\Omega_t \cdot M_t, \qquad M_0={\rm Id}_{{\cal V}},
 \label{eq.q.ode.lin}
\\
dN_t &=&  - N_t \cdot d\Omega_t, \qquad N_0={\rm Id}_{{\cal V}}.
 \label{eq.q.ode.lin2}
\end{eqnarray}
Here, $M$ and $N$ are also $L({\cal V},{\cal V})$-valued.
It is well-known that, if $\Omega \in {\rm BV} (L({\cal V},{\cal V}))$,
then unique solutions $M, N$ exist in ${\rm BV} (L({\cal V},{\cal V}))$
and $M_t N_t=N_tM_t={\rm Id}_{{\cal V}}$ for all $t$. 
This can be extended to the case of $q$-variational paths
($1 \le q<2$)
as in the following proposition.
\begin{prop}\label{prop.q.ode}
Let $1 \le q<2$ and $\Omega \in C_{0,q}(L({\cal V},{\cal V}))$.
Then, the unique solutions $M, N$ of (\ref{eq.q.ode.lin})
and (\ref{eq.q.ode.lin2}) exist in 
$C_{0,q} (L({\cal V},{\cal V}))+{\rm Id}_{{\cal V}}$.
It holds that $M_t N_t=N_tM_t={\rm Id}_{{\cal V}}$ for all $t$. 
Moreover, if there exists a control function $\omega$ such that
\begin{eqnarray}
|\Omega_t - \Omega_s|, \,
|\hat{\Omega}_t - \hat{\Omega}_s|  &\le& \omega(s,t)^{1/q},
\nn
\\
|(\Omega_t - \Omega_s)
- (\hat{\Omega}_t - \hat{\Omega}_s)|  &\le& \ve \omega(s,t)^{1/q},
\qquad
(s,t) \in \triangle,
\nn
\end{eqnarray}
then, there exists a constant $C$ depending only on 
$q$ and $\omega(0,1)$ such that
\begin{eqnarray}
|M_t - M_s|, \,
|\hat{M}_t - \hat{M}_s|  &\le& C \omega(s,t)^{1/q},
\label{ineq.M.1}
\\
|(M_t - M_s)
- (\hat{M}_t - \hat{M}_s)|  &\le& C\ve \omega(s,t)^{1/q},
\qquad
(s,t) \in \triangle.
\label{ineq.M.2}
\end{eqnarray}
Similar estimates also hold for $N$.
\end{prop}

\Proof 
By using the Young integration,
set $I^0_t= I^0 (\Omega)_{t} :={\rm Id}_{{\cal V}}$ and, 
for $n=1,2,\ldots$, 
\[
I^n_t= I^n (\Omega)_{t} :=
\int_{0< t_1<\cdots<t_n<t} d\Omega_{t_n} \cdots d\Omega_{t_1}.
\]
If we set 
$m_{t}=\sum_{n=0}^{\infty}I^n(\Omega)_{t}$,
then $t \mapsto  m_{t-t_0}A$ formally satisfies (\ref{eq.q.ode.lin})
with initial condition replaced with $m_{t_0}=A \in L({\cal V},{\cal V})$.
Therefore, we will verify the convergence.
(Note that in our construction of solutions, 
 the ``right invariance''
of the given differential equation (\ref{eq.q.ode.lin}) 
implicitly plays an important role.)

We will prove that,
for all $n \in {\mathbb N}$, $0<T<1$, and 
$(s,t) \in \triangle_{[0,T]}=\{(s,t)~|~ 0\le s \le t\le T  \}$, 
\begin{eqnarray}
|I^n_t - I^n_s | \le K^{n-1} \omega(s,t)^{1/q},
\qquad
\text{ where $K=\omega(0,T)^{1/q} \bigl(1+ 2^{2/q}\zeta(2/q) \bigr)$.}
\label{ineq.I}
\end{eqnarray}
Here, $\zeta$ denotes the $\zeta$-function.
Obviously, (\ref{ineq.I}) holds for $n=1$.

Suppose that (\ref{ineq.I}) holds for $n$.
Recall that 
$
I^{n+1}_t - I^{n+1}_s =  \int_s^t   d\Omega_u I^n_u
=
\lim_{ |{\cal P}| \to 0} S_{{\cal P}},
$
where, $S_{{\cal P}}$ is given by
$
S_{{\cal P}} = \sum_{i=1}^N 
( \Omega_{t_{i}}  -   \Omega_{t_{i-1}} )
I^n_{t_{i-1}} 
$
for a finite partition 
${\cal P}= \{ s=t_0<t_1< \cdots<t_N=t \}$ of $[s,t]$.
It is easy to see that 
\begin{eqnarray}
|S_{\{s,t  \}}| =  | ( \Omega_t - \Omega_s )I^{n}_s| 
\le K^{n-1} \omega(0,s)^{1/q} \omega(s,t)^{1/q}.
\label{ineq.S.shoko}
\end{eqnarray}
Let $t_j \in {\cal P}, (j=1,\ldots,N-1)$. 
Then, we have
\begin{eqnarray}
| S_{{\cal P}} - S_{ {\cal P} \setminus \{t_j\}  }|
&=&
\bigl| 
 ( \Omega_{t_{j}}  -   \Omega_{ t_{j-1} } ) I^n_{ t_{j-1} }
+
  ( \Omega_{t_{j+1}}  -   \Omega_{t_{j}} ) I^n_{t_{j}}
-
 ( \Omega_{t_{j+1}}  -   \Omega_{t_{j-1}} ) I^n_{t_{j-1}}
\bigr|
\nn
\\
&=&
\bigl| 
( \Omega_{t_{j+1}}  -   \Omega_{t_{j}} )(I^n_{t_{j}} - I^n_{t_{j-1}})
\bigr|
\le
K^{n-1} \omega(t_{j-1}  , t_{j+1} )^{2/q}.
\nn
%
\end{eqnarray}
From this and a routine argument,
\begin{eqnarray}
| S_{{\cal P}} - S_{ \{s,t\}  }|
\le
K^{n-1} 2^{2/q}  \zeta(2/q) \omega(s  , t )^{2/q}.
\label{ineq.S.middle}
\end{eqnarray}
From (\ref{ineq.S.shoko}) and  (\ref{ineq.S.middle}),
we see that
$
| S_{{\cal P}}| \le 
K^{n} \omega(s  , t )^{1/q}$.
This implies (\ref{ineq.I}) for $n+1$ (and, hence, for all $n$ by induction).

If $T$ is chosen so that $K <1$, 
then $t \mapsto m_t$ is convergent in $q$-variation topology 
on the restricted interval
and is a solution for  (\ref{eq.q.ode.lin})
which satisfies that
$|m_t - m_s| \le (1-K)^{-1} \omega(s,t)^{1/q}$
for $(s,t) \in \triangle_{[0,T]}$.

Take $0=T_0<T_1< \cdots<T_k$ such that
$\omega(T_{i-1},T_i)^{1/q} \bigl(1+ 2^{2/q}\zeta(2/q) \bigr)=1/2$
for $i=1,\ldots, k-1$ 
and 
$\omega(T_{k-1},T_k)^{1/q} \bigl(1+ 2^{2/q}\zeta(2/q) \bigr) \le 1/2$.
By the superadditivity of $\omega$, 
$k-1 \le 2^q  \bigl(1+ 2^{2/q}\zeta(2/q) \bigr)^q \omega(0,1)$.
Hence, $k$ is dominated by a constant which depends only $q$ and
$\omega(0,1)$.
On each time interval $[T_{i-1}, T_i]$, 
construct $M$ by
$M_t=  m_{t- T_{i-1}}M_{T_{i-1} }$.
By the facts we stated above this is a (global) solution 
of (\ref{eq.q.ode.lin}) with desired estimate (\ref{ineq.M.1}).
It is easy to verify the uniqueness.

Finally, we will prove the local Lipschitz continuity (\ref{ineq.M.2}).
In a similar way as above, we will show by induction that
\begin{eqnarray}
|(I^n_t - I^n_s)- ( \hat{I}^n_t - \hat{I}^n_s) | 
\le \ve n K^{n-1} \omega(s,t)^{1/q},
\qquad
\text{$(s,t) \in \triangle_{[0,T]}$, $n \in {\mathbb N}$.}
\label{ineq.II}
\end{eqnarray}
Obviously, (\ref{ineq.II}) holds for $n=1$.

In the same way as above, we see that
\begin{eqnarray}
|S_{\{s,t  \}} - \hat{S}_{\{s,t  \}}| 
&\le&  
| ( \Omega_t - \Omega_s )  
(I^{n}_s - \hat{I}^{n}_s )|
+ 
|  [ ( \Omega_t - \Omega_s ) 
- ( \hat{\Omega}_t - \hat{\Omega}_s ) ]
\hat{I}^{n}_s|
\nn
\\
&\le& 
\ve n K^{n-1} \omega(0,s)^{1/q} \omega(s,t)^{1/q}
+K^{n-1} \omega(0,s)^{1/q} \ve \omega(s,t)^{1/q}
\nn\\
&\le&
\ve (n+1) K^{n-1} \omega(0,1)^{1/q} \omega(s,t)^{1/q}.
\label{ineq.S.shoko2}
\end{eqnarray}
Let $t_j \in {\cal P}, (j=1,\ldots,N-1)$. 
Then,
\begin{eqnarray}
\lefteqn{
\bigl| 
(S_{{\cal P}} - S_{ {\cal P} \setminus \{t_j\}  } )
-
( \hat{S}_{{\cal P}} - \hat{S}_{ {\cal P} \setminus \{t_j\}  } )
\bigr|
}\nn\\
&=&
\bigl| 
  ( \Omega_{t_{j}}  -   \Omega_{ t_{j-1} } )
I^n_{ t_{j-1} }
+
 ( \Omega_{t_{j+1}}  -   \Omega_{t_{j}} )I^n_{t_{j}}
-
 ( \Omega_{t_{j+1}}  -   \Omega_{t_{j-1}} )I^n_{t_{j-1}}
\nn
\\
&&
\qquad
- 
  ( \hat{\Omega}_{t_{j}}  -  \hat{ \Omega}_{ t_{j-1} } )
\hat{I}^n_{ t_{j-1} }
-
( \hat{\Omega}_{t_{j+1}}  -   \hat{\Omega}_{t_{j}} )
\hat{I}^n_{t_{j}} 
+
 ( \hat{\Omega}_{t_{j+1}}  -   \hat{\Omega}_{t_{j-1}} )
\hat{I}^n_{t_{j-1}} 
\bigr|
\nn
\\
&\le&
\bigl| \Omega_{t_{j+1}}  -   \Omega_{t_{j}} 
\bigr|
\,
\bigl| 
(I^n_{t_{j}} - I^n_{t_{j-1}})
- (\hat{I}^n_{t_{j}} -  \hat{I}^n_{t_{j-1}})
\bigr|
\nn\\
&&
\qquad
+ 
\bigl| (\Omega_{t_{j+1}}  -   \Omega_{t_{j}} )
- (\hat{\Omega}_{t_{j+1}}  -   \hat{\Omega}_{t_{j}} )
\bigr|
\,
\bigl| 
(I^n_{t_{j}} - I^n_{t_{j-1}})
\bigr|
\nn\\
&\le&
\ve (n+1) K^{n-1} \omega(t_{j-1}  , t_{j+1} )^{2/q}.
\nn
%
\end{eqnarray}
From this and a routine argument,
\begin{eqnarray}
\bigl| 
(S_{{\cal P}} - S_{ \{s,t\}  })-(\hat{S}_{{\cal P}} - \hat{S}_{ \{s,t\}  })
\bigr|
\le
\ve (n+1)
K^{n-1} 2^{2/q}  \zeta(2/q) \omega(s  , t )^{2/q}.
\label{ineq.S.middle2}
\end{eqnarray}
From (\ref{ineq.S.shoko2}) and  (\ref{ineq.S.middle2}),
we see that
$
| S_{{\cal P}} -\hat{S}_{{\cal P}} | \le \ve (n+1)
K^{n} \omega(s  , t )^{1/q}$.
This implies (\ref{ineq.II}) for $n+1$ (and, hence, for all $n$ by induction).

In the same way as above, we can prolong the solutions 
and obtain (\ref{ineq.M.2}).
The proof for $N$ is essentially the same. So we omit it.
Take $q' \in (q,1)$ and apply Corollary \ref{cor.appro.cq}.
Then, because of the continuity we have just shown,
we see that $M_tN_t=N_tM_t ={\rm Id}_{{\cal V}}$ for all $t$.
\QED

%
%
%
%
%
From now on we will prove a lemma for Duhamel's principle 
in the context of the rough path theory.
When the operator-valued path $M$ below is of finite variation
and $[p]=2$, 
the principle was checked in \cite{inah-KB-2}.
Here, we will consider the case where $p \ge 2$,
$M$ is of finite $q$-variation
($1\le q <2$) with $1/p+1/q>1$.

We set
\begin{eqnarray}
{\cal C}_q (L({\cal V}, {\cal V}))
&:=&
\Big \{ (M,N) ~\big \vert~
M ,N
\in C_{0,q} \big (L({\cal V},{\cal V}) \big) +{\rm Id}_{{\cal V}}, 
\nonumber \\
&\mbox{ }& \qquad \qquad \qquad \qquad \qquad
M_{t}N_{t}=N_{t}M_{t}={\rm Id}_{{\cal V}}
\mbox{ for }
t \in [0,1]
\Big \}.
\nonumber
\end{eqnarray}
We say $M \in  {\cal C}_q(L( {\cal V},{\cal V} ))$ if 
$(M,M^{-1}) \in {\cal C}_q (L({\cal V},{\cal V}))$ for simplicity.
%

%
%
We define a map ${\bf {\Gamma}}:  
 C_{0,q}( {\cal V})\times  
{\cal C}_q (L( {\cal V},{\cal V})) \to C_{0,q}({\cal V})$
by 
\begin{equation}\label{def.Gamma}
{\bf{\Gamma}}(X,M)_t
=
{\bf{\Gamma}}\big ( X,(M,M^{-1}) \big )_{t}
:= 
M_{t}\int_0^t M_{s}^{-1}dX_s, \qquad  t \in [0,1]
\end{equation}
for $X \in C_{0,q}( {\cal V})$ and $M \in {\cal C}_q (L({\cal V},{\cal V}))$. 
Here, the right hand side is the Young integral.

\begin{lem}\label{lem_henkan}
Let ${\cal V}$ be a real Banach space,
 $p \ge 2$,  $1 \le q<2$ with $1/p+1/q>1$.
Let $\Gamma$ be as above.
Then, we have 
the following assertions:
\vspace{2mm}
\\
{\rm (1).~}
Assume that there exists a control function $\omega$ such that
\begin{eqnarray}
& &
| X^j_{s,t}| \le 
\omega(s,t)^{j/p}, \qquad
j=1,\ldots,[p], 
\label{ineq.Gam1}
\\
& &
| M_t -M_s|_{L( {\cal V},{\cal V})} + | M^{-1}_t -M^{-1}_s|_{L({\cal V},{\cal V})} 
\le  \omega(s,t)^{1/q}
\label{ineq.Gam2}
\end{eqnarray}
hold for all $(s,t) \in \triangle$. 
Then,
\begin{equation}
\big|  {\bf{\Gamma}}(X,M)_{s,t}^j \big| 
\le 
C \omega(s,t)^{j/p},
\quad
j=1,\ldots,[p], \quad (s,t) \in \triangle,
\label{ineq.Gam3}
\end{equation}
where $C$ is a
 positive constant
depending only on $p$, $q$, and $\omega(0,1)$.
\\ 
{\rm (2).~}
${\bf {\Gamma}}$  extends to a continuous map 
from $G\Omega_p ( {\cal V})  \times  
 {\cal C}_q (L({\cal V},{\cal V}))$ to $G\Omega_p ({\cal V})$.
(We denote it again by ${\bf {\Gamma}}$.)
Clearly, ${\bf{\Gamma}}( \ve X,M)  = \ve 
{\bf{\Gamma}}  (  X,M)$
holds
for any $X \in G\Omega_p ({\cal V}),
M\in {\cal C}_q (L({\cal V},{\cal V}))$ and $\ve \in {\mathbb R}$.
\end{lem}

\Proof
The proof is not very difficult. So we give a sketch of proof.
From (\ref{def.Gamma}), we have 
$
{\bf{\Gamma}}(X,M)_t
=
X_t - M_t \int_0^t (d M^{-1}_s) X_s.$
Note that the map that associates $(X,M)$ with the second term 
above is continuous from $G\Omega_p ( {\cal V})  \times  
 {\cal C}_q (L({\cal V},{\cal V}))$ to $C_{0,q} ({\cal V})$.
Using Corollary \ref{lem.shift.q}, we see that 
$(X,M) \mapsto \Gamma(X,M)$ is continuous from $G\Omega_p ( {\cal V})  \times  
 {\cal C}_q (L({\cal V},{\cal V}))$ to $G\Omega_p ({\cal V})$.
\QED

The following corollary is called (the rough path version of) 
Duhamel's principle
and will be used frequently below.  
\begin{cor}\label{cor.duhamel}
Let ${\cal V}$ be a real Banach space,
 $p \ge 2$,  $1 \le q<2$ with $1/p+1/q>1$.
For $\Omega \in C_{0,q}(L({\cal V},{\cal V}))$, 
define $M=M_{\Omega}$ and $N=M^{-1}_{\Omega}$ as in (\ref{eq.q.ode.lin})
and (\ref{eq.q.ode.lin2}).
For this $(M_{\Omega},M^{-1}_{\Omega})$ and $X \in C_{0,q}({\cal V})$ 
define $Y:=\Gamma(X,M_{\Omega})$ as in (\ref{def.Gamma}).
Then, the following (1)--(3) hold:

\noindent
{\rm (1).}~
$Y$ is clearly 
the unique solution of the following ODE (in the $q$-variational sense):
\[
dY_t-(d \Omega_t)\cdot Y_t =dX_t,
\qquad Y_0=0.
\]

\noindent
{\rm (2).}~
Assume that there exists a control function $\omega$ such that
\begin{eqnarray}
& &
| X^j_{s,t}| \le 
\omega(s,t)^{j/p}, \qquad
 j = 1,\ldots,[p], \quad (s,t) \in \triangle
\nn
\\
& &
| \Omega_t - \Omega_s|_{L( {\cal V},{\cal V})}  
\le  \omega(s,t)^{1/q},
\qquad (s,t) \in \triangle.
\nn
\end{eqnarray}
Then, 
$
\big|Y^j_{s,t} \big| 
\le 
C \omega(s,t)^{j/p},
\quad
j=1,\ldots,[p], \quad (s,t) \in \triangle,
$
where $C$ is a
 positive constant
depending only on $p$, $q$, and $\omega(0,1)$.

\noindent
{\rm (3).}~
$(X,\Omega) \mapsto Y=\Gamma(X,M_{\Omega})$ extends to a continuous map 
from $G\Omega_p ( {\cal V})  \times  
 C_{0,q} (L({\cal V},{\cal V}))$ to $G\Omega_p ({\cal V})$.

\end{cor}

\subsection{preliminary lemmas}
In this subsection we will prove several simple lemmas 
for later use.
Proofs are easy.
Let $p \ge 2$
and let ${\cal V}$ and ${\cal W}$ be real Banach spaces.

%
\begin{lem}\label{lem.lin.trans}
For $\alpha \in L({\cal V},{\cal W})$
and $X \in G\Omega_p({\cal V})$,
set 
$\bar{\alpha}(X)=(1,\bar{\alpha}(X)^1,\ldots, \bar{\alpha}(X)^{[p]})$
by
$\bar{\alpha}(X)^j_{s,t}= \alpha^{\otimes j}(X^j_{s,t}),~(j=1,\ldots,[p])$.

\noindent 
{\rm (1).}~
Then, $\bar{\alpha}(X) \in G\Omega_p({\cal W})$
and
$$|\bar{\alpha}(X)_{s,t}^j-\bar{\alpha}(\hat{X})^j_{s,t}|
\le |\alpha|_{L({\cal V},{\cal W})}^j 
|X_{s,t}^j-\hat{X}_{s,t}^j|$$
for any $X, \hat{X}  \in G\Omega_p({\cal V})$.
In particular, $\bar{\alpha}: G\Omega_p({\cal V}) \to G\Omega_p({\cal W})$
is Lipschitz continuous.

\noindent
{\rm (2).}~
If $\alpha, \beta \in L({\cal V},{\cal W})$,
then 
$$|\bar{\alpha}(X)_{s,t}^j-\bar{\beta}(X)^j_{s,t}|
\le j |\alpha-\beta|_{L({\cal V},{\cal W})}
(|\alpha|_{L({\cal V},{\cal W})} \vee |\beta|_{L({\cal V},{\cal W})} )^{j-1} 
|X_{s,t}^j|$$
for any $X \in G\Omega_p({\cal V})$.
\end{lem}

\Proof
By the basic property of the projective norm,
we can see that 
$|\alpha^{\otimes j}|_{L({\cal V}^{\otimes j},{\cal W}^{\otimes j} )}  
= |\alpha|_{L({\cal V},{\cal W})}^j$.
The rest is easy
\QED

%
%

The following is a slight modification of
Lemma 6.3.5 in p.171, \cite{lq}.
The proof is easy.
Note that the choice $\delta>0$ is independent of $\omega$.
In the following, 
$\Gamma_{a,b,c}: G\Omega_p ( {\cal V}  \oplus {\cal W}^{\oplus 2} )
\to  G\Omega_p ( {\cal V}  \oplus {\cal W}^{\oplus 2} )$
is defined by $\Gamma_{a,b,c}= 
\overline{ a {\rm Id}_{{\cal V}} 
\oplus b{\rm Id}_{{\cal W}} 
\oplus  c {\rm Id}_{{\cal W}} }$
($a,b,c \in {\mathbb R}$).
\begin{lem}\label{lem.delta.cont}
If $K, \hat{K} \in G\Omega_p ( {\cal V}  \oplus {\cal W}^{\oplus 2} )$ with 
$\pi_{{\cal V}} (K)=X,~\pi_{{\cal V}} (\hat{K})=\hat{X}$
and if $\omega$ is a control such that
\begin{eqnarray*}
|X^j_{s,t} |,  |\hat{X}^j_{s,t} | &\le& \frac12 \omega(s,t)^{j/p},
\quad
|K^j_{s,t} |,  |\hat{K}^j_{s,t} |  \le ( C_2 \omega(s,t) )^{j/p}
\\
|X^j_{s,t} -\hat{X}^j_{s,t} | &\le&  \frac{\ve}{2} \omega(s,t)^{j/p},
\quad
|K^j_{s,t} -\hat{K}^j_{s,t} |  \le \ve ( C_2 \omega(s,t) )^{j/p},
\end{eqnarray*}
for all $j=1,\ldots,[p]$ and $(s,t) \in \triangle$,
then there exists a constant $\delta \in (0,1]$ depending only
on $C_2$ and $[p]$ such that, for any $\delta_1, \delta_2 \in (0,\delta]$,
we have
\[
\bigl|
( \Gamma_{1, \delta_1, \delta_2}  K  )^j_{s,t}
\bigr|
\le 
\omega(s,t)^{j/p},
\quad
\bigl|
( \Gamma_{1, \delta_1, \delta_2}  K  )^j_{s,t}
-
( \Gamma_{1, \delta_1, \delta_2}  \hat{K}  )^j_{s,t}
\bigr|
\le 
\ve \omega(s,t)^{j/p},
\]
for all $j=1,\ldots,[p]$ and $(s,t) \in \triangle$.
\end{lem}

%

\begin{lem}\label{lem.prolong}
{\rm (1)}~
Let ${\cal V}$ be a real Banach space and
$\{0=T_0<T_1<\cdots<T_N=1\}$ be a partition of $[0,1]$.
For each $i=1,2,\ldots,N$, 
$A(i): \triangle_{[T_{i-1},T_i]} \to T^{([p])} ({\cal V})$ 
is a geometric rough path which 
satisfies that
\[
\bigl| 
A(i)^j_{s,t}
\bigr|
\le 
 \omega(s,t)^{j/p},
\qquad
\text{
$j=1,\ldots,[p]$,  $(s,t) \in \triangle_{[T_{i-1},T_i]}$.}
\]
We define $A: \triangle \to T^{([p])} ({\cal V})$ by
\[
A_{s,t}= A_{s,T_k}(k) \otimes A_{T_k,T_{k+1}} (k+1)
\otimes \cdots  \otimes A_{T_{l-1},t}(l)
\qquad
\text{in $T^{([p])} ({\cal V})$}
\]
for $(s,t) \in \triangle$
such that $s \in [T_{k-1}, T_k]$ and $t \in [T_{l-1}, T_l]$.
Then, $A$ is a geometric rough path such that
\[
\bigl| 
A^j_{s,t}
\bigr|
\le 
( C \omega(s,t))^{j/p},
\qquad
\text{
$j=1,\ldots,[p]$ and  $(s,t) \in \triangle$.}
\]
Here, $C>0$ is a constant which depends only on $p$ and $N$.
\\
\\
\noindent
{\rm (2)}~
Let $A(i)$ and $\hat{A}(i)$ be two such rough paths
on restricted intervals as above ($i=1,2,\ldots,N$).
In addition to the assumption of {\rm (1)} for both $A(i)$ and $\hat{A}(i)$,
we also assume that
\[
\bigl| 
A(i)^j_{s,t}-  \hat{A}(i)^j_{s,t}
\bigr|
\le 
 \ve \omega(s,t)^{j/p},
\qquad
\text{for 
$j=1,\ldots,[p]$ and  $(s,t) \in \triangle_{[T_{i-1},T_i]}$.}
\]
Then, $A$ and $\hat{A}$ defined as above satisfy that
\[
\bigl| 
A^j_{s,t}-  \hat{A}^j_{s,t}
\bigr|
\le \ve
( C' \omega(s,t))^{j/p},
\qquad
\text{for 
$j=1,\ldots,[p]$ and  $(s,t) \in \triangle$.}
\]
Here, $C'>0$ is a constant which depends only on $p$ and $N$.
\end{lem}

\Proof
We can show this by straight forward computation.
We prove (2) for example.
Assume $s \in [T_0, T_1]$ and $t \in [T_{N-1}, T_N]$, because the other cases are easier.
For the first level path,
\[
A^1_{s,t} = A(1)^1_{s,T_1}  + A(2)^1_{T_1,T_2}  +\cdots + A(N)^1_{T_{N-1},t}.
\]
Hence, 
\begin{eqnarray*}
|  A^1_{s,t}  - \hat{A}^1_{s,t}  |  &\le&    | A(1)^1_{s,T_1}  - \hat{A}(1)^1_{s,T_1}  |  +\cdots +  |  A(N)^1_{T_{N-1},t} -\hat{A}(N)^1_{T_{N-1},t}|
\\
&\le&
\ve \omega(s,T_1)^{1/p}  +\cdots +\ve \omega(T_{N-1},t)^{1/p} \le \ve N\omega(s,t)^{1/p}.
\end{eqnarray*}
For simplicity of notation, we write $s=\tau_0,  T_k = \tau_i ~(1\le i \le N-1), t=\tau_N$.
The second level path satisfies
\[
A^2_{s,t} = \sum_{i=1}^N  A(i)^2_{\tau_{i-1},  \tau_i} 
  +\sum_{ 1 \le   i<k \le N }  A(i)^1_{\tau_{i-1},  \tau_i}  \otimes  A(k)^1_{\tau_{k-1},  \tau_k}.
\]
From this, we see that 
\begin{eqnarray*}
|  A^2_{s,t}  - \hat{A}^2_{s,t}  |  
\le
 \ve (N  + 2 \frac{ N(N-1)}{2}  )\omega(s,t)^{2/p}.
\end{eqnarray*}
The higher level paths ($j \ge 3$) can be done in the same way.
\QED


The first assertion of the following lemma
is Corollary 3.2.1, \cite{lq}.
The second one is a weaker version 
of Theorem 3.2.2, \cite{lq}.
\begin{lem}\label{lem.arp}
Let ${\cal V}$ be a real Banach space and let 
$A, B: \triangle \to T^{([p])}({\cal V})$
be almost rough paths.
\\
\noindent
{\rm (1).}
If there exist $\theta >1$ and a control function
$\omega$ such that
\begin{eqnarray*}
|A^j_{s,t}| &\le& \omega(s,t)^{j/p},
\qquad
j=1,\ldots,[p], \quad  (s,t) \in \triangle,
\nn\\
|A^j_{s,t}- ( A_{s,u}  \otimes A_{u,t} )^j | 
&\le& \omega(s,t)^{\theta},
\qquad
j=1,\ldots,[p], \quad   (s,u), (u,t)\in \triangle,
\end{eqnarray*}
then,
there is a unique rough path $\hat{A}$
associated to $A \in \Omega_p ({\cal V})$ such that
\[
|\hat{A}^j_{s,t}| \le C\omega(s,t)^{j/p},
\qquad
j=1,\ldots,[p], \quad  (s,t) \in \triangle.
\]
Here, $C>0$ is a constant which 
depends only on $p, \theta, \omega(0,1)$.
\\
\noindent
{\rm (2).}
Assume that there exists $\theta >1$ and a control function
$\omega$ such that,
for all $j=1,\ldots,[p]$, $\ve >0$ and $(s,u), (u,t)\in \triangle$, 
\begin{eqnarray*}
|A^j_{s,t}|,  \,\,
|B^j_{s,t}| \le \omega(s,t)^{j/p},
\qquad
|A^j_{s,t}- B^j_{s,t}| &\le& \ve \omega(s,t)^{j/p},
\nn\\
|A^j_{s,t}- ( A_{s,u}  \otimes A_{u,t} )^j |,
\,\,\,
|B^j_{s,t}- ( B_{s,u}  \otimes B_{u,t} )^j |  &\le& \omega(s,t)^{\theta},
\nn\\
\Bigl| \bigl ( 
A^j_{s,t}- ( A_{s,u}  \otimes A_{u,t} )^j  \bigr)
-
\bigl ( 
B^j_{s,t}- ( B_{s,u}  \otimes B_{u,t} )^j   \bigr)
\Bigr|  &\le& \ve \omega(s,t)^{\theta}
\end{eqnarray*}
hold.
Then, the associated rough paths, $\hat{A}$ and $\hat{B}$, satisfy that
\[
|\hat{A}^j_{s,t}- \hat{B}^j_{s,t}| \le
\ve C \omega(s,t)^{j/p},
\qquad
j=1,\ldots,[p], \quad (s,t)\in \triangle.
\]
Here, $C>0$ is a constant which 
depends only on $p, \theta, \omega(0,1)$.
\end{lem}

\Proof
For the first assertion, see Corollary 3.2.1, \cite{lq}.
We can show the second assertion by modifying the proof of  Theorem 3.2.2, \cite{lq}.
In this proof, the positive constant $C$ may vary from line to line.
First let us consider the case $j=1$.
Let ${\cal P}= \{0=t_0 <t_1 <\cdots < t_r =t   \}$ be a partition of $[s,t]$.
Recall that
\[
\hat{A}^1_{s,t}  =\lim_{|{\cal P}| \to 0}  A({\cal P})^1_{s,t},   
\quad 
\mbox{where} \quad 
A({\cal P})^1_{s,t} =  \sum_{i=1}^r  A^1_{t_{i-1} ,t_{i}}.
\]
There exists a point  $t_l \in {\cal P} \setminus \{ s,t\}$ which satisfies $\omega(t_{l-1} ,t_{l+1} ) \le 2(r-1)^{-1} \omega(s,t)$.
Then,  denoting ${\cal P}' ={\cal P} \setminus \{ t_l \}$, we see that
\begin{eqnarray*}
\lefteqn{
\Bigl|  \{  A({\cal P})^1_{s,t} -  A({\cal P}')^1_{s,t}  \}
-
\{ B({\cal P})^1_{s,t} -  B({\cal P}' )^1_{s,t}  \}
\Bigr|
}
\\
&=&
\Bigl|
 ( A^1_{t_{l-1} ,t_{l}}  + A^1_{t_{l} ,t_{l+1}}  - A^1_{t_{l-1} ,t_{l+1}}  )
 -
  ( B^1_{t_{l-1} ,t_{l}}  + B^1_{t_{l} ,t_{l+1}}  - B^1_{t_{l-1} ,t_{l+1}}  )
   \Bigr|
   \\
   &\le&    \ve   \omega( t_{l-1} ,t_{l+1} )^{\theta}  \le   \frac{\ve 2^{\theta} }{(r-1)^{\theta} }  \omega(s,t)^{\theta}.
   \end{eqnarray*}
By routine argument, 
\begin{eqnarray}
\lefteqn{
\Bigl|  \{  \hat{A}^1_{s,t} -  A^1_{s,t}  \}
-
\{ \hat{B}^1_{s,t} -  B^1_{s,t}  \}
\Bigr|
}
\nn\\
&\le&
\lim_{|{\cal P}| \to 0}
\Bigl|  \{  A({\cal P})^1_{s,t} -  A^1_{s,t}  \}
-
\{ B({\cal P})^1_{s,t} -  B^1_{s,t}  \}
\Bigr|
\le   \ve 2^{\theta} \zeta(\theta)  \omega(s,t)^{\theta}.
\label{eq.endisnear}
\end{eqnarray}
Combining this with the estimate for $A^1_{s,t} - B^1_{s,t}$ and $\omega(s,t)^{\theta} \le C \omega(s,t)^{1/p} $, 
we prove the case for $j=1$.

Next consider the case $j=2$ ($[p] \ge 2$).
Recall  that 
\[
\hat{A}^2_{s,t}  =\lim_{|{\cal P}| \to 0}  A({\cal P})^2_{s,t},   
\quad 
\mbox{where} \quad 
A({\cal P})^2_{s,t} =  \sum_{i=1}^r  ( A^2_{t_{i-1} ,t_{i}} +   \hat{A}^1_{s,t_{i-1}}  \otimes \hat{A}^1_{t_{i-1} ,t_{i}} ).
\]

In a similar way, 
\begin{eqnarray*}
\lefteqn{
\Bigl|  \{  A({\cal P})^2_{s,t} -  A({\cal P}')^2_{s,t}  \}
-
\{ B({\cal P})^2_{s,t} -  B({\cal P}' )^2_{s,t}  \}
\Bigr|
}
\\
&=&
\Bigl|
 ( A^2_{t_{l-1} ,t_{l}}  + A^2_{t_{l} ,t_{l+1}}  - A^2_{t_{l-1} ,t_{l+1}}   + A^1_{t_{l-1}, t_l}  \otimes A^1_{t_{l} ,t_{l+1}}   )
 \\
 && \qquad \qquad
 -
  ( B^2_{t_{l-1} ,t_{l}}  + B^2_{t_{l} ,t_{l+1}}  - B^2_{t_{l-1} ,t_{l+1}}  + B^1_{t_{l-1}, t_l}  \otimes B^1_{t_{l} ,t_{l+1}}   )
   \Bigr|
  \\
  && 
  +  \Bigl|  (  \hat{A}^1_{t_{l-1}, t_l}  \otimes \hat{A}^1_{t_{l} ,t_{l+1}}  -A^1_{t_{l-1}, t_l}  \otimes A^1_{t_{l} ,t_{l+1}} )
  -   (  \hat{B}^1_{t_{l-1}, t_l}  \otimes \hat{B}^1_{t_{l} ,t_{l+1}}  -B^1_{t_{l-1}, t_l}  \otimes B^1_{t_{l} ,t_{l+1}} )
    \Bigr|.
   \end{eqnarray*}
The first term on the  right hand side is clearly dominated by $  \ve   \omega( t_{l-1} ,t_{l+1} )^{\theta}$ by the assumption.
From (\ref{eq.endisnear}) and the assumption, it is also easy to see that the second term is also dominated by 
$ \ve  C\omega( t_{l-1} ,t_{l+1} )^{\theta+(1/p)}$.
Hence, the right hand side is dominated by $\ve C \omega (t_{l-1} ,t_{l+1}   )^{\theta} \le \ve C (r-1)^{-\theta} \omega (s,t)^{\theta} $.
By using the same argument as above,  we can show that
\begin{eqnarray}
\lefteqn{
\Bigl|  \{  \hat{A}^2_{s,t} -  A^2_{s,t}  \}
-
\{ \hat{B}^2_{s,t} -  B^2_{s,t}  \}
\Bigr|
}
\nn\\
&\le&
\lim_{|{\cal P}| \to 0}
\Bigl|  \{  A({\cal P})^2_{s,t} -  A^2_{s,t}  \}
-
\{ B({\cal P})^2_{s,t} -  B^2_{s,t}  \}
\Bigr|
\le   \ve C  \omega(s,t)^{\theta}.
\label{eq.endisnear}
\end{eqnarray}
Thus, we can prove the case for $j=2$.
The case $j \ge 3$ can be shown similarly. 
\QED


%
\section
{A slight generalization of Lyons' continuity theorem}
The aim of this section is to slightly generalize
Lyons' continuity theorem 
(also known as ``the universal limit theorem'') 
for It\^{o} maps in the rough path theory.
This section is based on Sections 5.5 and 6.3 in Lyons and Qian \cite{lq}.
Notations and results in these sections will be referred frequently.
The following points seem new:
\begin{enumerate}
\item
We let the coefficient of an It\^o map also vary. 
\item
We give an explicit estimate for the ``local Lipschitz 
continuity'' of It\^o maps (Theorem 6.3.1, \cite{lq}).
\end{enumerate}

\subsection{A review of integration along a geometric rough path}

For a real Banach space ${\cal V}$, 
we denote by $G\Omega_p ({\cal V})$ the space of geometric 
rough paths over ${\cal V}$,
where $p \ge 2$ is the roughness
and the tensor norm is the projective norm.
Let ${\cal W}$ be another real Banach space
and let $f: {\cal V} \to L({\cal V},{\cal W})$ be $C^{[p]+1}$
in the sense of Fr\'echet differentiation.
We say $f \in C^{[p]+1}_{b,loc} ( {\cal V}, L({\cal V},{\cal W}))$
if $f$ is $C^{[p]+1}$ from $\cal V$ to $L({\cal V},{\cal W})$
such that $|D^j f|~(j=0,1,\ldots,[p]+1)$ are bounded
on any bounded set.
In this subsection we will consider $\int f(X)dX$ (for $X \in G\Omega_p (X)$).

For $n \in {\mathbb N}=\{1,2,\ldots\}$,
let $\Pi_n$ be the set of all permutations of $\{1,2,\ldots,n\}$.
We define the left action of $\pi \in \Pi_n$ on ${\cal V}^{\otimes n}$
by
$\pi(v_1 \otimes \cdots \otimes v_n)=
v_{\pi^{-1}(1)} \otimes \cdots \otimes v_{\pi^{-1}(n)}$.
(Note that this is different from the definition in \cite{lq},
where the right action is adopted.
However, this does not matter so much since no composition of
permutations will appear below.)
Given ${\bf l}=\{l_1,\ldots,l_i\}$
($l_1,\ldots,l_i \in {\mathbb N}$),
let $|{\bf l}|=l_1+\cdots+l_i$.
We say $\pi \in \Pi_{{\bf l}}$
if $\pi \in \Pi_{|{\bf l}|}$
satisfies the following conditions:
\begin{eqnarray}
\pi(1)< &\cdots&< \pi(l_1),
\nn\\
\pi(l_1+1)< &\cdots&< \pi(l_1+l_2),
\nn\\
 &\cdots&
\nn\\
\pi(l_1+\cdots+l_{i-1}+1)< &\cdots&< \pi(|{\bf l}|),
\nn\\
\pi(l_1)< \pi(l_1+l_2)<  &\cdots&< \pi(|{\bf l}|).
\label{ineq.pi}
\end{eqnarray}
(Note: The last condition in (\ref{ineq.pi})
is missing in p.138, \cite{lq}.)

Let $X \in G\Omega_p({\cal V})$ be a smooth rough path.
Then, for all $(s,t) \in \triangle$,
\begin{eqnarray}\label{eq.chikan}
\int_{s<u_1<\cdots<u_i<t}
dX_{s,u_1}^{l_1} \otimes \cdots
 \otimes  dX_{s,u_i}^{l_i}
=
\sum_{ \pi \in \Pi_{ {\bf l}} }
\pi  X^{|{\bf l}|}_{s,t}.
\end{eqnarray}
(See Lemma 5.5.1, \cite{lq}.
By Corollary \ref{cor.appro.cq} and the Young integration theory, 
(\ref{eq.chikan}) also holds for 
$X$ lying above an element of $C_{0,q}({\cal V})$.)

For $f \in C_{b,loc}^{[p]+1} ( {\cal V}, L({\cal V},{\cal W}))$ 
in Fr\'echet sense,
we denote $f^j (x)$ for $D^jf(x)$ for simplicity
($j=0,1,\ldots,[p]+1$). 
For $X \in G\Omega_p( {\cal V})$, 
we define $Y \in C_{0,p}(\triangle, T^{([p])}({\cal W}) )$ by
\begin{equation}\label{defY}
Y^i_{s,t}
=
\sum_{ {\bf l}=(l_1,\ldots,l_i), 1\le l_j, |{\bf l}| \le [p]   }
f^{l_1 -1}(X_s) \otimes \cdots \otimes f^{l_i -1}(X_s)
\Bigl\la
\sum_{\pi \in \Pi_{{\bf l}} }
\pi X^{|{\bf l}|}_{s,t}
\Bigr\ra
\end{equation}
for $(s,t) \in \triangle$ and $i=1,\ldots,[p]$.

For a smooth rough path $X$, it is easy to see
that, for any $(s,u),(u,t) \in \triangle$,
\begin{equation}\label{eq.sono1}
\sum_{l=1}^{[p]} f^{l -1}(X_s) \la dX^l_{s,t} \ra
=
\sum_{l=1}^{[p]}
\bigl(
 f^{l -1}(X_u) -R_l (X_s,X_u)
\bigr)
\la dX^l_{u,t} \ra
\end{equation}
holds.
Here, for $x,y \in {\cal V}$,
\begin{eqnarray}\label{eq.2no5}
R_l (x,y) &=&  f^{l-1}(y)
-\sum_{k=j-1}^{[p]-1}f^{k}(x)\la (y-x)^{\otimes (k-l+1)} \ra
\nn\\
&=&
\int_0^1 d\theta 
\frac{(1-\theta)^{[p]-l} }{ ([p]-l)!}
f^{[p]}(x+\theta (y-x)) \la (y-x)^{\otimes ([p]-l+1) } \ra
\end{eqnarray}
(See Lemma 5.5.2 in \cite{lq}.
A key fact is that the symmetric part of $X^l_{s,t}$
is $[X^1_{s,t}]^{\otimes l}/l!$.)

Using (\ref{eq.chikan}) and (\ref{eq.sono1}),
we see that, for a smooth rough path $X$
(and for a geometric rough path $X \in G\Omega_p ( {\cal V})$ by continuity),
\begin{equation}\label{semiyoung}
Y_{s,t}=Y_{s,u} \otimes M_{u,t},
\qquad
\text{ in $T^{([p])} ({\cal W} )$ for all $(s,u),(u,t) \in \triangle$}.
\end{equation}
Here, $M_{u,t}$ is given by $M^0_{u,t}=1$
and, for $i=1,\ldots,[p]$,
\begin{eqnarray}\label{eq.2no7}
M^i_{u,t}
&=&
\sum_{ {\bf l}=(l_1,\ldots,l_i), 1\le l_j,  |{\bf l}|\le [p]   }
\bigl(
 f^{l_1 -1}(X_u) -R_{l_1} (X_s,X_u)
\bigr)
\otimes 
\nn\\
&&
\qquad\qquad\qquad
\cdots \otimes 
\bigl(
 f^{l_i -1}(X_u) -R_{l_i} (X_s,X_u)
\bigr)
\Bigl\la
\sum_{\pi \in \Pi_{{\bf l}} }
\pi X^{|{\bf l}|}_{u,t}
\Bigr\ra.
\end{eqnarray}
(Actually, $M_{u,t}$ depends on $s$, too.
Eq. (\ref{semiyoung}) is Lemma 5.5.3 in \cite{lq}.)

Set $N^j_{t,u}=Y^j_{t,u}-M^j_{t,u}$.
Then, for a control function $\omega$ satisfying that
$|X^j_{s,t}| \le \omega(s,t)^{j/p}$
for $j=1,\ldots,[p]$,
it holds that
\begin{equation}\label{ineq.N}
|N^j_{u,t}| \le C M(f;[p],\omega(0,1))^j
 \omega(s,t)^{([p]+1)/p},
\qquad
j=1,\ldots,[p], \quad (s,t) \in \triangle.
\end{equation}
for some constant $C>0$ which depends only 
on $p, \omega(0,1)$.
Here, we set
\begin{equation}\label{def.M}
M(f;k,R) := \max_{0 \le j \le k}
\sup\{ |f^j(x)| \,:\, |x|\le R\}
\qquad
\text{ for $R>0$ and $k \in {\mathbb N}$}
\end{equation}
and $M(f;k) :=M(f;k,\infty)$.
Then, $Y$ defined by (\ref{defY}) is an almost rough path.
Indeed, noting that
$(Y_{s,u} \otimes N_{u,t})^k
=\sum_{i+j=k, i \ge 0, j\ge 1}  Y^i_{s,u} \otimes N^j_{u,t}$,
we can easily see from above that
there exists a constant $C>0$ which depends only 
on $p, \omega(0,1)$ such that
\begin{equation} \label{ineq.arp}
|Y^j_{s,t} -(Y_{s,u} \otimes Y_{u,t})^j| 
\le C 
M(f;[p],\omega(0,1))^j
\omega(s,t)^{([p]+1)/p}
\end{equation}
for all $j=1,\ldots,[p]$ and $(s,u),(u,t) \in \triangle$.

We denote by $\int f(X)dX$
the unique rough path which is associated to $Y$.
It is well-known that,
if $X$ is a smooth rough path lying above $t \mapsto X_t$, 
then $\int f(X)dX$ 
is a smooth rough path lying above $t \mapsto \int_0^t f(X_u)dX_u$.
By the next proposition, $X \mapsto \int f(X)dX$
is continuous, 
which implies that $\int f(X)dX \in G\Omega_p ( {\cal W})$.

The following is essentially Theorem 5.5.2 
in \cite{lq}.
Varying the coefficient $f$ 
and giving an explicit estimate for 
the local Lipschitz continuity are newly added.
We say $f_n \to f$ 
as $n \to \infty$ in $C_{b,loc}^{k} ({\cal V}, L({\cal V},{\cal W}))$
if $M(f-f_n;k,R) \to 0$ as $n \to \infty$
for any $R>0$.

\begin{prop}\label{thm.int.cont}
We assume 
$f, g \in  C_{b,loc}^{[p]+1} ({\cal V}, L({\cal V},{\cal W}))$.
%
%
\\
\noindent
{\rm (1)}.~
Let $f$ be as above.
If $X \in G\Omega_p({\cal V})$ and a control function $\omega$
satisfy that
\[
|X^j_{s,t}| \le \omega(s,t)^{j/p},
\qquad
j=1,\ldots,[p], \quad (s,t) \in \triangle,
\]
then, for a constant $C>0$ which depends only on 
$p$ and $\omega(0,1)$, it holds that
\[
| \int_s^t f(X)dX^j | \le 
C
 M(f;[p],\omega(0,1))^j  \omega(s,t)^{j/p},
\qquad
j=1,\ldots,[p], \quad (s,t) \in \triangle.
\]
%
%
%
\\
\noindent
{\rm (2)}.~
Let $f, g$ be as above and 
let $X, \hat{X} \in G\Omega_p({\cal V})$ such that, for a control 
function $\omega$,
\begin{eqnarray*}
|X^j_{s,t}|,  \,\,
|\hat{X}^j_{s,t}| &\le& \omega(s,t)^{j/p},
\nn\\
|X^j_{s,t}- \hat{X}^j_{s,t}| &\le& \ve \omega(s,t)^{j/p},
\qquad
j=1,\ldots,[p], \quad  (s,t) \in \triangle.
\end{eqnarray*}
Then, 
\begin{eqnarray*}
\bigl| \int_s^t f(X)dX^j -  \int_s^t g(\hat{X})d\hat{X}^j  \bigr| 
&\le&
C  M(f-g;[p],\omega(0,1))  \tilde{M}_{[p]}^{j-1}
  \omega(s,t)^{j/p}
\\
&&+
\ve  C
 \tilde{M}_{[p]+1}^j  \omega(s,t)^{j/p}
\end{eqnarray*}
for $j=1,\ldots,[p]$ and $(s,t) \in \triangle$.
Here, $\tilde{M}_k=M(f;k,\omega(0,1)) \vee  M(g;k,\omega(0,1))$
and $C>0$ is a constant which depends only on 
$p$ and $\omega(0,1)$
.
%
%
%
%
\\
\noindent
{\rm (3)}.~
In particular,
the following map is continuous:
\[
(f, X) \in  C_{b,loc}^{[p]+1} ({\cal V}, L({\cal V},{\cal W}))
 \times G\Omega_p({\cal V})
\mapsto
\int f(X)dX
\in G\Omega_p({\cal W}).
\]
%
%
\end{prop}

\Proof
In this proof the constant $C>0$ may vary from line to line.
To show the first assertion, 
note that
\[
|Y^j_{s,t}| \le C M(f;[p],\omega(0,1))^{j}
\omega(s,t)^{j/p}\qquad
j=1,\ldots,[p], \quad (s,t) \in \triangle.
\]
From this and (\ref{ineq.arp})
we may apply Lemma \ref{lem.arp} to 
$M(f;[p],\omega(0,1))^{-1} \cdot Y$
to obtain the first assertion.

Now we prove the second assertion.
First we consider $\int_s^t g(X)dX^j-\int_s^t g(\hat{X})d\hat{X}^j$.
It is easy to see from (\ref{defY}) that
\[
|Y^j_{s,t} -\hat{Y}^j_{s,t} |
\le
\ve C M(g;[p]+1,\omega(0,1))^j 
 \omega(s,t)^{j/p}, 
\qquad
j=1,\ldots,[p], \quad (s,t) \in \triangle.
\]
Similarly, we have from (\ref{ineq.N})
that
\[
|N^j_{u,t}  - \hat{N}^j_{u,t} | 
\le \ve C M(g;[p]+1,\omega(0,1))^j
 \omega(s,t)^{([p]+1)/p},
\qquad
j=1,\ldots,[p], \quad (s,t) \in \triangle.
\]
This implies that
\begin{eqnarray*}
\lefteqn{
\Bigl|
[ Y^j_{s,t} -(Y_{s,u} \otimes Y_{u,t})^j ]
-[ \hat{Y}^j_{s,t} -(\hat{Y}_{s,u} \otimes \hat{Y}_{u,t})^j ]
\Bigr| 
}
\nn\\
&\le& \ve C  M(g;[p]+1,\omega(0,1))^j\omega(s,t)^{([p]+1)/p}
\qquad
j=1,\ldots,[p], \quad (s,t) \in \triangle.
\nn
\end{eqnarray*}
Now, setting $A=M(g;[p]+1,\omega(0,1))^{-1}\cdot Y$
and $B=M(g;[p]+1,\omega(0,1))^{-1}\cdot \hat{Y}$,
we may use Lemma \ref{lem.arp} 
to obtain that, for all $j=1,\ldots,[p], \quad (s,t) \in \triangle$,
\[
\Bigl|
\int_s^t g(X)dX^j-\int_s^t g(\hat{X})d\hat{X}^j
\Bigr|
\le
\ve  C
 \tilde{M}_{[p]+1}^j  \omega(s,t)^{j/p}.
\]

Next we consider $\int_s^t f(X)dX^j-\int_s^t g(X)dX^j$.
It is easy to see from (\ref{defY}) that
\[
|Y(f)^j_{s,t} -Y(g)^j_{s,t} |
\le
 C M(f-g;[p],\omega(0,1)) \tilde{M}_{[p]}^{j-1} 
 \omega(s,t)^{j/p}, 
\quad
j=1,\ldots,[p], \quad (s,t) \in \triangle.
\]
Similarly, we have from (\ref{ineq.N})
that, for $j=1,\ldots,[p]$ and $(s,t) \in \triangle$, 
\[
|N(f)^j_{u,t}  - N(g)^j_{u,t} | 
\le  C M(f-g;[p],\omega(0,1))  \tilde{M}_{[p]}^{j-1}
 \omega(s,t)^{([p]+1)/p}.
\]
This implies that
\begin{eqnarray*}
\lefteqn{
\Bigl|
[ Y(f)^j_{s,t} -(Y(f)_{s,u} \otimes Y(f)_{u,t})^j ]
-[ Y(g)^j_{s,t} -( Y(g)_{s,u} \otimes Y(g)_{u,t})^j ]
\Bigr| 
}
\nn\\
&\le& C  M(f-g;[p],\omega(0,1))
\tilde{M}_{[p]}^{j-1}\omega(s,t)^{([p]+1)/p}
\qquad
j=1,\ldots,[p], \quad (s,t) \in \triangle.
\nn
\end{eqnarray*}
Now, setting $A=\tilde{M}_{[p]}^{-1}\cdot Y$,
$B=\tilde{M}_{[p]}^{-1}\cdot \hat{Y}$ and 
$\ve=M(f-g;[p],\omega(0,1))\tilde{M}_{[p]}^{-1}$,
we may use Lemma \ref{lem.arp} below
to obtain that, for all $j=1,\ldots,[p], \quad (s,t) \in \triangle$,
\[
\Bigl|
\int_s^t f(X)dX^j-\int_s^t g(X)dX^j
\Bigr|
\le  CM(f-g;[p],\omega(0,1))
 \tilde{M}_{[p]}^{j-1}  \omega(s,t)^{j/p}.
\]
This proves the second assertion.

The third assertion is trivial from the second.
\QED

%

\subsection{Existence of solutions of differential equations}

In this subsection we check existence and 
uniqueness of the differential equation
in the rough path sense.
Essentially, everything in this subsection
is taken from Section 6.3, Lyons and Qian \cite{lq}.

Let ${\cal V}, {\cal W}$ be real Banach spaces
and let $f \in C_{b}^{[p]+1} ({\cal W}, L( {\cal V},{\cal W}))$.
%
%
We define 
$F \in C_{b}^{[p]+1} ( {\cal V}\oplus {\cal W} , 
L( {\cal V}\oplus {\cal W},  {\cal V}\oplus {\cal W} ))$
by
\[
F(x,y) \la (\xi,\eta ) \ra
=
( \xi, f( y)\xi),
\qquad
(x,y), (\xi,\eta) \in {\cal V}\oplus {\cal W}.
\]
For given $f$ as above and $\beta >0$,
we set 
\begin{equation}\label{def.Psi}
\Psi_{\beta}(y,k)=
\beta \bigl(
f(y) -f(y  -\beta^{-1}k )
\bigr),
\qquad y, k \in {\cal W}.
\end{equation}
Clearly, $\Psi_{\beta}$ is a map from ${\cal W}\oplus {\cal W}$ 
to $L({\cal V},{\cal W})$.

In this section
we consider the following $\cal W$-valued differential equation for given $X$
in the rough path sense:
\[
dY_t = f(Y_t) dX_t,
\qquad
Y_0=0 \in {\cal W}.
\]
Note that, by replacing $f$ with $f(\,\cdot\,+y_0)$,
we can treat the same differential equation
with an arbitrary initial condition $Y_0=y_0 \in {\cal W}$.
By a solution of the above differential equation,
we mean a solution of the following 
integral equation:
\begin{eqnarray}\label{eq.ode1}
Z_{s,t}^j
=
\int_s^t F (Z)dZ^j,
\qquad
j=1,\ldots,[p], 
(s,t) \in \triangle,
\text{ and
$\pi_V (Z)=X$.}
\end{eqnarray}
Note that $Z \in G\Omega_p(V\oplus W)$
and we also say $\pi_W (Z)=Y$ is a solution for given $X$.
(Here, $\pi_{{\cal V}}$ and $\pi_{{\cal W}}$ are the projections 
from ${\cal V} \oplus {\cal W}$
onto ${\cal V}$ and ${\cal W}$, respectively.)

As usual we use the Picard iteration:
\[
Z(n+1)= \int F(Z(n)) dZ(n),
\qquad 
\text{ with $Z(0)=(X,0)$.}
\]
For a smooth rough path $X$, this is equivalent to
\begin{eqnarray}
dX &=& dX,
\nn\\
dY(n+1)  &=&  f(Y(n)) dX,  \qquad \text{$Y(n+1)_0=0$}.
\nn
\end{eqnarray}
We may include the difference 
$D(n)=Y(n)-Y(n-1)$ in the equations:
for $n\in {\mathbb N}$,
\begin{eqnarray}
dX &=& dX,
\nn\\
dY(n+1)  &=&  f(Y(n)) dX,  
\nn\\
dD(n+1)  &=& \Psi_1 (Y(n), D(n)) dX.
\label{eq.ode2}
\end{eqnarray}
By scaling by $\beta >0$, we see that 
(\ref {eq.ode2}) is equivalent to the following:
for $n \in {\mathbb N},$
\begin{eqnarray}
dX &=& dX,
\nn\\
dY(n+1)  &=& f(Y(n)) dX, 
\nn\\
d \beta D(n+1)  &=& 
\Psi_{\beta} (Y(n),   \beta D(n)) dX.
\label{eq.ode3}
\end{eqnarray}
Set $\Phi_{\beta}:  {\cal V}\oplus {\cal W}^{\oplus 2} \to 
L( {\cal V}\oplus {\cal W}^{\oplus 2},
{\cal V}\oplus {\cal W}^{\oplus 2} )$ by
\[
\Phi_{\beta}(x,y,z)\la (\xi,\eta,\zeta)  \ra
=
\bigl( \xi,  f(y) \xi, \Psi_{\beta} ( y,z)\xi
\bigr),
\qquad
(x,y,z),(\xi,\eta,\zeta) \in   {\cal V}\oplus {\cal W}^{\oplus 2}.
\]
Clearly, $\Phi_{\beta}$ is the coefficient for (\ref{eq.ode3}).


\begin{lem}\label{lem.est.Phi}
Let $\beta \ge 1$ and $f$ be as above and define
$\Phi_{\beta}$. 
Then, for any $n$,
there exists a positive constant $c_n$ 
independent of $\beta \ge 1$ and $R>0$ such that
\begin{eqnarray*}
M( \Phi_{\beta};n,R) &\le& c_n(1+R)M(f;n+1,R),
\qquad
n \in {\mathbb N}, R>0.
\\
M( \Phi_{\beta};0,R) &\le& 1+(1+R)M(f;1,R),
\qquad R>0,
\end{eqnarray*}
\end{lem}

\Proof
First note that if $|(x,y,z)|=|x|+|y|+|z| \le R$,
then $|y-\theta \beta^{-1}z| \le R$ for any $\theta \in [0,1]$.
Clearly, $|\Phi_{\beta} (y,z)|=1+|f(y)|+|\Psi_{\beta} (y,z)|$.
By the mean value theorem,
$
\Psi_{\beta} (y,z)=\int_0^1 d\theta Df(y-\theta \beta^{-1}z) \la z \ra.
$
Hence,  we have the second inequality.

Denoting by $\pi_2, \pi_3$ the projection from 
${\cal V}\oplus {\cal W}^{\oplus 2}$
onto the second and the third component respectively,
we have 
\[
D\Psi_{\beta} (y,z)
=
\beta \bigl( Df(y)-Df(y-  \beta^{-1}z)
\bigr) \circ \pi_{2}
+Df(y-  \beta^{-1}z) \circ \pi_{3}.
\]
We can deal with the first term on the right hand side
in the same way to prove the first inequality 
of the lemma for $n=1$.
By continuing straight forward computation like this,
we can prove the rest.
\QED

%

For $\ve, \delta, \beta \in {\mathbb R}$
and a smooth rough path $K$ lying above 
$(k,l,m) \in {\rm BV}({\cal V}\oplus {\cal W}^{\oplus 2})$,
we define $\Gamma_{\ve, \delta, \beta} K$
by a smooth rough path lying above $(\ve k, \delta l, \beta m)$.
Then,
$K \mapsto \Gamma_{\ve, \delta, \beta} K$
extends to a continuous map from 
$G\Omega_p ({\cal V} \oplus {\cal W}^{\oplus 2})$
to itself. 
(Note that $\Gamma_{\ve, \delta, \beta}
=\overline{ \ve {\rm Id}_{\cal V} 
\oplus \delta {\rm Id}_{\cal W}
\oplus \beta {\rm Id}_{\cal W}  }$.)
%
%
The following is essentially Lemma 6.3.4, \cite{lq}.
\begin{lem}\label{lem.relat.Phi}
For any $r, \rho, \beta \in {\mathbb R}\setminus \{0\}$
and a geometric rough path 
$K \in G\Omega_p ({\cal V}\oplus {\cal W}^{\oplus 2})$,
we have
\begin{eqnarray}
\Gamma_{\rho, \rho, \beta \rho} \int \Phi_{r} (K)  dK
&=&
\int \Phi_{\beta r} (\Gamma_{\rho,1, \beta }  K)  
d  \Gamma_{\rho, 1, \beta }K.
\nn
\end{eqnarray}
\end{lem}

Now we consider the following iteration procedure
for given $X \in  G\Omega_p ({\cal V})$:
\begin{eqnarray}
K(n+1) =\int \Phi_1(  K(n)) dK(n)
\qquad
\mbox{ for $n \in {\mathbb N}$}
\label{eq.iter1}
\end{eqnarray}
with $K(0)=(X,0,0)$ and $K(1)=(X,f(0)X,f(0)X)$.
Note that $K(0)$ and $K(1)$ are well-defined 
not only for a smooth rough path $X$,
but also for any geometric rough path $X$.
Since $ \Phi_1$ is the coefficient for (\ref{eq.ode2}),
this corresponds to (\ref{eq.ode2}) 
at least if $X$ is a smooth rough path.

We also set, for $n \in {\mathbb N}$,
\begin{eqnarray}
Z(n+1) =\int F (  Z(n)) dZ(n)
\qquad
\text{with $Z(0)=(X,0)$.}
\label{eq.iter2}
\end{eqnarray}
Then, we have $\pi_{{\cal V}} (K(n))=X$
and $\pi_{{\cal V}  \oplus {\cal W}}
(K(n))=Z(n)$ for all $n \in {\mathbb N}$.
These relations are trivial when $X$ is a smooth rough path 
and can be shown by continuity for general $X \in  G\Omega_p ( {\cal V})$.

Now we set $H(n)=\Gamma_{1,1,\beta^{n-1}} K(n)$
for $\beta \neq 0$ and $n \in {\mathbb N}$.
If $X$ is a smooth rough path, then $H(n)$ is lying above
$(X,Y(n), \beta^{n-1} D(n))$.
From Lemma \ref{lem.relat.Phi}, we easily see that
\begin{eqnarray}\label{eq.H(n)}
H(n+1)=
\Gamma_{1,1,\beta} \int \Phi_{\beta^{n-1}}  (H(n))  dH(n),
\qquad
\text{ for $\beta \neq 0$ and $n \in {\mathbb N}$}.
\end{eqnarray}

From Proposition \ref{thm.int.cont} and Lemma \ref{lem.est.Phi}  
we see the following: 
If $K \in G\Omega_p (  {\cal V}\oplus {\cal W}^{\oplus 2})$ 
satisfies that, for some control $\omega$ with $\omega(0,1) \le 1$,
\begin{equation}
|K^j_{s,t}| \le \omega(s,t)^{j/p}
\qquad
\text{
for  
$j=1,\ldots,[p]$ and $(s,t) \in \triangle$,}
\nn
\end{equation}
then, 
for any $\beta >1$,
\begin{equation}\label{eq.def.C_1}
\bigl|
\int_s^t  \Phi_{\beta} (K)dK^j 
\bigr|
\le (C_1 \omega(s,t))^{j/p}
\qquad
\text{
for  
$j=1,\ldots,[p]$ and $(s,t) \in \triangle$.}
\end{equation}
Here, the constant $C_1>0$ can be chosen so that it depends only on 
$p, M(f;[p]+1)$. (Note that $M(f;[p]+1,1) \le M(f;[p]+1)$.
See Lemma \ref{lem.est.Phi}.
Note that (i) $C_1$ is independent of $\beta > 1$,
(ii) we can take the same $C_1$ 
even if we replace $f$ with $f(\,\cdot\, +y_0)$ for any $y_0 \in {\cal W}$.)

\begin{prop}\label{prop.choice.omega}
Let $C_1$ as in (\ref{eq.def.C_1}) and for this $C_1$
define $\delta$ 
as in Lemma \ref{lem.delta.cont}.
Choose $\beta >1$ arbitrarily and set $\rho =\beta/\delta$.
Let $X \in G\Omega_p ({\cal V})$ such that 
$
| 
X^j_{s,t}
|
\le 
 \hat{\omega}(s,t)^{j/p}$
for 
$j=1,\ldots,[p]$ and  $(s,t) \in \triangle$
for some control function $\hat\omega$.
%
%
%
%
\\
\noindent 
{\rm (1).}~
Set 
\begin{equation}\label{def.ome}
\omega(s,t) = \bigl( 2+ \frac{\rho +2|f|_{\infty}}{\rho}\bigr)^p 
  \hat\omega(s,t).
\end{equation}
Then, for
all $j=1,\ldots,[p]$ and  $(s,t) \in \triangle$,
\begin{eqnarray}
|X^j_{s,t}| \le  \frac12 \omega(s,t)^{j/p},
\qquad
| \Gamma_{\rho,1,1} K(1)^j_{s,t} |
\le  \bigl(  \rho^p \omega(s,t)  \bigr)^{j/p}.
\label{ineq.rho.ome1}
\end{eqnarray}
%
%
%
%
%
%
\\
\noindent 
{\rm (2).}~
Take $T_1>0$ so that $\rho^p \omega(0,T_1) \le 1$.
Then, on the restricted time interval $[0,T_1]$,
we have the following estimate:
\begin{equation}
\bigl|\Gamma_{\rho,1,1}  H(n)^j_{s,t}  \bigr| 
\le \bigl(  \rho^p \omega(s,t)  \bigr)^{j/p},
\qquad
j=1,\ldots,[p], \quad (s,t) \in \triangle_{[0,T_1]}.
\label{ineq.rho.ome2}
\end{equation}
%
%
%
%
%
%
\\
\noindent 
{\rm (3).}~
It holds that
\begin{equation}\label{ineq.rho.ome3}
\bigl|  H(n)^j_{s,t}  \bigr| 
\le \bigl( \rho^p \omega(s,t)  \bigr)^{j/p},
\qquad
j=1,\ldots,[p], \quad (s,t) \in \triangle_{[0,T_1]}.
\end{equation}
%
%
\end{prop}


\Proof
We prove the first assertion for the first and the second level paths.
Proofs for higher level paths are essentially the same.
It is obvious that 
$|X^j_{s,t}| \le 2^{-j}\omega(s,t)^{j/p}$ for $j=1,2$.
Since $\rho {\rm Id}_V \oplus f(0) \oplus f(0): 
{\cal V} \to  {\cal V} \oplus  {\cal W}^{\oplus 2}$
is a bounded linear map,
it naturally extends to 
$\overline{\rho {\rm Id}_{{\cal V}} \oplus f(0) \oplus f(0)}:
G\Omega_p (V)
\to G\Omega_p ( {\cal V}\oplus {\cal W}^{\oplus 2})$
and 
$\Gamma_{\rho,1,1} K(1)
= \overline{\rho {\rm Id}_{{\cal V}} \oplus f(0) \oplus f(0)}  \, X $.
Therefore, 
\[
| \Gamma_{\rho,1,1} K(1)^1_{s,t} |
\le (\rho +2|f|_{\infty})    |X^1_{s,t} |
\le
 \bigl(  \rho^p \omega(s,t)  \bigr)^{1/p}
\]
and, in a similar way, 
\[
| \Gamma_{\rho,1,1} K(1)^2_{s,t} |
\le (\rho^2  +4\rho |f|_{\infty}+ 4 |f|^2_{\infty} )    |X^2_{s,t} |
\le
 \bigl(  \rho^p \omega(s,t)  \bigr)^{2/p}.
\]

Now we prove the second assertion by induction.
Assume the inequality is true for $n$.
Then, using (\ref{eq.def.C_1}) with a 
new control function $\rho^p \omega$,
\begin{eqnarray}
\bigl|
\int_s^t  \Phi_{\beta^{n-1}} ( \Gamma_{\rho,1,1}  H(n))
d\Gamma_{\rho,1,1}  H(n)^j 
\bigr|
\le 
(C_1 \rho^p \omega(s,t))^{j/p}
%
\nn
\end{eqnarray}
for 
$j=1,\ldots,[p]$ and  $(s,t) \in \triangle_{[0,T_1]}$.
By Lemma \ref{lem.relat.Phi}, 
\begin{eqnarray}
\bigl|
\Gamma_{\rho,\rho,\rho}
\int_s^t  \Phi_{\beta^{n-1}} (   H(n))
d H(n)^j 
\bigr|
\le 
(C_1 \rho^p \omega(s,t))^{j/p}
\label{ineq.kk1}
\end{eqnarray}
for 
$j=1,\ldots,[p]$ and  $(s,t) \in \triangle_{[0,T_1]}$.
Note that
\begin{eqnarray}
\bigl| \pi_{{\cal V}}  \Gamma_{\rho,\rho,\rho}
\int_s^t  \Phi_{\beta^{n-1}} (  H(n))
d H(n)^j 
\bigr|
=
|\rho^j X^j_{s,t}|
\le 
\frac12 ( \rho^p \omega(s,t))^{j/p}
%
\label{ineq.kk2}
\end{eqnarray}
for 
$j=1,\ldots,[p]$ and  $(s,t) \in \triangle_{[0,T_1]}$.
Remembering that $\delta$ (in Lemma \ref{lem.delta.cont})
is independent of the control function,
we may use Lemma \ref{lem.delta.cont}
for (\ref{ineq.kk1}) and (\ref{ineq.kk2})
to obtain 
\begin{eqnarray}
\bigl| \Gamma_{1,\delta\beta^{-1} ,\delta }
\Gamma_{\rho,\rho,\rho}
\int_s^t  \Phi_{\beta^{n-1}} (  H(n))
d H(n)^j 
\bigr|
\le 
( \rho^p \omega(s,t))^{j/p}
%
\nn
\end{eqnarray}
for 
$j=1,\ldots,[p]$ and  $(s,t) \in \triangle_{[0,T_1]}$.
Note that $ \Gamma_{1,\delta\beta^{-1} ,\delta }
\Gamma_{\rho,\rho,\rho}
=\Gamma_{\rho,1,1}\Gamma_{1,1,\beta}.$
Then, by (\ref{eq.H(n)}), we have
 \begin{eqnarray}
\bigl| 
\Gamma_{\rho,1,1} H(n+1)^j_{s,t}
\bigr|
\le 
( \rho^p \omega(s,t))^{j/p}
%
\nn
\end{eqnarray}
for 
$j=1,\ldots,[p]$ and  $(s,t) \in \triangle_{[0,T_1]}$.
Thus, the induction was completed.

The third assertion is easily verified from the second
and Lemma \ref{lem.lin.trans}, 
since $\rho^{-1} \le 1$.
\QED

%
Set $Z^{\prime}(n) =(Z(n),0)$ when $X$ is (hence, $Z(n)$ is ) 
a smooth rough path.
Clearly, this naturally extends to the case of
geometric rough paths
and we use the same notation for simplicity.
Remember that
$K(n)= \Gamma_{1,1,\beta^{-(n-1)}} H(n)$
and 
$Z^{\prime}(n)= \Gamma_{1,1,0} H(n)$.
%
%
%
%
Therefore, it is easy to see from Lemma \ref{lem.lin.trans} that
 \begin{eqnarray}
\bigl| 
 K(n)^j_{s,t} -  Z^{\prime}(n)^j_{s,t}
\bigr|
\le 
j \beta^{-(n-1)}
( \rho^p \omega(s,t))^{j/p}
%
%
\label{ineq.kk3}
\end{eqnarray}
for 
$j=1,\ldots,[p]$ and  $(s,t) \in \triangle_{[0,T_1]}$.

Let $\alpha: {\cal V} \oplus {\cal W}^{\oplus 2}
\to {\cal V} \oplus {\cal W}$
be a bounded linear map defined by $\alpha\la(x,y,d) \ra= (x,y-d)$.
Then, $|\alpha|_{L(   {\cal V} \oplus {\cal W}^{\oplus 2} ,    
{\cal V} \oplus {\cal W})
} = 1$.
It is obvious that
$\overline{\alpha} K(n)=Z(n-1)$
and $\overline{\alpha} Z^{\prime} (n)=Z(n)$.
Combined with
(\ref{ineq.kk3}) and Lemma \ref{lem.lin.trans}, 
these imply that 
\begin{eqnarray}
\bigl| 
 Z(n)^j_{s,t} -Z(n-1)^j_{s,t}
\bigr|
\le 
j \beta^{-(n-1)}
( \rho^p \omega(s,t))^{j/p}
\label{ineq.conv1}
\end{eqnarray}
for 
$j=1,\ldots,[p]$ and  $(s,t) \in \triangle_{[0,T_1]}$.

Since $\beta >1$, the inequality above implies that
there exists $Z \in  G\Omega_p ( {\cal V}\oplus {\cal W} )$
such that
$\lim_{n \to \infty } Z(n)=Z$ in $G\Omega_p ({\cal V}\oplus {\cal W})$.
In particular,
\begin{eqnarray}
\bigl| 
 Z(n)^j_{s,t} -Z^j_{s,t}
\bigr|
\le 
j \sum_{m=n}^{\infty } \beta^{-(m-1)}
( \rho^p \omega(s,t))^{j/p}
\label{ineq.conv2}
\end{eqnarray}
for 
$j=1,\ldots,[p]$ and  $(s,t) \in \triangle_{[0,T_1]}$.
This $Z$ is the desired solution of (\ref{eq.ode1})
 on the restricted time interval $[0,T_1]$.

\begin{remark}\label{rem_def_const}
Let us recall how the constants are defined.
First, $C_1$ depends only on $p$ and $M(f;[p]+1)$.
$\delta$ depends only on $C_1$ and $p$ and 
so does $\rho:=\beta/\delta$, where $\beta >1$ is arbitrary chosen.
Therefore, the constants on the right hand side of (\ref{ineq.conv1})
and (\ref{ineq.conv2})
depends only on $p$ and $M(f;[p]+1)$ (and the choice of $\beta >1$).
In particular, 
the constants $C_1$, $\delta$, $C_3$ (below) etc. can be chosen
independent of $y_0$
even if we replace $f$ with $f(\,\cdot\, +y_0)$.
\end{remark}


From  (\ref{ineq.conv2}) and Remark \ref{rem_def_const},
\begin{eqnarray}
\bigl| 
Z^j_{s,t}
\bigr|
\le 
( C_3  \omega(s,t))^{j/p},
\qquad
j=1,\ldots,[p], \quad (s,t) \in \triangle_{[0,T_1]}.
\label{ineq.conv3}
\end{eqnarray}
for some constant $C_3>0$ which depends only on 
$p$ and $M(f;[p]+1)$ (and the choice of $\beta> 1$).

Now we will consider prolongation of solutions.
Take $0=T_0<T_1<\cdots<T_N=1$ so that
$\rho^p \omega(T_{i-1},T_i)=1$ for $i=1,\ldots,N-1$
and
$\rho^p \omega(T_{N-1},T_N) \le 1$.
By the superadditivity of $\omega$, 
$N-1 \le \rho^p \omega(0,1)= (3\rho +2|f|_{\infty})^p \hat{\omega}(0,1)$.
Hence, $N$ is dominated by a constant which depends only on 
$\hat{\omega}(0,1)$,  $p$, and $M(f;[p]+1)$
(and the choice of $\beta >1$).

On $[T_{i-1}, T_{i}]$, we solve the differential equation 
(\ref{eq.ode1}) for a initial condition $Y_{T{i-1}}$
instead of $Y_0=y_0$.
By Remark \ref{rem_def_const} and
Remark \ref{rem_initial},
we see that (\ref{ineq.conv3}) holds on each time interval 
$[T_{i-1}, T_{i}]$ with the same $C_3>0$.
Then, we prolong them by using Lemma \ref{lem.prolong}.
Thus, we obtain a solution on the whole interval $[0,1]$.


\begin{remark}\label{rem_initial}
By the definition of $\omega$ in (\ref{def.ome}), 
we can take the same $T_i$'s
even if we replace $f$ by $f(\,\cdot\, +y_0)$.
This is the reason why we assume the boundedness of $|f|$.
This fact enables us to use a prolongation method as above.
If $|f|$ is of linear growth, 
then the prolongation of solution may fail.
The author does not know whether Lyons' continuity theorem
still holds or not in such a case.
\end{remark}
Summing up the above arguments, 
we have the following existence theorem.
\begin{thm}\label{thm.exist.sol}
Consider the differential equation (\ref{eq.ode1}).
Then, for any $X \in G\Omega_p({\cal V})$ 
there exists a unique solution $Z  \in G\Omega_p({\cal V} \oplus {\cal W})$
of  (\ref{eq.ode1}).
Moreover, if $X$ satisfies that
\begin{eqnarray}
\bigl| 
X^j_{s,t}
\bigr|
\le 
 \hat{\omega}(s,t)^{j/p},
\qquad
j=1,\ldots,[p],  \quad  (s,t) \in \triangle
\nn
\end{eqnarray}
for some control function $\hat\omega$,
then $Z$ satisfies that
\begin{eqnarray}
\bigl| 
Z^j_{s,t}
\bigr|
\le 
 (L \hat{\omega} (s,t))^{j/p},
\qquad
j=1,\ldots,[p], \quad  (s,t) \in \triangle.
\nn
\end{eqnarray}
Here, $L>0$ is a constant which depends only on $\hat{\omega}(0,1)$,
$p$, and $M(f;[p]+1)$.
\end{thm}

\Proof
All but uniqueness have already been shown.
Since we are mainly interested in estimates of solutions,
we omit a proof of uniqueness.
See pp. 177--178 in \cite{lq}.
\QED

\subsection{Local Lipschitz continuity of It\^o maps}
In this section we will prove the local Lipschitz continuity 
of It\^o maps.
In \cite{lq} the coefficient of It\^o maps is fixed. 
Here, we will let the coefficient vary.
This kind of generalization of Lyons' continuity theorem for the case $[p]=2$
was done by Coutin, Friz, and Victoir \cite{cfv}.

Let $X, \hat{X} \in G\Omega_p({\cal V})$ and $\hat{\omega}$
be a control function such that
\begin{eqnarray}
|X^j_{s,t}|, |\hat{X}^j_{s,t}|
\le
 \hat{\omega}(s,t)^{j/p},
\quad
|X^j_{s,t}-\hat{X}^j_{s,t}| &\le&
 \ve \hat{\omega}(s,t)^{j/p}
\label{ineq.XhatX}
\end{eqnarray}
for $j=1,\ldots,[p]$ and  $(s,t) \in \triangle$.
Let $y_0, \hat{y}_0 \in {\cal W}$ be initial points such that
\begin{eqnarray}
|y_0|, |\hat{y}_0| \le r_0, 
\qquad |y_0-\hat{y}_0| \le \ve'
\label{ineq.yhaty2}
\end{eqnarray}
Let $f, \hat{f} \in C_b^{[p]+2} ( {\cal W}  ,L({\cal V},{\cal W}))$.
For this $f$ and $\hat{f}$, 
we assume that, for any $R>0$,
\begin{eqnarray}
M(f;[p]+1), M(\hat{f};[p]+1) \le M,
\qquad
M(f-\hat{f};[p], R)  \le \ve_R''.
\label{ineq.fhatf2}
\end{eqnarray}
Essentially, $\ve'$ and  $\ve_R''$ vary only on 
$0 \le \ve' \le 2r_0$ and $0 \le \ve_R'' \le 2M$, respectively.


As in Theorem \ref{thm.exist.sol} and its proof,
we can solve the differential equation $f$ and $\hat{f}$
with the initial point $y_0$
and $\hat{y}_0$, respectively,
as in the previous subsection.
Set 
\begin{equation}\label{def.R0}
R_0:=1+r_0+ L \hat\omega(0,1),
\end{equation}
where $L>0$ is the constant in Theorem \ref{thm.exist.sol}.
Then, 
the first level paths of the solutions  satisfy
that, for any $n=1,2,\ldots,  t \in [0,1]$,
\begin{eqnarray}
|y_0+Y(n)^1_{0,t}|, |y_0+\hat{Y}(n)^1_{0,t}|
|y_0+Y^1_{0,t}|, |\hat{y}_0+ \hat{Y}^1_{0,t}|
\le
R_0-1.
\label{ineq.hankei}
\end{eqnarray}
So, under (\ref{ineq.XhatX})--(\ref{ineq.fhatf2}), 
we use information of $f$ and $\hat{f}$
only on $\{y \in {\cal W} ~|~ |y| \le R_0  \}$. 

The following is 
on the  local Lipschitz continuity of It\^o maps
and is the main theorem in this section.
\begin{thm}\label{thm.loc.Lip}
\begin{sloppypar}
Let $\hat{\omega}$, $r_0$, $R_0$,$M$, $\ve, \ve', \ve''_{R}$,
$X, \hat{X} \in G\Omega_p ({\cal V})$,
$y_0, \hat{y}_0 \in {\cal W}$, and $f, \hat{f} \in 
C_b^{[p]+2} ({\cal W},L({\cal V},  {\cal W}))$
satisfy (\ref{ineq.XhatX})--(\ref{ineq.fhatf2}).
Set $R_0$ and $\ve''=\ve''_{R_0}$ as in (\ref{def.R0}).
We denote by $Z, \hat{Z}$
be the solutions of It\^o maps corresponding to $f, \hat{f}$
with initial condition $y_0, \hat{y}_0$, respectively.
Then, in addition to Theorem \ref{thm.exist.sol},
we have the following;
there is a positive constant $L'$
such that
\end{sloppypar}
\begin{eqnarray}\label{eq.shimesu}
| Z^j_{s,t}- \hat{Z}^j_{s,t}|
\le
(\ve+\ve'+\ve'') (L' \hat{\omega} (s,t))^{j/p}
\end{eqnarray}
for $j=1,\ldots,[p]$ and  $(s,t) \in \triangle$.
Here, $L'$ depends only on $\hat{\omega}(0,1)$,
$p$, $r_0$ and $M$.
\end{thm}

For $k \in {\mathbb N}$ and $M>0$,
set ${\cal C}_M^{k}({\cal W}, L({\cal V},{\cal W}))
=\{f \in C_b^{k}(  {\cal W}, L({\cal V},{\cal W}) ) ~|~  M(f;k) \le M\}$.
We say $f_n \to f$ in ${\cal C}_M^{k}(  {\cal W}, L({\cal V},{\cal W}) )$ 
as $n \to \infty$
if $M(f-f_n;k, R) \to 0$ as $n \to \infty$ for any $R>0$.
\begin{cor}\label{cor.loc.Lip}
Let $Z$ be the solution corresponding to 
$X, f, y_0$ given as above.
Then, the map 
\[
(f, X, y_0) \in {\cal C}_M^{[p]+2}( {\cal W}, L({\cal V},{\cal W}))
 \times G\Omega_p ({\cal V})
\times {\cal W}
\mapsto
Z \in G\Omega_p ({\cal V} \oplus {\cal W})
\]
is continuous for any $M>0$. 
\end{cor}


\begin{demo}{{\it Proof.\ }}
Instead of (\ref{ineq.yhaty2}) and (\ref{ineq.fhatf2}), 
we assume that
\begin{eqnarray}
|y_0|, |\hat{y}_0| \le R_0 -1, 
\qquad |y_0-\hat{y}_0| \le \ve'
\label{ineq.yhaty}
\end{eqnarray}
and that
\begin{eqnarray}
M(f;[p]+1), M(\hat{f};[p]+1) \le M,
\qquad
M(f-\hat{f};[p], R_0)  \le \ve'' (:=\ve_{R_0}'').
\label{ineq.fhatf}
\end{eqnarray}
Clearly, if (\ref{ineq.XhatX})--(\ref{ineq.fhatf2})
are satisfied, 
then so are (\ref{ineq.XhatX}), (\ref{ineq.yhaty}),
and (\ref{ineq.fhatf}).
Now we will start with (\ref{ineq.XhatX}), (\ref{ineq.yhaty}),
and (\ref{ineq.fhatf}).


First, let us consider the case $f=\hat{f}$, $y_0=\hat{y}_0$,
$X \neq \hat{X}$.
(The argument of the previous subsection should be 
modified by substitution $f=f(\cdot +y_0)$, etc.) 
We will prove by induction that
there exists $T_1 \in (0,1]$ such that, 
\begin{eqnarray}\label{eq.new1}
|Z(n)^j_{s,t} -  \hat{Z}(n)^j_{s,t}|
\le 2\ve (1+M)^{[p]} \hat\omega(s,t)^{i/p},
\quad
n \in {\bf N},
j=1,\ldots,[p], (s,t)\in \triangle_{[0,T_1]}.
\end{eqnarray}
When $n=0$, (\ref{eq.new1}) clearly holds
since $Z(0)=(X,0)$ and $\hat{Z}(0)=(\hat{X},0)$
(with $T_1$ chosen arbitrarily).
Now we assume that (\ref{eq.new1}) is true for $n-1$.

Let $A(n)$ and $\hat{A}(n)$ be the almost rough path
that approximate $Z(n)$ and $\hat{Z}(n)$ as in (\ref{eq.iter2}), respectively.
By straight forward computation, we see that
\begin{eqnarray}\label{eq.new2}
| A(n)^j_{s,t} - \hat{A}(n)^j_{s,t}| 
\le 
2\ve(1+M)^{[p]} \bigl[
\frac12 +  \hat\omega(0,T_1)^{1/p} q_j (M, \hat\omega(0,T_1)^{1/p})
\bigr]
\hat\omega(s,t)^{j/p}
\end{eqnarray}
for any 
$n \in {\bf N},
j=1,\ldots,[p], (s,t)\in \triangle_{[0,T_1]}.$
Here, $q_j$ is a polynomial of two variables
with positive coefficients which is independent of $n$.
(In the sequel, $q_j$ may vary from line to line.)
Similarly,
\begin{eqnarray}\label{eq.new3}
\bigl|
( A(n)^j_{s,t} -   (A(n)_{s,u} \otimes A(n)_{u,t}^j )
-   ( \hat{A}(n)^j_{s,t} -   (\hat{A}(n)_{s,u} \otimes \hat{A}(n)_{u,t}^j )
\bigr| 
\nn\\
\le
2\ve(1+M)^{[p]}  q_j (M, \hat\omega(0,T_1)^{1/p})
\hat\omega(s,t)^{([p]+1)/p}
\end{eqnarray}
for any 
$n \in {\bf N},
j=1,\ldots,[p], (s,u), (u,t)\in \triangle_{[0,T_1]}.$
Using the same argument as in Lemma \ref{lem.arp} or 
Proposition \ref{thm.int.cont}, we have
\begin{eqnarray}\label{eq.new2}
| Z(n)^j_{s,t} - \hat{Z}(n)^j_{s,t}| 
\le 
2\ve(1+M)^{[p]} \bigl[
\frac12 +  \hat\omega(0,T_1)^{1/p} q_j (M, \hat\omega(0,T_1)^{1/p})
\bigr]
\hat\omega(s,t)^{j/p}
\end{eqnarray}
for any 
$n \in {\bf N},
j=1,\ldots,[p], (s,t)\in \triangle_{[0,T_1]}.$
Therefore, we can choose $T_1$  sufficiently small 
so that (\ref{eq.new1}) holds.
Note that the choice of $T_1$ independent of $n$.

The case
$f \neq \hat{f}$, $y_0=\hat{y}_0$,
$X =\hat{X}$ can be done in the same way.
Note that the differece of the $[p]+1$th derivative 
$|f^{[p]+1}(y)-\hat{f}^{[p]+1}(y)|$
is not involved in the argument.

The case
$f= \hat{f}$, $y_0 \neq \hat{y}_0$,
$X =\hat{X}$ can be reduced to the previous one
by setting $g=f(\cdot +y_0)$ and $\hat{g}=f(\cdot +\hat{y}_0)$.
Note that 
$|f^{j}(y +y_0)-f^{j}(y +\hat{y}_0)|$
for $1\le j \le [p]$
is dominated in terms of $|y_0- \hat{y}_0|$ and  $|f^{j+1}(y)|$.
Thus, we have shown the theorem on the restricted 
time interval $[0,T_1]$.


Now we consider the prolongation of solutions.
We have obtained that
$|Y^1_{0,T_1}- \hat{Y}^1_{0,T_1}| \le c_1 \ve_1$,
where $\ve_1 =\ve+\ve'+\ve''$ and $c_1$ is a positive constant.
Therefore, the difference of the initial values 
on the second interval $[T_1,T_2]$ is
dominated by 
$
| (y_0+ Y^1_{0,T_1})- (\hat{y}_0 +\hat{Y}^1_{0,T_1} )| 
\le 
\ve' + \ve_1 \le (1+c_1)\ve_1.$

Therefore, from the above computation and (\ref{ineq.hankei}),
on the second time interval $[T_1,T_2]$,
(\ref{ineq.XhatX}), (\ref{ineq.yhaty}),
and (\ref{ineq.fhatf})
are again satisfied,
with $\ve'$ in (\ref{ineq.yhaty}) being replaced with $(1+c_1)\ve_1$.
(Note that $\hat\omega$, $R_0$, and $M$ are not changed.)

Thus, we can do the same argument on $[T_1,T_2]$
with $\ve'$ being replaced  with $(1+c_1) \ve_1$
to obtain (\ref{eq.shimesu}) on the second interval.
Similarly, we obtain that
$
| (y_0+ Y^1_{0,T_2})- (\hat{y}_0 +\hat{Y}^1_{0,T_2} )| 
\le c_2 \ve_1.$
Repeating this argument finitely many times 
and use Lemma \ref{lem.prolong}, we can prove the theorem. 
\toy\end{demo}

\subsection{Estimate of difference of higher level paths of
two solutions}
Let $p \ge 2$ and $1/p+1/q >1$
and let $f$ be as in the previous subsection.
In this subsection we only consider the case $f=\hat{f}$.
For given $X$, the solution $Z=(X,Y)$ of (\ref{eq.ode1}) 
is denoted by $Z_X=(X,Y_X)$.
In Theorem \ref{thm.loc.Lip} we estimated the ``difference''
of $Z_X^j$ and $Z_{\hat{X}}^j$.
If $\hat{X} \in C_{0,q}({\cal V}) \subset G\Omega_p({\cal V})$, 
Then, $Z_{\hat{X}}$ is an element of $C_{0,q}({\cal W})$
and, therefore, the ``difference'' 
$Z_X-Z_{\hat{X}}$ is a ${\cal W}$-valued geometric rough path. 
The purpose of this section is to give an estimate
for $(Z_X-Z_{\hat{X}})^j$ in such a case.
Note that $(Z_X-Z_{\hat{X}} )^j$ and $Z_X^j-Z_{\hat{X}}^j$
is not the same if $j \neq 1$.

Roughly speaking, 
we will show that
$\| (Z_{X+ \Lm}-Z_{\Lm} )^j \|_{p/j} \le C(c_1,\kappa_0,f) \xi(X)^j$
for $\Lm \in C_{0,q}( {\cal V} )$ and $X \in G\Omega_p({\cal V})$
with $\|\Lm\|_q \le c_1$ and $\xi(X) \le \kappa_0$.
Since this is a continuous function of 
$(X, \Lm)$, we may only think of $X$ lying above 
an element of $C_{0,q}({\cal V}).$
Note also that
if we set 
\[
\omega (s,t):= \|\Lm\|_{q,[s,t]}^q +
\sum_{j=1}^{[p]} \kappa^{-p} \|X^j\|_{p/j,[s,t]}^{p/j},
\qquad 
( \text{here, we set $\kappa := \xi(X)$}),
\]
then this control function satisfies that
$\omega (0,1) \le c_1^q+[p]$,
$|\Lm^1_{s,t}| \le \omega (s,t)^{1/q}$,
and 
 $|X^j_{s,t}| \le \kappa^{j}\omega (s,t)^{j/p}$.

%
First we prove the following lemma.
Heuristically, $\kappa >0$ is a small constant. 
\begin{lem}\label{lem.geki}
Let $\kappa_0 >0$. 
Assume that a control function $\hat\omega$,
 $\Lm \in C_{0,q}({\cal V})$ and $X \in G\Omega_p({\cal V})$,
and $\kappa \in [0, \kappa_0]$
satisfy that
\[
|\Lm^1_{s,t}| \le \hat\omega (s,t)^{1/q},
\qquad
|X^j_{s,t}| \le (\kappa \hat\omega (s,t))^{j/p}
\qquad
j=1,\ldots,[p], (s,t) \in \triangle.
\]
Then, there is a positive constant $C$ which depends only on 
$\kappa_0, \hat\omega(0,1),p, M(f;[p]+1)$
such that
\[
\Bigl|
\bigl(
\kappa^{-1}X, \Lm, Y_{X+\Lm}
\bigr)^j_{s,t}
\Bigr|
\le C \hat\omega(s,t)^{p/j},
\qquad
s<t.
\]
Here, the left hand side denotes the $j$th level path 
of a ${\cal V}^{\oplus 2} \oplus {\cal W}$-valued geometric rough path
and $Y_{X+\Lm}$ denotes (the ${\cal W}$-component of) the solution of 
(\ref{eq.ode1}) for $X+\Lm$.
\end{lem}

\Proof 
We will proceed in a similar way as
in the previous subsection.
 We define 
$\tilde{F} \in C_{b}^{[p]+1} (  {\cal V}^{\oplus 2}\oplus {\cal W} ,
 L({\cal V}^{\oplus 2}\oplus {\cal W}  , 
{\cal V}^{\oplus 2}\oplus {\cal W} ))$
by
\[
\tilde{F}(x,x', y) \la (\xi,\xi',\eta ) \ra
=
( \xi, \xi',f( y)(\xi+\xi') ),
\quad
\text{ for $(x,x',y), (\xi,\xi',\eta) 
\in {\cal V}^{\oplus 2}\oplus {\cal W}$}
\]
and define 
$\tilde{\Phi}_{\beta}:   {\cal V}^{\oplus 2}\oplus {\cal W}^{\oplus 2} \to 
L(  {\cal V}^{\oplus 2}\oplus {\cal W}^{\oplus 2},
{\cal V}^{\oplus 2}\oplus {\cal W}^{\oplus 2} )$ by
\[
\tilde{\Phi}_{\beta}(x,x',y,z)\la (\xi,\xi',\eta,\zeta)  \ra
=
\bigl( \xi, \xi', f(y) (\xi+\xi'), \Psi_{\beta} ( y,z)(\xi+\xi')
\bigr),
\]
for $\beta >0$ and
$(x,x', y,z),(\xi,\xi',\eta,\zeta) 
\in {\cal V}^{\oplus 2}\oplus {\cal W}^{\oplus 2}.$
Note that $\tilde{\Phi}_{\beta}$ satisfies a similar 
estimates as in Lemma \ref{lem.est.Phi}.

Instead of (\ref{eq.ode1}), we now consider 
\begin{eqnarray}\label{eq.ode6}
\tilde{Z}_{s,t}^j
=
\int_s^t \tilde{F} ( \tilde{Z})d\tilde{Z}^j,
\qquad
j=1,2,\ldots,[p], 
(s,t) \in \triangle,
\text{ and
$\pi_{  {\cal V}^{\oplus 2}} (\tilde{Z})=(X,\Lm)$.}
\end{eqnarray}
It is easy to check that the solution of this equation
is $(X, \Lm, Y_{X+\Lm})$.

In order to solve  (\ref{eq.ode6}),
we use the iteration method as in (\ref{eq.ode2}) or (\ref{eq.ode3}).
More explicitly,
\begin{eqnarray}
dX &=& dX,
\qquad
d\Lm = d\Lm
\nn\\
d\tilde{Y}(n+1)  &=& f(\tilde{Y}(n)) d(X+\Lm), 
\nn\\
d \beta \tilde{D} (n+1)  &=& 
\Psi_{\beta} (\tilde{Y}(n),   \beta \tilde{D}(n)) d(X+\Lm).
\label{eq.ode7}
\end{eqnarray}

Now consider the iteration procedure 
for given $(X,\Lm) \in  G\Omega_p ({\cal V}) \times C_{0,q}({\cal V})$ and 
for $\tilde{F}$
and $\tilde{\Phi}_1$ as in (\ref{eq.iter1})
and (\ref{eq.iter2}).
Also define $\tilde{K}(n)$ and $\tilde{Z}(n)$
as in (\ref{eq.iter1})
and (\ref{eq.iter2}) with 
$ \tilde{K}(0)=(X,\Lm,0,0)$ and $\tilde{K}(1)=(X,\Lm,f(0)(X+\Lm),f(0)(X+\Lm))$.

Set $\tilde{H}(n)=\Gamma_{1,1,1,\beta^{n-1}}\tilde{K}(n)
=(X,\Lm, \tilde{Y}(n),\beta^{n-1}\tilde{D}(n))$.
Then, in the same way as in (\ref{eq.H(n)}), we have
\begin{eqnarray}\label{eq.H(n)2}
\tilde{H} (n+1)=
\Gamma_{1,1,1,\beta} \int \tilde{\Phi}_{\beta^{n-1}}  
(\tilde{H}(n))  d\tilde{H}(n),
\qquad
\text{ for $\beta \neq 0$ and $n \in {\mathbb N}$}.
\end{eqnarray}

Slightly modifying 
Proposition \ref{thm.int.cont} and Lemma \ref{lem.est.Phi}, 
we see from the estimates for $\tilde{\Phi}_{\beta}$
the following: 
If $\tilde{K} \in G\Omega_p ({\cal V}^{\oplus 2} 
 \oplus {\cal W}^{\oplus 2})$ 
satisfies that, for some control $\omega_0$ with $\omega_0(0,1) \le 1$,
\begin{equation}
|\Gamma_{1/\kappa,1,1,1} \tilde{K}^j_{s,t}| \le \omega_0(s,t)^{j/p}
\qquad
\text{
for  
$j=1,\ldots,[p]$,  \quad $(s,t) \in \triangle$}
\nn
\end{equation}
then, 
for any $\beta >1$,
\begin{equation}\label{eq.def.C_1_2}
\bigl|
\Gamma_{1/\kappa,1,1,1}
\int_s^t  \tilde{\Phi}_{\beta} (  \tilde{K})d\tilde{K}^j 
\bigr|
\le ( \tilde{C}_1 \omega_0 (s,t))^{j/p}
\qquad
\text{
for  
$j=1,\ldots,[p]$, $(s,t) \in \triangle$}
\end{equation}
Here, $\tilde{C}_1>0$ is a constant which depends only on 
$\kappa_0,p, M(f;[p]+1)$. (See Lemma \ref{lem.est.Phi}.
Note that (i) $\tilde{C}_1$ is independent of $\beta > 1$,
(ii) we may take $\tilde{C}_1$ independent of $y_0$
even if we replace $f$ with $f(\,\cdot\, +y_0)$.)

Let $\tilde{C}_1$ as in (\ref{eq.def.C_1_2}) and for this $\tilde{C}_1$
define $\delta$ 
as in Lemma \ref{lem.delta.cont}.
Choose $\beta >1$ arbitrarily and set $\rho =\beta/\delta$,
where $\delta$ is given in Lemma \ref{lem.delta.cont}.
As in (\ref{prop.choice.omega}), 
there exists a constant $c=c(\kappa_0,\rho,p,|f|_{\infty})$
such that
$\omega(s,t) =c\hat\omega(s,t)$ satisfies
that
, for
all $j=1,\ldots,[p]$ and  $(s,t) \in \triangle$,
\begin{eqnarray}
|(\kappa^{-1}X, \Lm)^j_{s,t}| \le  \frac12 \omega(s,t)^{j/p},
\qquad
| \Gamma_{\rho/\kappa,\rho,1,1} \tilde{K}(1)^j_{s,t} |
\le  \bigl(  \rho^p \omega(s,t)  \bigr)^{j/p}.
\label{ineq.rho.ome66}
\end{eqnarray}

Now we will show that, for $T_1 \in (0,1]$
such that $\omega(0,T_1) \le 1$,
it holds
on the restricted time interval $[0,T_1]$ that
\begin{equation}
\bigl|\Gamma_{\rho/\kappa,\rho,1,1}  \tilde{H}(n)^j_{s,t}  \bigr| 
\le \bigl(  \rho^p \omega(s,t)  \bigr)^{j/p},
\qquad
j=1,\ldots,[p], \quad (s,t) \in \triangle_{[0,T_1]}.
\label{ineq.rho.ome7}
\end{equation}

We use induction.
The case $n=1$ was already shown since $\tilde{H}(1)=\tilde{K}(1)$.
Using (\ref{eq.def.C_1_2}) for $\rho^p \omega$,
we have
\begin{eqnarray}
\bigl|
\Gamma_{1/\kappa,1,1,1}
\int_s^t  \tilde{\Phi}_{\beta^{n-1}} ( \Gamma_{\rho,\rho,1,1}  \tilde{H}(n))
d\Gamma_{\rho,\rho,1,1}  \tilde{H}(n)^j 
\bigr|
&=&
\bigl|
\Gamma_{\rho/\kappa,\rho,\rho,\rho}
\int_s^t  \tilde{\Phi}_{\beta^{n-1}} (  \tilde{H}(n))
d  \tilde{H}(n)^j 
\bigr|
\nn\\
&\le& 
(\tilde{C}_1 \rho^p \omega(s,t))^{j/p}.
\label{ineq.kk5}
\end{eqnarray}
By projection onto the ${\cal V}^{\oplus 2}$-component,
\begin{eqnarray}
\bigl| \pi_{ {\cal V}^{\oplus 2}}
\Gamma_{\rho/\kappa,\rho,\rho,\rho}
\int_s^t  \tilde{\Phi}_{\beta^{n-1}} (  \tilde{H}(n))
d \tilde{H}(n)^j 
\bigr|
=
|\rho^j (\kappa^{-1}X, \Lm)^j_{s,t}|
\le 
\frac12 ( \rho^p \omega(s,t))^{j/p}
\label{ineq.kk6}
\end{eqnarray}
We may use Lemma \ref{lem.delta.cont}
for (\ref{ineq.kk1}) and (\ref{ineq.kk2})
to obtain 
\begin{eqnarray}
\bigl| \Gamma_{1,1,\delta\beta^{-1} ,\delta }
\Gamma_{\rho/\kappa,\rho,\rho,\rho}
\int_s^t  \tilde{\Phi}_{\beta^{n-1}} (  \tilde{H}(n))
d  \tilde{H}(n)^j 
\bigr|
\le 
( \rho^p \omega(s,t))^{j/p}
\nn
\end{eqnarray}
From this, we see that (\ref{ineq.rho.ome7}) for $n+1$.
Hence we have shown (\ref{ineq.rho.ome7}) for any $n$.
From (\ref{ineq.rho.ome7}), it is easy to see that
\begin{equation}
\bigl|\Gamma_{1/\kappa,1,1,1}  \tilde{H}(n)^j_{s,t}  \bigr| 
\le \bigl(  \rho^p \omega(s,t)  \bigr)^{j/p},
\qquad
j=1,\ldots,[p], \quad (s,t) \in \triangle_{[0,T_1]}.
\label{ineq.rho.ome88}
\end{equation}

In the same way as in (\ref{ineq.kk3})--(\ref{ineq.conv3}),
we obtain from (\ref{ineq.rho.ome88}) that
\begin{eqnarray}
\bigl| 
\Gamma_{1/\kappa,1,1}  \tilde{Z}^j_{s,t}
\bigr|
\le 
( \tilde{C}_3  \omega(s,t))^{j/p},
\qquad
j=1,\ldots,[p], \quad (s,t) \in \triangle_{[0,T_1]}.
\label{ineq.conv4}
\end{eqnarray}
for some constant $\tilde{C}_3>0$ which depends only on 
$\kappa_0, p$, and $M(f;[p]+1)$ (and the choice of $\beta> 1$).

Note that (\ref{ineq.conv4}) is the desired inequality
(on the restricted interval).
By prolongation of solution we can prove the lemma. 
\QED

For $X$ and $\Lm$ as above,
set $Q=Q(X,\Lm):=Y_{X+\Lm}-Y_{\Lm}$. 
Clearly, 
\[
dQ= f(Y_{X+\Lm})dX + [f(Y_{X+\Lm})-f(Y_{\Lm})] d\Lm.
\]
\begin{lem}\label{lem.geki2}
Let $\kappa_0 >0$. 
Assume that a control function $\hat\omega$,
 $\Lm \in C_{0,q}({\cal V})$ and $X \in G\Omega_p({\cal V})$,
and $\kappa \in [0, \kappa_0]$
satisfy that
\[
|\Lm^1_{s,t}| \le \hat\omega (s,t)^{1/q},
\qquad
|X^j_{s,t}| \le (\kappa \hat\omega (s,t))^{j/p}
\qquad
j=1,\ldots,[p], \quad (s,t)\in \triangle.
\]
Then, there is a positive constant $C$ which depends only on 
$\kappa_0, \hat\omega(0,1),p, M(f;[p]+2)$
such that
\[
\Bigl|
\bigl(
\kappa^{-1}X, \Lm, Y_{X+\Lm}, Y_{\Lm},\kappa^{-1}Q
\bigr)^j_{s,t}
\Bigr|
\le (C \hat\omega(s,t))^{j/p},
\qquad
j=1,\ldots,[p],
\quad
(s,t)\in \triangle.
\]
Here, the left hand side denotes the $j$th level path 
of a ${\cal V}^{\oplus 2} \oplus {\cal W}^{\oplus 3}$-valued 
geometric rough path.
\end{lem}

\Proof
In this proof, the positive constant $C$ may change
from line to line.
As before we may assume that $X \in C_{0,q}({\cal V})$.
Recall that $\Lm \in C_{0,q}({\cal V}) \mapsto Y_{\Lm} \in C_{0,q}({\cal W})$
is continuous 
and 
there exists a constant $C>0$ such that
$|(Y_{\Lm})^1_{s,t} | \le C\hat\omega(s,t)^{1/q}$.
Combining this with Lemma \ref{lem.geki}, we have 
\[
\Bigl|
\bigl(
\kappa^{-1}X, \Lm, Y_{X+\Lm}, Y_{\Lm}
\bigr)^j_{s,t}
\Bigr|
\le (C \hat\omega(s,t) )^{j/p},
\qquad
j=1,\ldots,[p],
\quad
(s,t)\in \triangle.
\]

From Lemma \ref{lem.geki}, we easily see that, for some constant $C>0$,
\[
\Bigl|
\bigl(
\kappa^{-1}X, \Lm, Y_{X+\Lm}, Y_{\Lm}, \kappa^{-1}\int f(Y_{X+\Lm})dX
\bigr)^j_{s,t}
\Bigr|
\le (C \hat\omega(s,t))^{j/p},
\quad
j=1,\ldots,[p],
(s,t)\in \triangle.
\]

Note that from the local 
Lipschitz continuity of (the first level path of)
the It\^o map (Theorem \ref{thm.loc.Lip}), 
we see that, for some constant $C>0$,
\[
\Bigl|  \kappa^{-1}
\int_s^t [f(Y_{X+\Lm})-f(Y_{\Lm})] d\Lm
\Bigr|
\le C  \hat\omega(s,t)^{1/q},
\qquad
j=1,\ldots,[p],
\quad
(s,t) \in \triangle.
\]
Here, the left hand side is the Young integral.
From these, we can easily obtain the theorem.
\QED


%

%
\section
{A stochastic Taylor-like expansion}

\subsection
{Estimates for ordinary  terms in the expansion}
In this section we will estimate ordinary terms
in the stochastic Taylor-like expansion for It\^o maps.
Let $p \ge 2$ and $1 \le q <2$ with $1/p+1/q>1$
and let ${\cal V}, \hat{ {\cal V}}, {\cal W}$ be real Banach spaces.
Let $\sigma \in C_b^{\infty}([0,1]\times {\cal W}, L({\cal V},{\cal W}))$
and 
$b \in C_b^{\infty}([0,1]\times {\cal W}, L(\hat{ {\cal V}},{\cal W}))$.
Here, $[0,1]\times {\cal W}$ is considered as a subset 
of the direct sum ${\bf R} \oplus {\cal W}$.
We will consider the following ODE:
for $\ve >0$,
$X \in G\Omega_p ( {\cal V} )$, and $\Lambda \in C_{0,q}(\hat{{\cal V}})$,
\begin{equation}\label{eq.ode_formal}
dY^{(\ve)}_t
=
\sigma(\ve,Y^{(\ve)}_t) \ve dX_t
+
b(\ve, Y^{(\ve)}_t) d\Lambda_t,
\qquad 
Y^{(\ve)}_0=0
\end{equation}
Note that if $X$ is lying above an element of $C_{0,q}( {\cal V})$,
then (\ref{eq.ode_formal}) makes sense in the $q$-variational setting.

More precisely, the above equation  (\ref{eq.ode_formal})
can be formulated as follows.
Define $\tilde\sigma$
by
$\tilde\sigma= \sigma \circ p_1 +b \circ p_2$,
where $p_1$ and $p_2$ are canonical projection from 
${\cal V} \oplus \hat{{\cal V}}$
onto the first and the second component, respectively.
Then, $\tilde{\sigma} 
\in C_b^{\infty}([0,1]\times
 {\cal W}, L( {\cal V} \oplus \hat{{\cal V}},{\cal W}))$.
We consider the It\^o map $\Phi^{\ve}: 
G\Omega_p ( {\cal V} \oplus \hat{{\cal V}}) \to G\Omega_p ( {\cal W})$
which corresponds to 
the coefficient $\tilde{\sigma}(\ve, \,\cdot\,)$
with the initial condition $0$.
If $X \in C_{0,q}( {\cal V})$ and $\Lambda \in C_{0,q}(\hat{\cal V})$, 
then $(X, \Lambda) \in C_{0,q}( {\cal V} \oplus \hat{{\cal V}})$.
This map naturally extends to a continuous map from 
$G\Omega_p (  {\cal V}) \times C_{0,q}(\hat{{\cal V}})$ 
to $G\Omega_p (  {\cal V} \oplus \hat{ {\cal V}})$
(see Corollary \ref{cor.loc.Lip}). 
The precise meaning of  (\ref{eq.ode_formal})
is that $Y^{(\ve)}=\Phi^{\ve} ( (\ve X,\Lambda))_1$.

\begin{remark}
In Azencott \cite{az},
he treated differential equations with the coefficients
of the form $\sigma (\ve, t, y)$ and $b(\ve, t, y)$.
ODE (\ref{eq.ode_formal}), however, includes such cases.
In order to see this, set ${\cal W}' ={\cal W} \oplus {\mathbb R}$,
$\hat{\cal V}' =\hat{\cal V} \oplus {\mathbb R}$,
$\Lambda'_t =(\Lambda_t,t)$,
and add to (\ref{eq.ode_formal}) the following trivial equation;
$dY'_t =dt$.
\end{remark}

We set $Y^0=\Phi^{0} ( ({\bf 1},\Lambda))_1$,
where ${\bf 1}=(1,0,\ldots,0)$ 
is the unit element in the truncated tensor algebra
(which is regarded as a constant rough path).
Note that $\Lambda \in C_{0,q}(\hat{{\cal V}}) 
\mapsto Y^0 \in C_{0,q}({\cal W})$
is locally Lipschitz continuous (see \cite{lq}).

We will expand $Y^{(\ve)}$ in the following form:
\[
Y^{(\ve)} \sim Y^0+\ve Y^1+ \ve^2 Y^2+ \ve^3 Y^3+\cdots
\qquad
\text{as $\ve \searrow 0$.}
\]
(Note that $Y^k$ does NOT denote the $k$th level path of $Y$.
The  $k$th level path of $Y^j$ will be denoted by $(Y^j)^k$.
Similar notations will be used for $I^k$ and $J^k$ below.
This may be a little confusing. Sorry.)
By considering a (formal) Taylor expansion 
for $\sigma(\ve,Y^{(\ve)}_t)$ and $b(\ve,Y^{(\ve)}_t)$,
we will find explicit forms of $Y^n~(n \in {\mathbb N})$ as follows. 
(Or equivalently, 
we may formally operate $(n!)^{-1} (d/d\ve)^n$ at $\ve =0$
on the both sides of (\ref{eq.ode_formal})).
\begin{eqnarray}
dY^n_t -\de_y b(0, Y^0_t) \la Y^n_t, d \Lambda_t\ra
&=&  dI^n_t +dJ^n_t,
\qquad (n \in {\mathbb N})
\label{def.ord.term}
\end{eqnarray}
where $I^n=I^n(X,\Lm)$ and $J^n=J^n(X,\Lm)$ are given by
\begin{eqnarray}\label{def.iandj2}
dI^1_t &=& \sigma(0,Y^0_t) dX_t,
\qquad
dJ^1_t = \de_{\ve} b(0, Y^0_t) d\Lm_t,
\nn\\
\end{eqnarray}
with $I^1_0 =J^1_0 =0$,  and, for $n=2,3,\ldots$,
\begin{eqnarray}
dI^n_t &=&
\sum_{j=0}^{n-2}\sum_{k=1}^{n-1-j}
\sum_{(i_1, \ldots, i_k) \in S^{n-1-j}_{k}}
\frac{1}{j!k!}
\de_{\ve}^j\de_{y}^k \sigma(0, Y^0_t)
\la  Y^{i_1}_t, \ldots, Y^{i_k}_t, dX_t \ra
\nn\\
&&
+
\frac{1}{(n-1)!} \de_{\ve}^{n-1}\sigma(0, Y^0_t)dX_t
\nn\\
dJ^n_t &=&
\sum_{j=1}^{n-1}\sum_{k=1}^{n-j}
\sum_{(i_1, \ldots, i_k) \in S^{n-j}_{k}}
\frac{1}{j!k!}
\de_{\ve}^j\de_{y}^k b(0, Y^0_t)
\la  Y^{i_1}_t, \ldots, Y^{i_k}_t, d\Lm_t \ra
\nn\\
&&
+
\sum_{k=2}^{n}
\frac{1}{k!}\de_{y}^k b(0, Y^0_t)
\la  Y^{i_1}_t, \ldots, Y^{i_k}_t, d\Lm_t \ra
+\frac{1}{n!} \de_{\ve}^{n} b(0, Y^0_t )d\Lm_t
\label{def.iandj}
\end{eqnarray}
with $I^n_0=J^n_0=0$.
Here, $\de_{\ve}$ and $\de_y$ denote the partial Fr\'echet derivatives
in $\ve$ and in $y$, respectively,
and 
\[
S^{m}_{k} =\{ (i_1, \ldots, i_k) \in  {\mathbb N}^k=\{1,2,\ldots\}^k ~|~
i_1+ \cdots +i_k=m
\}.
\]

Now we define functions which appear on the right hand sides
of (\ref{def.iandj2}) and  (\ref{def.iandj}).
Let ${\cal X}_{n-1} ={\cal V} \oplus \hat{\cal V} \oplus {\cal W}^{\oplus n}$.
An element in ${\cal X}_{n-1}$ is denoted by
$v=(x,\hat x ; y^0,y^1, \ldots,y^{n-1})$.
Partial Fr\'echet derivatives are denoted by $\de_{x}$, $\de_{y^1}$, etc.
and the projection from ${\cal X}_{n-1}$ onto each components
are denoted by $p_{x}$, $p_{y^1}$, etc.
Set $f_1, g_1 \in C^{\infty}_{b,loc} ({\cal X}_0, L({\cal X}_0, {\cal W}))$ by   
\begin{equation}\label{def.fncFG2}
f_1( y^0)=\sigma(0,y^0) \circ p_{x},
\qquad
g_1( y^0)=\de_{\ve} b(0, y^0) \circ p_{\hat{x}}.
\end{equation}
For $n=2,3,\ldots$,
set $f_n, g_n 
\in C^{\infty}_{b,loc} ({\cal X}_{n-1}, L({\cal X}_{n-1}, {\cal W}))$ by
\begin{eqnarray}
f_n(  y^0,\ldots,y^{n-1})
&=&
\Bigl[
\sum_{j=0}^{n-2}\sum_{k=1}^{n-1-j}
\sum_{(i_1, \ldots, i_k) \in S^{n-1-j}_{k}}
\frac{1}{j!k!}
\de_{\ve}^j\de_{y}^k \sigma(0, y^0)
\la  y^{i_1}, \ldots, y^{i_k}, \,\cdot\, \ra
\nn\\
&&
+
\frac{1}{(n-1)!} \de_{\ve}^{n-1}\sigma(0, y^0)
\Bigr] \circ p_x,
\nn\\
g_n(  y^0,\ldots,y^{n-1})
&=&
\Bigl[
\sum_{j=1}^{n-1}\sum_{k=1}^{n-j}
\sum_{(i_1, \ldots, i_k) \in S^{n-j}_{k}}
\frac{1}{j!k!}
\de_{\ve}^j\de_{y}^k b(0, y^0)
\la  y^{i_1}, \ldots, y^{i_k}, \,\cdot\, \ra
\nn\\
&&
+
\sum_{k=2}^{n}
\frac{1}{k!}\de_{y}^k b(0, y^0)
\la  y^{i_1}, \ldots, y^{i_k}, \,\cdot\, \ra
+
\frac{1}{n!} \de_{\ve}^{n} b(0, y^0 )
\Bigr] \circ p_{\hat{x}}.
\label{def.fncFG}
\end{eqnarray}
Clearly, the functions $f_n$ and $g_n$ are actually 
independent of $x$ and $\hat x$.

%
%
\begin{lem}\label{lem.est.derFG}
Let $f_n$ and $g_n$ be as above and $\delta>0, C>0, r \in {\mathbb N}$.
For ${\bf \xi}=(\xi_1, \dots, \xi_r ) \in \{0,\ldots,n-1\}^r$,
we set $|\xi|=\sum_{k=1}^r \xi_k$.
Then, on the following set 
$$
\bigl\{(y^0, \ldots,y^{n-1})~|~ 
\text{ $|y^i| \le C(1+\delta)^i$ for $0 \le i \le n-1$} \bigr\},$$
it holds that, for any  
${\bf \xi}$ such that 
$|\xi| \le n-1$,
\begin{eqnarray}
\bigl| \de^{r}_{\xi }
f_n( x,\hat x ; y^0,\ldots,y^{n-1}) \bigr|
&\le&
C' (1+\delta)^{n-1- |\xi| },
\nn
\end{eqnarray}
Here, $\de^{r}_{\xi }= \de_{y^{\xi_1}} \cdots \de_{y^{\xi_r}}$
and $C'$ is a positive constant 
independent of $\delta$. 
If $|\xi| > n-1$, then 
$\de^{r}_{\xi}  f_n =0$.

Similarly, it holds on the same set that, for any 
$\xi $ such that 
$|\xi| \le n$, 
\begin{eqnarray}
\bigl| \de^{r}_{ \xi }
g_n( x,\hat x ; y^0,\ldots,y^{n-1}) \bigr|
&\le&
C' (1+\delta)^{n- |\xi|}.
\nn
\end{eqnarray}
If $|\xi| > n$, then 
$\de^{r}_{\xi} g_n =0$.
\end{lem}

\Proof
This lemma can be shown by straight forward computation
since $f_n$ and $g_n$ are (i) $C_b^{\infty}$ in $y^0$-variable
and (ii) ``polynomials'' in $(y^1, \ldots, y^{n-1})$-variables.
\QED

%
%
In fact, we can compute
$\de^{r}_{\xi}  f_n
= \de_{y^{\xi_1}} \cdots \de_{y^{\xi_r}} f_n$  explicitly as follows.
\begin{lem}\label{lem.diffFG}
Let $f_n~(n=1,2,\ldots)$ be as above and $r=1,2,\ldots$.
For ${\bf \xi}=(\xi_1, \dots, \xi_r ) \in \{0,\ldots,n-1\}^r$,
set
$\mu=\sharp \{ k~|~1\le k \le r, ~\xi_k \neq 0 \}$. 
%
%
%
Then, if $|\xi| \le n-1$,
\begin{eqnarray*}
\lefteqn{
(\de_{\xi}^r   f_n)(y^0, \ldots,y^{n-1})
}
\\
&=&
\Bigl[
\sum_{j=0}^{n-2-|\xi|}
\sum_{k=\mu+1}^{n-1-j-|\xi| }
\sum_{(i_1, \ldots, i_{k-\mu}) \in S^{n-1-j-|\xi|}_{k}}
\frac{1}{j!(k-\mu)!}
\de_{\ve}^j
\de_{y}^{k+r-\mu} \sigma(0, y^0)
\la  y^{i_1}, \ldots, y^{i_{k-\mu}}, \,\cdot\, \ra
\nn\\
&&
+
\frac{1}{(n-1-|\xi|)!} \de_{\ve}^{n-1-|\xi|}
\de_{y}^r
\sigma(0, y^0)
\Bigr] 
\circ p_{(\xi_1,  \ldots, \xi_r, x)}.
\end{eqnarray*}
Here, the right hand side is regarded as in 
$L^{r+1}( {\cal X}_{n-1},\ldots,{\cal X}_{n-1}; {\cal W})$.
Note that if $|\xi|= n-1$
the first term on the right hand is regarded as zero.
If $|\xi| > n-1$,
$\de_{\xi}^r   f_n=0$.
\end{lem}

\Proof
We give here a slightly heuristic proof.
However, since the difficulty of this lemma lies only in algebraic part,
it does not cause a serious trouble.

From the Taylor expansion for $\sigma$ 
\begin{eqnarray}
\ve \sigma(\ve, y^0+\Delta y)
&\sim&
\ve \sum_{j,k} 
\frac{ \ve^{j} }{j!k!}
\de_{\ve}^j\de_{y}^{k}\sigma(0, y^0) 
\la\overbrace{ \Delta y,\ldots,\Delta y }^{k} , \,\cdot\, \ra,
\nn\\
\Delta y &\sim& \ve^1 y^1+ \ve^2 y^2+ \ve^3 y^3+\cdots,
\qquad
\text{ as $\ve \searrow 0$.}
\label{eq.form.tay}
\end{eqnarray}
Then, we get a linear combination of the terms of the form
\[
\ve^{j+1+ i_1+\cdots +i_k}
\de_{\ve}^{j}\de_{y}^{k}\sigma(0, y^0) 
\la y^{i_1} ,\ldots, y^{i_k} , \,\cdot\,\ra.
\]
(In this proof we say  
the above term is of order $j+1+i_1+\cdots+i_k$.)
Recall that the definition of $f_n$ is the sum of terms of order $n$
in the right hand side of (\ref{eq.form.tay}).

Let us first consider $\de_{y^s}f_n~(s \neq 0)$.
Then, if $n-s \ge 0$,
$\de_{y^s}f_n(y^0,\ldots,y^{n-1})$
is the sum of terms of order $n-s$ of the $\de_{y^s}$-derivative 
of (\ref{eq.form.tay}), which is given by 
\[
 \sum_{j,k} 
\frac{ \ve^{j+1+s} }{j!(k-1)!}
\de_{\ve}^j\de_{y}^{k}\sigma(0, y^0) 
\la \overbrace{\Delta y,\ldots,\Delta y}^{k-1}, \,\cdot\, \ra,
\qquad
\Delta y \sim \ve^1y^1+\ve^2 y^2+\cdots.
\]
Picking up terms of order $n$, we easily see 
$\de_{y^s}f_n(y^0,\ldots,y^{n-1})$ is given as in the statement of
this lemma.
The case for $\de_{y^0}f_n$ is easier.

Thus, we have shown the lemma for $r=1$.
Repeating this argument, we can show the general case ($r \ge 2$).
\QED

We set some notations for iterated integrals.
Let ${\cal A}^i$ be real Banach spaces
and let $\phi^i$ be  ${\cal A}^i$-valued paths
 ($1\le i \le n$).
Define 
\[
{\cal I}^n
[\phi^1, \ldots, \phi^n ]_{s,t}
= 
\int_{s < u_1 < \cdots < u_n  < t}   
d\phi^1_{u_1} \otimes \cdots \otimes d\phi^n_{u_n}
\in
{\cal A}_1 \otimes \cdots \otimes {\cal A}_n
\]
whenever possible. 
(For example, when $\phi^i \in C_{0,q}({\cal A}_i)~(1\le i \le n)$ for 
some $1\le q<2$.)

Let ${\cal B}^i$ be real Banach spaces ($1 \le i \le n$) and 
${\cal B}=\otimes_{i=1}^n {\cal B}_i$.
For $\pi \in \Pi_n$ and 
$b_1  \otimes \cdots \otimes b_n \in {\cal B}$, we write
\[
\pi(b_1  \otimes \cdots \otimes b_n) 
= (b_{\pi^{-1}(1)}, \ldots, b_{\pi^{-1}(n)})
\in {\cal B}_{\pi^{-1}(1)} \otimes \cdots \otimes {\cal B}_{\pi^{-1}(n)}
\]

Let ${\cal C}=\oplus_{i=1}^m {\cal A}_i$
and consider ${\cal C}^{\otimes n}$.
The $(i_1, \ldots,i_n)$-component of $\eta \in {\cal C}^{\otimes n}$
is denoted by $\eta^{(i_1, \ldots,i_n)}
\in {\cal A}_{i_1} \otimes \cdots\otimes {\cal A}_{i_n}$.
Clearly, $\sum_{1\le i_1, \ldots,i_n \le m} 
\eta^{(i_1, \ldots,i_n)}=\eta$.
Let $\psi =(\psi^1,\ldots,\psi^m)$ be a nice path in ${\cal C}$.
The $n$th level path of the rough path lying above
$\psi$ is ${\cal I}^n
[\psi, \ldots, \psi ]$.
The action of $\pi \in \Pi_n$ in component form is given by 
\begin{equation}\label{eq.compo}
\bigl(
\pi {\cal I}^n[\psi, \ldots, \psi ]_{s,t}
\bigr)^{(i_1, \ldots,i_n)}
=
\pi {\cal I}^n[\psi^{i_{\pi(1)}}, \ldots, \psi^{i_{\pi(n)}} ]_{s,t}.
\end{equation}
This equality can be verified by straightforward computation.

\begin{sloppypar}
Now we state our main theorem in this subsection.
In the following we set
\[
\xi(X) :=\sum_{j=1}^{[p]}  \|X^{j}\|_{p/j}^{1/j}
\qquad
\mbox{
for $X \in C_{0,q}({\cal V})$ (or $X \in G\Omega_{p}({\cal V})$).
}
\]
Here, $X^j$ denotes the $j$th level path of (the  rough path 
lying above) $X$.
Clearly, $\xi (rX)=|r| \xi (X)$ for $r \in {\mathbb R}$.
In the following we set $\nu(-2)=1$, $\nu(-1)=0$, and
$\nu(i)=i$ for $i \ge 0$.
\begin{thm}\label{thm.main.estI}
The map $(X, \Lm) \mapsto (X, \Lm,Y^0,\ldots,Y^n)$
extends to a continuous map from 
$G\Omega_p({\cal V})\times C_{0,q}(\hat{\cal V})$
to $G\Omega_p({\cal X}_{n})$.
Moreover, for any $X \in C_{0,q}({\cal V})$ and $\Lm 
\in C_{0,q}(\hat{\cal V})$,
there exists a control function $\omega=\omega_{X,\Lm}$
such that the following {\rm (i)} and {\rm (ii)}
hold:
\\
\noindent
{\rm (i)}
For any $(s,t) \in \triangle$, 
$X \in C_{0,q}({\cal V})$, $\Lm \in C_{0,q}(\hat{\cal V})$,
$j=1,\ldots,[p]$,
and $i_1,\ldots,i_j \in \{-2,-1,\ldots,n\}$, 
it holds that
\begin{eqnarray}
\bigl|
{\cal I}^j
[Y^{i_1}, \ldots, Y^{i_j} ]_{s,t}
\bigr|
\le 
(1+ \xi(X))^{\nu(i_1)+\cdots+\nu(i_j)} \omega(s,t)^{j/p}
\label{ineq.interI}
\end{eqnarray}
Here, for notational simplicity,
we set $Y^{-2}=X$, $Y^{-1}=\Lm$.
\\
\noindent
{\rm (ii)}
For any $r>0$, there exists a constant $c=c(r)>0$ such that
\[
\sup\{ \omega(0,1) ~|~ 
\mbox{ $X \in C_{0,q}({\cal V})$,
$\Lm \in C_{0,q}(\hat{\cal V})$ with $\|\Lm \|_q \le r$}  
\}
\le c.
\]
\end{thm}
\end{sloppypar}

\Proof
In this proof, $c$ and $\omega$ may change from line to line
and we will denote $\delta=\xi(X)$.
We will use induction.
The case $n=0$ is easy, since the map
 $\Lm \in C_{0,q}(\hat{\cal V}) \mapsto Y^0 \in C_{0,q}({\cal W})$
is locally Lipschitz continuous (see \cite{lq}).
Now we assume the statement of the theorem holds for $n-1$
and will prove the case for $n$.

Set $Z^{n-1}=(X, \Lm,Y^0,\ldots,Y^{n-1})$ for $n \in {\mathbb N}$.
Then, it is obvious that
\[
(X, \Lm,Y^0,\ldots,Y^{n-1},I^n)=
\int ({\rm Id}_{{\cal X}^{n-1}} \oplus f_n) (Z^{n-1})dZ^{n-1}.
\]
By using (\ref{defY})--(\ref{eq.2no7}),
we will estimate the almost rough path $\Xi \in A\Omega_p({\cal X}^n)$, 
which defines the integral on the right hand side.

Let $1 \le k\le [p]$.
We consider the ${\bf i}=(i_1, \ldots,i_k)$-component of $\Xi^k$, 
where $-2 \le i_j \le n$ for all $j=1,\ldots,k$.
Set ${\cal N}({\bf i})=\{ j~|~i_j =n\}$.
Note that, if $j \notin {\cal N}({\bf i})$ (equivalently, if $i_j \neq n$), 
then ${\bf l}=(l_1,\ldots,l_k)$ in the sum of \ref{defY}
must satisfy $l_j=1$.

We will fix such an $\bf l$.
For $j \notin {\cal N}({\bf i})$,
set $L_j:=p_{i_j}$ (the projection onto the $i_j$-component)
and ${\bf m}_j={\bf m}'_j=i_j$.
For $j \in {\cal N}({\bf i})$, set
${\bf m}_j=(m^1_{j}, \ldots, m^{l_j -1}_{j})$,
${\bf m}'_j=(m^1_{j}, \ldots, m^{l_j -1}_{j},-2)$,
and 
\begin{equation}\label{eq.defLm}
L_j= \sum_{ {\bf m}_j \in \{0,\ldots,n-1  \}^{l_j -1} }  
\de^{l_j -1}_{ {\bf m}_j}
f_n (Z^{n-1}_s)
=:
\sum_{ {\bf m}_j \in \{0,\ldots,n-1  \}^{l_j -1} }  L_j^{{\bf m}_j}.
\end{equation}
Then, from (\ref{defY}), 
\begin{eqnarray}
[\Xi_{s,t}^k]^{({\bf i})}
&=&
\sum_{{\bf l}} 
\Bigl[
\sum_{ {\bf m}_j \in \{0,\ldots,n-1  \}^{l_j -1} }
(L_1 \otimes \cdots \otimes L_k )
\Bigl\la
\sum_{\pi \in \Pi_{{\bf l}}}
[\pi (Z^{n-1})^{|{\bf l}|}_{s,t} ]^{({\bf m}'_1,\ldots,{\bf m}'_k)}
\Bigr\ra
\Bigr].
\label{eq.tsumeta}
\end{eqnarray}
The sum is over such ${\bf l}$'s.
From (\ref{eq.compo}) and the assumption of induction, we easily see that
\begin{eqnarray}
\bigl|
[\pi (Z^{n-1})^{|{\bf l}|}_{s,t} ]^{({\bf m}'_1,\ldots,{\bf m}'_k)}
\bigr|
\le
c (1+\delta )^{\nu ( {\bf m}'_1,\ldots,{\bf m}'_k )}  \omega(s,t)^{k/p},
\end{eqnarray}
where 
$\nu ( {\bf m}'_1,\ldots,{\bf m}'_k ) 
:=\sum_{j=1}^k \sum_{r=1}^{l_j} \nu(m_j^{r})$.
Combining this with Lemma \ref{lem.est.derFG},
we have that
\begin{equation}\label{ineq.apXi1}
\bigl| [\Xi_{s,t}^k]^{({\bf i})} \bigr|
\le 
c (1+\delta )^{\nu ({\bf i} )}  \omega(s,t)^{k/p},
\qquad
k=1,\ldots,[p], \,
(s,t) \in \triangle, 
\end{equation}

Next we estimate ${\bf i}$-component of 
$\Xi^k_{s,t}- ( \Xi_{s,u} \otimes \Xi_{u,t})^k$
for $s<u<t$.
For that purpose, it is sufficient to estimate $R_l(x,y)$ 
in (\ref{eq.2no5}) for $f={\rm Id}_{{\cal X}^n-1} \oplus f_n$.
If $i \neq n$, $i$th component of $R_l(x,y)$
(that is equal to $R_l(x,y)$ for the projection $p_i$)
vanishes.
From Lemma \ref{lem.est.derFG}, the $n$th component 
$R(f_n)_l(Z^{n-1}_s, Z^{n-1}_u)$ satisfies the following:
for ${\bf m}_j$ as above, 
\begin{eqnarray}
\Bigl|
\int_0^1 d\theta
\frac{(1-\theta)^{[p] -l} }{ ([p] -l)!}
D^{[p] -l+1} \de^{l-1}_{{\bf m}_{j} }
f_n \bigl( Z^{n-1}_s +\theta (Z^{n-1})^1_{s,u} \bigr)
\la (Z^{n-1})_{s,u}^1)^{\otimes ([p] -l+1)}  \ra
\Bigr|
\nn\\
\le
c(1+\xi(X))^{n-1-(m^1_{j}+\cdots +m_{j}^{l-1} ) },
\qquad
\mbox{if $ m^1_{j}+\cdots +m_{j}^{l-1} \le n-1$.} 
\nn
\end{eqnarray}
Here, the left hand side is regarded as a multilinear map
from $Y^{m^1_{j}} \times \cdots \times Y^{m^{l-1}_{j}} \times X$.
If $ m^1_{j}+\cdots +m_{j}^{l-1} > n-1$, then the left hand side vanish.

Denoting 
by $\hat{L}_j$ the left hand side of the above inequality,
we can do the same argument as in (\ref{eq.defLm})--(\ref{ineq.apXi1})
to obtain that
\begin{equation}\label{ineq.apXi2}
\Bigl| 
\bigl[
\Xi_{s,t}^k -(\Xi_{s,u} \otimes \Xi_{u,t})^k
\bigr]^{({\bf i})} 
\Bigr|
\le 
c (1+\delta )^{\nu ({\bf i} )}  \omega(s,t)^{([p] +1)/p}
\end{equation}

From (\ref{ineq.apXi1})--(\ref{ineq.apXi2}), 
$(X, \Lm,Y^0,\ldots,Y^{n-1},I^n)$ satisfies
a similar inequality to (\ref{ineq.interI})
(with $Y^n$ in (\ref{ineq.interI}) 
being replaced with $I^n$).

Now, by the Young integration theory, we see that
$
| J^n_t -J^n_s| \le  c (1+\delta)^{n}  \omega(s,t)^{1/q}
$,
which implies that
$(X, \Lm,Y^0,\ldots,Y^{n-1},I^n+J^n)$ satisfies
a similar inequality to (\ref{ineq.interI})
(with $Y^n$ in (\ref{ineq.interI})
being replaced with $I^n+J^n$).

Set $d\Omega_t=\de_yb(0,Y^0_t ) \la \,\cdot\, ,d\Lm_t\ra$
and set 
$M=M(\Lm)$ 
by $dM_t= d\Omega_t \cdot M_t$ with $M_0={\rm Id}_{{\cal W}}$.
It is easy to see that $\Lm \in C_{0,q}({\cal V}) 
\mapsto M\in {\cal C}_q (L({\cal W},{\cal W}))$
is continuous
(See Proposition \ref{prop.q.ode}).
Also set $\hat{M}={\rm Id}_{{\cal X}_{n-1} } 
\oplus M \in {\cal C}_q (L( {\cal X}_{n-1},{\cal X}_{n-1}  ))$.
Then, applying Duhamel's principle (Corollary \ref{cor.duhamel})
for $\hat{M}$ and $(r X, Y^0,r Y^1,\ldots,r^{n-1} Y^{n-1},r^n(I^n+J^n))$
with $1/r=(1+\delta )$,
we obtain 
$(r X, Y^0,r Y^1,\ldots,r^{n-1} Y^{n-1},r^n Y^n)$.
This completes the proof.
%
%
\QED

\subsection{Estimates for remainder terms in the expansion}
In this subsection we give estimates for remainder terms
in the stochastic Taylor-like expansion.
We keep the same notations as in the previous subsection.
If $X \in C_{0,q}({\cal V})$ and $\Lm \in C_{0,q}(\hat{\cal V})$,
then 
\[
Q^{n+1,(\ve)}=Q^{n+1,(\ve)}(X,\Lm)
:=
Y^{(\ve)}- (Y^0+\ve Y^1+ \cdots+ \ve^n Y^n),
\qquad
\ve \in [0,1] 
\]
is clearly well-defined.
We prove that the correspondence $(X,\Lm) \mapsto Q^{n+1,(\ve)}(X,\Lm)$
extends to a continuous map from $G\Omega_p ({\cal V}) 
\times C_{0,q}(\hat{\cal V})$
to $G\Omega_p ({\cal W})$
and that $Q^{n+1,(\ve)}(X,\Lm)$ is a term of ``order $n+1$''.

For simplicity we assume that
$X \in C_{0,q}({\cal V})$ and $\Lm \in C_{0,q}(\hat{{\cal V}})$.
From (\ref{eq.ode_formal}) and (\ref{def.ord.term}),
we see that
\begin{eqnarray}
\lefteqn{
dQ^{n+1,(\ve)}_t-\de_y b(0, Y^0_t) \la  Q^{n+1,(\ve)}_t , d \Lambda_t\ra
= 
\sigma(\ve,Y^{(\ve)}_t) \ve dX_t
-\sum_{k=1}^n   \ve^k dI^k_t
}
\nn\\
&&+ 
b(\ve, Y^{(\ve)}_t) d\Lambda_t
-b(0, Y^{0}_t) d\Lambda_t
-\de_y b(0, Y^0_t) \la Y^{(\ve)}_t- Y^0_t , d \Lambda_t\ra
-\sum_{k=1}^n \ve^k  dJ^k_t.
\label{rep.rem.term}
\end{eqnarray}
Note that $I^k$ and $J^k~(k=1,\ldots,n)$ depends only on 
 $Y^0,Y^1,\ldots,Y^{n-1}$, but not on $Y^n$.

Set $\hat{\cal X}_{n} ={\cal V} \oplus \hat{\cal V}  
\oplus {\cal W}^{\oplus n+2}$.
An element in ${\cal X}_{n}$ is denoted by
$(x,\hat x ; y^{-1},y^0, \ldots,y^{n})$.
Then, in a natural way, 
$f_n, g_n 
\in C^{\infty}_{b,loc} 
(\hat{\cal X}_{n}, L( \hat{\cal X}_{n}, {\cal W}))$
and Lemmas \ref{lem.est.derFG} and \ref{lem.diffFG}
hold with trivial modification.

In the following theorem, we set for simplicity
\[
Z^{n,(\ve)}=
Z^{n,(\ve)}(X,\Lm) :=
(\ve X,\Lm,Y^{(\ve)},
Y^0,\ve Y^1, \ldots, \ve^{n-1} Y^{n-1}, Q^{n,(\ve)}).
\]
This is a $\hat{\cal X}_{n}$-valued path.
We also set $\hat{Y}^{i,(\ve)}= \ve^{i}Y^{i}$ for $0 \le i \le n-1$,
$\hat{Y}^{n,(\ve)}=Q^{n,(\ve)}$,
$\hat{Y}^{-1,(\ve)}=Y^{(\ve)}$, $\hat{Y}^{-2,(\ve)}=\Lm$,
 and $\hat{Y}^{-3,(\ve)}=\ve X$.
Then, 
$Z^{n,(\ve)}=(\hat{Y}^{-3},\ldots,\hat{Y}^{n})$.
We define 
$\nu(i)=i$ for $i \ge 0$, $\nu(-1)=\nu(-2)=0$, and
$\nu(-3)=1$.

\begin{thm}\label{thm.rem.est}
For each $n \in {\mathbb N}$ and $\ve \in [0,1]$, 
the map $(X,\Lm) \mapsto Z^{n,(\ve)}(X,\Lm)$
extends to a continuous map from 
$G\Omega_p( {\cal V})\times C_{0,q}(\hat{{\cal V}})$
to $G\Omega_p(\hat{\cal X}_{n})$.
Moreover, for any $X \in C_{0,q}({\cal V})$ and 
$\Lm \in C_{0,q}(\hat{\cal V})$,
there exists a control function $\omega=\omega_{X,\Lm}$
such that the following {\rm (i)} and {\rm (ii)}
hold:
\\
\noindent
{\rm (i)}
For any $(s,t) \in \triangle$, 
$X \in C_{0,q}( {\cal V})$, $\Lm \in C_{0,q}(\hat{{\cal V}})$,
$j=1,\ldots,[p]$,
and $i_s \in \{-3,-2,\ldots,n\}$
($1 \le s \le j$), 
it holds that
\begin{eqnarray}
\bigl|
{\cal I}^j
[\hat{Y}^{i_1,(\ve)}, \ldots, \hat{Y}^{i_j,(\ve)} ]_{s,t}
\bigr|
\le 
(\ve + \xi(\ve  X))^{\nu(i_1)+\cdots+\nu(i_j)} \omega(s,t)^{j/p}
\label{ineq.interI_rem}
\end{eqnarray}
\\
\noindent
{\rm (ii)}
For any $r_1, r_2>0$, there exists a constant $c=c(r_1,r_2)>0$ 
depending only on $r_1,r_2,n$
such that
\[
\sup\{ \omega (0,1) ~|~ 
\mbox{ $X \in C_{0,q}( {\cal V})$ with $\xi(\ve X) \le r_1$,
$\Lm \in C_{0,q}(\hat{\cal V})$ with $\|\Lm \|_q \le r_2$}  
\}
\le c.
\]
\end{thm}

\Proof
In this proof the constant $c$ and the control function $\omega$
may change from line to line.
As before we write $\delta=\xi(X)$.
We use induction.
First we consider the case $n=1$.
The estimate 
of the difference of $\Phi^{\ve}(X,\Lm)$
and $\Phi^{\ve}( {\bf 1},\Lm)$ is essentially shown 
in Lemma \ref{lem.geki2}.
The estimate of the difference of $\Phi^{\ve}({\bf 1},\Lm)$
and $\Phi^{0}( {\bf 1},\Lm)$ is a simple exercise
for ODEs in $q$-variational sense.
Thus, combining these, we can easily show the case $n=1$.

Now we assume the statement of the theorem is true for $n$,
and will prove the case for $n+1$. 
First we give estimates for the first term on the right hand side
of (\ref{rep.rem.term}). 
Slightly modifying the definition of $Z^{n,(\ve)}$,
we set 
$\tilde{Z}^{n,(\ve)}
=(X,\Lm,Y^{(\ve)},Y^0,\ldots,Y^{n-1},Q^{n,(\ve)})$.
We also set 
\begin{equation}\label{eq.defFn}
F_n^{(\ve)} (y^{-1},y^0,\ldots,y^{n-1})
=
\ve \sigma(\ve, y^{-1})
-\sum_{k=1}^n \ve^k f_k (y^0,\ldots,y^{k-1}).
\end{equation}
Then, 
$F_n^{(\ve)} \in C^{\infty}_{b,loc} 
(\hat{\cal X}_{n},L(\hat{\cal X}_{n},{\cal W}  ))$
and
$\tilde{F}_n^{(\ve)} :=
h_n^{(\ve)} \oplus  F_n^{(\ve)}
\in C^{\infty}_{b,loc} (\hat{\cal X}_{n},L(\hat{\cal X}_{n},
\hat{\cal X}_{n+1} )$.
Here, $h_n^{(\ve)} \in L(\hat{\cal X}_{n},\hat{\cal X}_{n})$
is defined by
\[
h_n^{(\ve)} (y^{-3},\ldots,y^{n})
=
 (\ve y^{-3},  y^{-2}, y^{-1},y^0, \ve y^1,\ldots,\ve^{n-1}y^{n-1} ,y^{n}).
\]
We now consider $\int \tilde{F}_n^{(\ve)} (\tilde{Z}^{n,(\ve)}) 
d\tilde{Z}^{n,(\ve)}$.

Since it is too complicated to give at once estimates like 
(\ref{defY}) for  all the components of all the level paths of
the above integral,
we first consider (\ref{defY}) 
for the first level path of the last 
component of the above integral, i.e., 
$[\int F_n^{(\ve)} (\tilde{Z}^{n,(\ve)}) 
d\tilde{Z}^{n,(\ve)}]^1$.
Substitute $i=1$,
$f=F_n^{(\ve)}$, and $X=\tilde{Z}^{n,(\ve)}$ in  (\ref{defY}).
Then, 
(the first level of) the almost rough path $\Theta$
which approximates $[\int F_n^{(\ve)} (\tilde{Z}^{n,(\ve)}) 
d\tilde{Z}^{n,(\ve)}]^1$ satisfies that
\begin{eqnarray}
[\Theta^1]^{n+1}
:=
\sum_{l=1}^{[p]} D^{l-1}F_n^{(\ve)} (\tilde{Z}^{n,(\ve)}_s)
\la
[\tilde{Z}^{n,(\ve)}]^l_{s,t}
\ra.
\label{eq.huku}
\end{eqnarray}
Here, $D$ denotes the Fr\'echet derivative on $\hat{\cal X}_{n}$.

Now we estimate the right hand side of (\ref{eq.huku}).
Choose $l$ and fix it. 
The contribution from the first term on the right hand side of
(\ref{eq.defFn}) (i.e., $\ve \sigma(\ve, y^{-1})$) is given as follows;
\begin{eqnarray}
\del_y^{l-1} \sigma(\ve, Y^{(\ve)}_s)
\la  
{\cal I}^{l}[ Y^{(\ve)}, \ldots, Y^{(\ve)}, \ve X]_{s,t}
\ra.
\label{eq.nishi}
\end{eqnarray}
By the Taylor expansion for $\sigma$,
\begin{eqnarray*}
\lefteqn{
\ve  \del_y^{l-1} \sigma(\ve, Y^{(\ve)}_s)
}
\nn\\
&=&
\sum_{j,k;j+k \le n-1} 
\frac{\ve^{j+1}}{j! k!}
\del_y^{l-1+k} \del_{\ve}^{j} \sigma(\ve, Y^{0}_s)
\la  
\overbrace{
Y^{(\ve)}_s-Y^{0}_s, \ldots, Y^{(\ve)}_s-Y^{0}_s
}^{k},
\,\cdot\,
\ra
+
\ve S_{n},
\nn\\
Y^{(\ve)}
&=&
Y^{0}+\ve Y^{1}+\cdots +\ve^{n-1} Y^{n-1}+ Q^{n, (\ve)}.
\end{eqnarray*}
with $|\ve S_{n}| \le 
c (\ve + \ve \delta)^{n+1}$.
Therefore, (\ref{eq.nishi}) is equal to a sum of 
terms of the following form;
\begin{eqnarray}
\frac{\ve^{j}}{j! k!}
\del_y^{l-1+k} \del_{\ve}^{j} \sigma(\ve, Y^{0}_s)
\la  
\hat{Y}^{i_1,(\ve)}_s, \ldots, \hat{Y}^{i_k,(\ve)}_s;
{\cal I}^{l}[  \hat{Y}^{i_{k+1},(\ve)},
\ldots, \hat{Y}^{i_{l-1+k}, (\ve)} , \ve X]_{s,t}
\ra.
\label{eq.huku2}
\end{eqnarray}
Here, $1\le i_1, \ldots, i_k \le n$
and $0 \le i_{k+1}, \ldots, i_{l-1+k} \le n$.
We say this term is of order $1+j+\sum_{a=1}^{l-1+k} \nu(i_a)$,
since this is dominated by
$c (\ve + \ve \delta)^{1+j+\sum_{a=1}^{l-1+k} \nu(i_a)} \omega(s,t)^{l/p}$.
Note that, if a term of this form involves 
$\hat{Y}^{n,(\ve)} =Q^{n,(\ve)}$,
then its order is larger than $n$.

Let $m \le n$.
It is sufficient to show that the terms of order $m$
in (\ref{eq.huku2}) cancel off with 
\begin{eqnarray}
\lefteqn{
\ve^m D^{l-1}f_m^{(\ve)} (\tilde{Z}^{n,(\ve)}_s)
\la
[\tilde{Z}^{n,(\ve)}]^l_{s,t}
\ra
}
\nn\\
&=&
\ve^m
\sum_{ \xi \in \{0,\ldots,m-1 \}^{l-1} }
\de_{\xi}^{l-1}f_m^{(\ve)} 
(  Y^{0}_s, \ldots, Y^{m-1}_s)
\la
[ {\cal I}^l [ Y^{\xi_1}, \ldots, Y^{\xi_{m-1}}  ,X ]_{s,t}
\ra.
\nn
\end{eqnarray}
By Lemma \ref{lem.diffFG} and its proof,
this cancels off with all the terms of order $m$ in (\ref{eq.huku2}).
(Note that 
$l-1, m, k$, and $(i_{k+1}, \ldots, i_{l-1+k})$ in (\ref{eq.huku2})
correspond to
$r, n, k-\mu$, and $(\xi_1, \ldots,\xi_{r})$ in Lemma \ref{lem.diffFG},
respectively.)
Hence, 
$D^{l-1}F_n^{(\ve)} (\tilde{Z}^{n,(\ve)}_s)
\la
[\tilde{Z}^{n,(\ve)}]^l_{s,t}
\ra$
is dominated by 
$c (\ve + \ve \delta )^{n+1} \omega(s,t)^{l/p}$.
Thus, we have obtained an estimate for (\ref{eq.huku}).

Let $1 \le k\le [p]$.
We consider the ${\bf i}=(i_1, \ldots,i_k)$-component of $\Theta^k$, 
where $-3 \le i_j \le n+1$ for all $j=1,\ldots,k$.
Set ${\cal N}({\bf i})=\{ j~|~i_j =n+1\}$.
Note that, if $j \notin {\cal N}({\bf i})$ (equivalently, if $i_j \neq n+1$), 
then ${\bf l}=(l_1,\ldots,l_k)$ in the sum in \ref{defY}
must satisfy $l_j=1$.

We will fix such an $\bf l$.
For $j \notin {\cal N}({\bf i})$,
set $L_j:=\ve^{\nu(i_j)} p_{i_j}$ if $i_j \neq n$
and $L_j:= p_{i_j}$ if $i_j = n$
and also set
${\bf m}_j={\bf m}'_j=i_j$.
For $j \in {\cal N}({\bf i})$, set
${\bf m}_j=(m^1_{j}, \ldots, m^{l_j -1}_{j})$,
${\bf m}'_j=(m^1_{j}, \ldots, m^{l_j -1}_{j},-3)$,
and 
\[
L_j
= D^{l_j -1}F_n^{(\ve)} (\tilde{Z}^{n,(\ve)}_s)
=
\sum_{ {\bf m}_j \in \{0,\ldots,n-1  \}^{l_j -1} }  
\de^{l_j -1}_{ {\bf m}_j}
F_n^{(\ve)} (\tilde{Z}^{n,(\ve)}_s)
\circ
p_{{\bf m}'_j }.
\]

Then, from (\ref{defY}), 
\begin{eqnarray}
[\Theta_{s,t}^k]^{({\bf i})}
&=&
\sum_{{\bf l}} (L_1 \otimes \cdots \otimes L_k )
\bigl\la
\sum_{\pi \in \Pi_{{\bf l}}}
\pi ( \tilde{Z}^{n,(\ve)}  )^{|{\bf l}|}_{s,t} 
\bigr\ra
\nn\\
&=&
\sum_{{\bf l}} (L_1 \otimes \cdots \otimes L_k )
\bigl\la
{\cal I}^k
[  ( \tilde{Z}^{n,(\ve)}  )^{l_1}_{s,\cdot}, 
\ldots,  ( \tilde{Z}^{n,(\ve)}  )^{l_k}_{s,\cdot}
]
\bigr\ra.
\label{eq.tsumeta2}
\end{eqnarray}
The sum is over such ${\bf l}$'s as in (\ref{defY}).

Now we estimate the right hand side of (\ref{eq.tsumeta2}).
Fix ${\bf l}$ for a while.
For $j \notin {\cal N}({\bf i})$, 
Pairing of $L_j$ and $ ( \tilde{Z}^{n,(\ve)}  )^{l_j}$
is clearly of order $\nu(i_j)$.
For $j \in {\cal N}({\bf i})$, 
consider the pairing of $L_j$ and $ ( \tilde{Z}^{n,(\ve)}  )^{l_j}$.
Then, we can see that the same cancellation takes place as in 
(\ref{eq.defFn})--(\ref{eq.huku2})
and that this is of order $n+1$.
If we notice that, for nice paths $\psi^1, \ldots, \psi^{|{\bf l}|}$,
\begin{eqnarray*}
{\cal I}^k
\bigl[ 
 {\cal I}^{l_1}[ \psi^{1}, \ldots, \psi^{l_1}]_{s,\cdot}, 
{\cal I}^{l_2}[ \psi^{l_1+1}, \ldots, \psi^{l_1+l_2}]_{s,\cdot}, 
\ldots,  
{\cal I}^{l_k}[  \psi^{l_1+\cdots +l_{k-1}+1}, \ldots, 
\psi^{|{\bf l}|}]_{s,\cdot}
\bigr]_{s,t}
\nn\\
=
\sum_{\pi \in \Pi_{{\bf l}}}
\pi {\cal I}^{|{\bf l}|}
[  \psi^{\pi(1)} , \ldots, \psi^{\pi(|{\bf l}|)} ]_{s,t},
\end{eqnarray*}
then, by using (\ref{ineq.interI_rem}) for $n$,  
we obtain from the above observation for (\ref{eq.tsumeta2}) that
\begin{eqnarray}
\bigl|
[\Theta_{s,t}^k]^{({\bf i})}
\bigr|
\le
c (\ve+ \ve \delta )^{ \nu(i_1) +\cdots + \nu(i_k)}  \omega(s,t)^{k/p}.
\label{eq.est.theta}
\end{eqnarray}

Next we estimate $R_l(X_s,X_u),~(s<u<t)$ in (\ref{eq.2no5})
with $f=F_n^{(\ve)}$ and $X= \tilde{Z}^{n,(\ve)}$
(since $R_l$ of
other components in $\tilde{F}_n^{(\ve)}$ clearly vanish).
Fix $1 \le l \le [p]$.
From (\ref{eq.defFn}),
the first term in $R_l(F_n^{(\ve)}) 
( \tilde{Z}^{n,(\ve)}_s,\tilde{Z}^{n,(\ve)}_u  )
\la ( \tilde{Z}^{n,(\ve)})^l_{u,t}
\ra$
is given by 
\begin{eqnarray}
\int_0^1 d\theta  \frac{ (1-\theta)^{[p]-l} }{([p]-l)!}
\de_{y}^{[p]} \sigma( \ve, Y^{\ve}_{s,u;\theta})
\la
[( Y^{\ve} )^1_{s,u}]^{\otimes [p]-l+1},
{\cal I}^l
[Y^{(\ve)}, \ldots, Y^{(\ve)},\ve X]_{u,t}
\ra,
\label{eq.huku3}
\end{eqnarray}
where $Y^{\ve}_{s,u;\theta}:= Y^{\ve}_{s}+ \theta ( Y^{\ve} )^1_{s,u}$.
Then, by expanding this as in the previous section,
we can see that the same cancellation takes place as in 
(\ref{eq.defFn})--(\ref{eq.huku2})
and that
\[
|R_l(F_n^{(\ve)}) 
( \tilde{Z}^{n,(\ve)}_s,\tilde{Z}^{n,(\ve)}_u  )
\la ( \tilde{Z}^{n,(\ve)})^l_{u,t}
\ra| 
\le 
c (\ve +\ve \delta)^{n+1} \omega(s,t)^{([p]+1)/p}.
\]
This implies that
\begin{eqnarray*}
|[\Theta^1_{s,t}  -\Theta^1_{s,u}-
\Theta^1_{u,t}]^{n+1}|
&\le&
\sum_{l=1}^{[p]} | R_l(F_n^{(\ve)}) 
( \tilde{Z}^{n,(\ve)}_s,\tilde{Z}^{n,(\ve)}_u  )
\la ( \tilde{Z}^{n,(\ve)})^l_{u,t}
\ra|
\nn\\
&\le&
c (\ve + \ve\delta )^{n+1} \omega(s,t)^{([p]+1)/p},
\quad
s<u<t.
\end{eqnarray*}

From (\ref{eq.2no5})--(\ref{eq.2no7}) and the above estimate for $R_l$, 
we may compute in the same way as in (\ref{eq.tsumeta2})--(\ref{eq.est.theta})
to obtain that
\begin{eqnarray}
\bigl|
[\Theta_{s,t}^k -(\Theta_{s,u}\otimes \Theta_{u,t}  )^k  ]^{({\bf i})}
\bigr|
\le
c (\ve+ \ve \delta )^{ \nu(i_1) +\cdots + \nu(i_k)} 
 \omega(s,t)^{([p]+1)/p}.
\label{eq.est.theta2}
\end{eqnarray}

From (\ref{eq.est.theta}) and (\ref{eq.est.theta2})
we see that the map
\[
(X ,\Lm) \mapsto 
\bigl( \ve X,\Lm,Y^{(\ve)},
Y^0,\ve Y^1, \ldots, \ve^{n-1} Y^{n-1}, Q^{n,(\ve)}, 
\int \sigma(\ve,Y^{(\ve)}) \ve dX
-\sum_{k=1}^n   \ve^k I^k
\bigr)
\]
is continuous and satisfies the inequality (\ref{ineq.interI_rem}) for $n+1$
(with $\hat{Y}^{n,(\ve)}$ and $\hat{Y}^{n+1,(\ve)}$
being replaced with
$Q^{n,(\ve)}$ and $\int \sigma(\ve,Y^{(\ve)}) \ve dX
-\sum_{k=1}^n   \ve^k I^k$, respectively).

By expanding $\int b(\ve, Y^{(\ve)}) d\Lm$ in (\ref{rep.rem.term})
with $Y^{(\ve)}=Y^0+\cdots+\ve^{n-1}Y^{n-1}+Q^{n,(\ve)}$,
we see that
\begin{eqnarray*}
\Bigl| 
\int_s^t \bigl(
b(\ve, Y^{(\ve)}_u) d\Lambda_u
-b(0, Y^{0}_u) d\Lambda_u
-\de_y b(0, Y^0_u) \la Y^{(\ve)}_u - Y^0_u , d \Lambda_u \ra
-\sum_{k=1}^n \ve^k  dJ^k_u \bigr)
\Bigr|
\\
\le c  (\ve+ \ve \delta )^{n+1} 
 \omega(s,t)^{1/q}.
\end{eqnarray*}
(Thanks to the third term 
$\de_y b(0, Y^0_u) \la Y^{(\ve)}_u - Y^0_u , d \Lambda_u \ra$, 
all the terms that involve $Q^{n,(\ve)}$ in the expansion
are of order $n+1$ or larger.)

Let $K^{n+1,(\ve)}$ be the right hand side of (\ref{rep.rem.term}).
From these we see that the map
\[
(X ,\Lm) \mapsto 
\bigl( \ve X,\Lm,Y^{(\ve)},
Y^0,\ve Y^1, \ldots, \ve^{n-1} Y^{n-1}, Q^{n,(\ve)}, K^{n+1,(\ve)}
\bigr)
\]
is continuous and satisfies the inequality (\ref{ineq.interI_rem}) for $n+1$
(with $\hat{Y}^{n,(\ve)}$ and $\hat{Y}^{n+1,(\ve)}$
being replaced with
$Q^{n,(\ve)}$ and $K^{n+1,(\ve)}$, respectively).

By applying Duhamel's principle (Lemma \ref{cor.duhamel}) 
in the same way
as in the proof of Theorem \ref{thm.main.estI} in the previous subsection, 
 we see that the map
\begin{equation}\label{eq.stve}
(X ,\Lm) \mapsto 
\bigl( \ve X,\Lm,Y^{(\ve)},
Y^0,\ve Y^1, \ldots, \ve^{n-1} Y^{n-1}, Q^{n,(\ve)}, Q^{n+1,(\ve)}
\bigr)
\end{equation}
is continuous and satisfies the inequality (\ref{ineq.interI_rem}) for $n+1$
(with $\hat{Y}^{n,(\ve)}$ and $\hat{Y}^{n+1,(\ve)}$
being replaced with
$Q^{n,(\ve)}$ and $Q^{n+1,(\ve)}$, respectively).

Finally define
a bounded linear map $\alpha \in L( \hat{\cal V}_{n+1} ,\hat{\cal V}_{n+1})$
by
\[
\alpha( y^{-3}, \ldots, y^{n}, y^{n+1})
=( y^{-3}, \ldots, -y^{n}+y^{n+1}, y^{n+1})
\]
and apply $\overline{\alpha}$ to (\ref{eq.stve}),
which completes the proof of Theorem \ref{thm.rem.est}.
\QED


%
\section
{Laplace approximation for It\^o functionals of Brownian rough paths}

In this section,
by using the expansion for It\^o maps in the rough path sense
in the previous sections,
we generalize the Laplace approximation
for It\^o functionals of Brownian rough paths,
which was shown in 
Aida \cite{aida} or Inahama and Kawabi \cite{inah-KB-2} (Theorem 3.2).
In this section we always assume $2<p<3$.

\subsection
{Setting of the Laplace approximation}
Let $({\cal V}, {\cal H}, \mu)$ be an abstract Wiener space,
that is, 
${\cal V}$ is a real separable Banach space,
${\cal H}$ is a  real separable Hilbert space
embedded continuously and densely in ${\cal V}$,
and $\mu$ is a Gaussian measure on ${\cal V}$ such that
\[
\int_{{\cal V}}  \exp( \sqrt{-1} \la \phi ,x\ra )   \mu(dx)
=
\exp (-\| \phi \|^2_{{\cal H}^*} /2 ),
\qquad
\text{  for any $\phi \in {\cal V}^*$.}
\]
By the general theory of abstract Wiener spaces, 
there exists a ${\cal V}$-valued Brownian motion $w=(w_t)_{ t \ge 0}$
associated with $\mu$.
The law of the scaled Brownian motion $\ve w$
on $P({\cal V})=\{ y:[0,1]\to {\cal V} | \text{continuous and $y_0=0$}\}$
is denoted by $\hat{\mathbb P}_{\ve}$ ($\ve >0$).

We assume the {exactness condition} {\bf (EX)} below for 
the projective norm on ${\cal V} \otimes {\cal V}$
and $\mu$.
This condition implies the existence of the Brownian rough paths $W$.
(See Ledoux, Lyons, and Qian \cite{llq}.)
The law of the scaled Brownian rough paths $\ve W$
is a probability measure on $G\Omega_p({\cal V})$ 
and is denoted by ${\mathbb P}_{\ve}$ ($\ve >0$).
\vspace{2mm}
\\
{\bf (EX):}
We say that the Gaussian measure $\mu$ 
and the projective norm on $X\otimes X$
satisfies the {\it exactness condition} if 
there exist $C>0$ and $1/2 \le \alpha <1$
such that, for all $n=1,2,\ldots$,
\[
{\mathbb E} \Bigl[ \Bigl|
\sum_{i=1}^{n}  \eta_{2i-1} \otimes \eta_{2i}
 \Bigr|\Bigr]
\le C n^{\alpha}.
\]
Here, $\{\eta_i\}_{i=1}^{\infty}$ are an independent 
and identically distributed  
random variables on $X$
such that the law of $\eta_i$ is $\mu$.
\vspace{2mm}

We consider an ODE in the rough path sense in 
the following form.
Let ${\cal W}$ be another real Banach space.
For $\sigma \in C_b^{\infty} ([0,1] \times {\cal W}, 
L({\cal V},{\cal W}) )$
and $b \in C_b^{\infty} ([0,1] \times {\cal W}, L({\bf R},{\cal W}) )
=
C_b^{\infty} ([0,1] \times {\cal W}, {\cal W})$,
\begin{eqnarray*} 
dY^{(\ve)}_t =\sigma(\ve, Y^{(\ve)}_t) \ve dW +b(\ve, Y^{(\ve)}_t) dt,
\qquad 
Y^{(\ve)}_0 =0.
\end{eqnarray*}
Using the notation of the previous section,
we may write $Y^{(\ve)}= \Phi^{\ve}(\ve W, \Upsilon)$.
Here $\Phi^{\ve}$ is the It\^o map corresponding to $(\sigma, b)$ 
and $\Upsilon$ is the $\bf R$-valued path defined by $\Upsilon_t =t$.

We will study the asymptotics of $(Y^{(\ve)})_1= \Phi^{\ve}(\ve W, \Upsilon)_1$
as $\ve \searrow 0$.
We impose the following conditions 
on the functions $F$ and $G$. In what follows, we especially denote by $D$
the Fr\'echet derivatives on $L_2^{0,1} ({\cal H})$ and $P({\cal W})$.
The Cameron-Martin space for $\hat{\mathbb P}_1$ is denoted by
$L_2^{0,1} ({\cal H})$, which is a linear subspace of $P({\cal W})$.
Note that 
$L_2^{0,1} ({\cal H}) \subset {\rm BV} ({\cal V}) \subset
G\Omega_p ({\cal V})$.
Set $\Psi^0:L_2^{0,1} ({\cal H}) \to P({\cal W})$
by $\Psi^0 (\Lm)= \Phi^{0}(\Lm, \Upsilon)_1$.
\vspace{2mm}
\\
{\bf (H1):}~$F$ and $G$ are real-valued bounded continuous functions defined on $P({\cal W})$.
%
\vspace{2mm}
\\
{\bf (H2):}~The function
$\tilde {F} := F\circ \Psi^{0}+\Vert \cdot
 \Vert_{ L_2^{0,1} ({\cal H}) }^{2}/2$ defined on
$L_2^{0,1} ({\cal H})$ attains its minimum 
at a unique point $\Lm \in L_2^{0,1} ({\cal H})$.
For this $\Lm$, we write $\phi :=\Psi^{0}(\Lm)$. 
%
\vspace{2mm}
\\
{\bf (H3):}~The functions $F$ and $G$ are $n+3$ and $n+1$ 
times Fr\'echet differentiable 
on a neighbourhood $B(\phi)$
of $\phi \in P({\cal W})$, respectively. Moreover there exist positive 
constants $M_{1}, \ldots, M_{n+3}$
such that 
\begin{eqnarray}
\big \vert D^{k}F(\eta)\big[ y,\ldots, y \big ]
\big \vert
&\leq & M_{k} \Vert y \Vert_{P({\cal W})}
^{k}, \quad k=1, \ldots, n+3,
\nonumber \\
\big \vert D^{k}G(\eta)\big[ y,\ldots, y \big ]
\big \vert
&\leq & M_{k} \Vert y \Vert_{P({\cal W})}
^{k}, \quad k=1, \ldots, n+1,
\nonumber
\end{eqnarray}
hold for any $\eta \in B(\phi)$ and $y\in P({\cal W})$.
\vspace{2mm} 
\\
{\bf (H4):}~
At the point $\Lm \in L_2^{0,1} ({\cal H})$, consider the Hessian 
$A:=D^{2}(F\circ \Psi^{0})(\Lm)|_{
L_2^{0,1} ({\cal H})\times L_2^{0,1} ({\cal H})  }$. 
As a bounded self-adjoint operator on $ L_2^{0,1} ({\cal H})$,
the operator $A$ 
is strictly larger than $-
{\rm Id}_{  L_2^{0,1} ({\cal H})}$ in the form sense.
%
\vspace{2mm}

Now we are in a position to state our main theorem. 
This can be considered as a rough path version of 
Azencott \cite{az} or Ben Arous \cite{ben}.
The key of the proof is the stochastic Taylor-like
expansion of the It\^o map around the minimal point $\Lm$,
which will be explained in the next subsection.
(There are many other nice results on this topic 
in the conventional SDE theory.
See Section 5.2 of Pitarbarg and Fatalov \cite{pf}.
For results in the Malliavin calculus, 
see Kusuoka and Stroock \cite{KS1, KS2},
and Takanobu and Watanabe \cite{taka-wata}.)
%
%
%
%
%
\begin{thm}
\label{tenkai-1}
Under conditions {\rm {\bf (EX)}}, {\rm {\bf (H1)}}--{\rm {\bf (H4)}}, 
we have the following asymptotic expansion:
(${\mathbb E}$ is the integration 
with respect to ${\mathbb P}_1$ or $\hat{\mathbb P}_1$.)
\begin{eqnarray}
\lefteqn{
{\mathbb E}\Big [ G(Y^{(\ve)} ) 
\exp \big( -F(Y^{(\ve)})/\ve^{2} \big) \Big]}
\nonumber \\
&=&
\exp \big(- \tilde{F} (\Lm)/\ve^{2}\big)
\exp \big (-c(\Lm)/\ve
\big)\cdot 
\big(\alpha_{0}+\alpha_{1}\ve+ \cdots +\alpha_{n}\ve^{n}+O(\ve^{n+1})\big),
\nonumber \\
\label{thm.expansion}
\end{eqnarray}
where the constant $c(\Lm)$
is given by $c(\Lm):=DF(\phi)[\Xi (\Lm)]$.
Here $\Xi(\Lm) \in P({\cal W})$ is the 
unique solution of the differential equation 
\begin{eqnarray}
d\hspace{0.5mm} 
\Xi_{t}
-\de_{y} \sigma (0,\phi_{t})[\Xi_{t},d\Lm_{t}]
- \de_{y} b(0, \phi_{t})[\Xi_{t}]dt
= 
\de_{\ve}  \sigma(0, \phi_{t})d\Lm_t
+
\de_{\ve}  b(0, \phi_{t})dt 
\nn
\end{eqnarray}
with $\Xi_{0}=0$.
Note that $\Xi$ is non-random.
\end{thm}

\begin{remark}\label{rem.Laplace}
In \cite{inah-KB-2},
only equations of the following form were discussed.
\begin{equation}
dY^{(\ve)}_{t}
={\sigma}( Y^{(\ve)}_{t} )\ve dW_{t}
+\sum_{i=1}^n a_i(\ve) b_i(Y^{(\ve)}_{t} )dt,
\qquad Y^{(\ve)}_{0}=0.
\nn
\end{equation}
Here, $a_i:[0,1]\to {\bf R}$ are ``nice'' functions.
This may be somewhat unnatural.
However,
since we extended the stochastic Taylor-like expansion
to the ``$\ve$-dependent case'' in the previous section,
we are able to slightly generalize the Laplace asymptotics 
(and the large deviation)
as in the above theorem.
\end{remark}

\subsection{Sketch of proof for Theorem \ref{thm.expansion}}
The proof for Theorem \ref{thm.expansion} is essentially 
the same as the one for Theorem 3.2, \cite{inah-KB-2},
once the stochastic Taylor-like expansion is obtained.
Therefore, we only give a sketch of proof in this subsection.

\vspace{2mm}
\noindent
Roughly speaking, there are three steps in the proof:
\\
\noindent
{\bf Step 1:}~A large deviation principal 
for the laws of $Y^{(\ve)}$ as $\ve \searrow 0$.
\\
\noindent
{\bf Step 2:}~The stochastic Taylor expansion 
around the maximal point.
We expand $\hat{Y}^{(\ve)}$  
as $\ve \searrow 0$ as in the previous sections, 
where $\hat{Y}^{(\ve)}$
is given by the following differential equation:
\begin{equation}
d\hat{Y}^{(\ve)}_t
=
{\sigma}(\ve, \hat{Y}^{(\ve)}_{t} ) ( \ve dW_{t} +d\Lm_t)
+
 b(\ve, \hat{Y}^{(\ve)}_{t} )dt,
\qquad \hat{Y}^{(\ve)}_{0}=0.
\label{eq.trans.ode}
\end{equation}
Here, $\Lm\in L_2^{0,1} ({\cal H})$ is given in Assumption {\rm {\bf (H2)}}.
Note that $\phi =\hat{Y}^0$ if we use the notation in the previous section.
\\
\noindent
{\bf Step 3:}~
Combine the expansion for $\hat{Y}^{ (\ve)}$
with the Taylor expansion for $F,G,$ and $\exp$.

\vspace{2mm}

Firstly, we explain  Step 1.
We use the large deviation for Brownian rough paths
(Theorem 1 in Ledoux, Qian, and Zhang \cite{lqz}. 
The infinite dimensional case is in \cite{inah-KB})
and then use the contraction principle of It\^o map,
which is continuous.
This strategy was established in \cite{lqz}.

\begin{prop}
The law of $Y^{(\ve)}$ on $P({\cal W})$
satisfies a large deviation principle as $\ve \searrow 0$
with the following rate function $I$:
\begin{eqnarray*}
I(y) = \left\{
    \begin{array}{ll}
  \inf
\{  \|X\|^2_{L_2^{0,1} ({\cal H}) }  /2  ~|~ y=\Psi^{0}(X) \}  &  
\text{ if $y=\Phi^{0}(X)$  for some $X \in L_2^{0,1} ({\cal H})$,}  \\
   \infty &  \text{ otherwise.} 
    \end{array}
\right.
\end{eqnarray*}
\end{prop}

\Proof
First recall that ${\mathbb P}_{\ve}$ on $G\Omega_p ({\cal V})$
satisfies a large deviation principle as $\ve \searrow 0$
with the following rate function $J$
(see Theorem 1, \cite{lqz} or Theorem 3.2, \cite{inah-KB}):
\begin{eqnarray*}
J( X ) = \left\{
    \begin{array}{ll}
  
\|X\|^2_{L_2^{0,1} ({\cal H}) }  /2  
&  
\text{ if  $X \in L_2^{0,1} ({\cal H})$,}  
\\
   \infty &  \text{ otherwise.} 
    \end{array}
\right.
\end{eqnarray*}
Then, from the slight extension of Lyons' continuity theorem 
(Theorem \ref{cor.loc.Lip})
and the  slight extension of the contraction principle
(Lemma 3.9, \cite{inah-KB}, for instance),
we can prove the proposition.
\QED

Secondly, we explain  Step 2.
We can use Theorems \ref{thm.main.estI} and \ref{thm.rem.est},
if we set ${\cal V}={\cal V}$, $\hat{\cal V}={\cal V}\oplus {\bf R}$
and regard the equation (\ref{eq.trans.ode}) as follows:
\[
d\hat{Y}^{(\ve) }
=
{\sigma}(\ve, \hat{Y}^{(\ve)}_{t} ) \ve dW_{t} 
+\bigl[
 \sigma (\ve, \hat{Y}^{(\ve)}_{t} ) d\Lm_t
+
 b(\ve, \hat{Y}^{(\ve)}_{t} )dt
\bigr],
\qquad \hat{Y}^{(\ve)}_{0}=0.
\]

Finally, we explain  Step 3.
This step is essentially the same as in Section 6, \cite{inah-KB}.
In this step, a Fernique type theorem 
and a Cameron-Martin type theorem 
for the Brownian rough paths
are used.
(See, for instance, Theorem 2.2 and Lemma 2.3, \cite{inah}.)



\begin{thebibliography}{99}

\bibitem{aida2}
S. Aida: {\it 
Notes on proofs of continuity theorem 
in rough path analysis}, preprint, 2006.
%
\bibitem{aida}
S. Aida: {\it Semi-classical limit of the bottom of spectrum of a 
Schr\"odinger operator on a path space over a compact Riemannian
manifold}, to appear in J. Funct. Anal.
%
%
%
\bibitem{az}
R. Azencott:
{\it Formule de Talyor stochastique et d\'{e}veloppement asymptotique
d'int\'{e}grales de Feynman}, in \lq \lq Seminar on Probability, XVI,
Supplement \rq \rq,
pp. 237--285, Lecture Notes in Math., {\bf 921}, 
Springer, Berlin-New York, 1982.
%


\bibitem{bc}
F. Baudoin and L. Coutin,
{\it Operators associated with a stochastic 
differential equation driven by fractional Brownian motions,}
Stochastic Processes and their Applications, 
{\bf 117}, (2007), 5, pp. 550-574



\bibitem{ben} G. Ben Arous:
{\it Methods de Laplace et de la phase stationnaire sur l'espace de Wiener},
Stochastics {\bf 25}, (1988), no.3, pp. 125--153.
%
%
%
%


\bibitem{cq}
L. Coutin and Z. Qian:
{\it Stochastic analysis, rough path analysis and fractional 
Brownian motions.} 
Probab. Theory Related Fields 122 (2002), no. 1, 108--140. 



\bibitem{cfv}
L. Coutin, P. Friz and N. Victoir:
{\it Good Approximations to rough paths and applications.}
Annals of Probability 35 (2007), no. 3, 1172--1193. 

\bibitem{du}
J. Diestel and J. J. Uhl Jr:
{\it Vector measures,} 
Mathematical Surveys, No. 15. American Mathematical Society, P
rovidence, R.I., 1977.


\bibitem{fv}
P. Friz and N. Victoir:
{\it Euler estimates for rough differential equations. }
J. Differential Equations 244 (2008), no. 2, 388--412.



\bibitem{fv2}
P. Friz and N. Victoir:
Differential equations driven by Gaussian signals.
To appear in Annales de l'Institut Henri Poincar\'e (B) Probability and Statistics (2009).



%
\bibitem{inah}
Y. Inahama:
{\it{Laplace's method for the laws of heat processes on loop spaces}},
J. Funct. Anal. {\bf 232} (2006), pp. 148--194.
%


\bibitem{inah2}
Y. Inahama:
{\it  Laplace approximation for rough differential equation driven by fractional Brownian motion},
submitted (2009).


\bibitem{inah-KB}
Y. Inahama and H. Kawabi:
{\it{Large deviations for heat kernel measures on loop spaces 
via rough paths}},
J. London Math. Society {\bf 73} (2006), no. 3, pp. 797--816.
%
\bibitem{inah-KB-1.5}
Y. Inahama and H. Kawabi:
{\it{On asymptotics of certain Banach space-valued It\^o functionals
of Brownian rough paths}},
 in Stochastic Analysis and Applications; The Abel Symposium 2005, 
pp. 415--434,
Stockholm,
Springer, 2007.
%


\bibitem{inah-KB-2}
Y. Inahama and H. Kawabi:
{\it Asymptotic expansions for the Laplace approximations
for It\^{o} functionals of Brownian rough paths,}
J. Funct. Anal. {\bf 243} (2007), no. 1, 270--322.


%
\bibitem{KS1}
S. Kusuoka and D.W. Stroock:
{\it Precise asymptotics of certain Wiener functionals},
J. Funct. Anal. {\bf 99} (1991), no. 1, pp. 1--74.

\bibitem{KS2}
S. Kusuoka and D.W. Stroock:
{\it Asymptotics of certain Wiener functionals with degenerate extrema},
Comm. Pure Appl. Math. {\bf 47} (1994), no. 4, pp. 477--501.
%
%
\bibitem{llq}
M. Ledoux, T. Lyons and Z. Qian:
{\it
L\'{e}vy area of Wiener processes in Banach spaces}, 
Ann. Probab. {\bf 30} (2002), no. 2, pp. 546--578. 
%
\bibitem{lqz} 
M. Ledoux, Z. Qian and T.S. Zhang:
{\it
Large deviations and support theorem for diffusion processes via rough paths}, 
Stochastic Process. Appl. {\bf 102} (2002), no. 2, pp.265--283. 
%
\bibitem{lq}
T. Lyons and Z. Qian:
{\rm System control and rough paths}, 
Oxford University Press, Oxford, 2002. 
%
%
\bibitem{ms}
A. Millet, M. Sanz-Sol\'e:
{\it Large deviations for rough paths of the fractional Brownian motion.} 
Ann. Inst. H. Poincare Probab. Statist. 42 (2006), no. 2, 245--271. 


\bibitem{pf}
V. I. Piterbarg, V. R. Fatalov:
{\it The Laplace method for probability measures in Banach spaces.} 
English translation in Russian Math. Surveys 50 (1995), no. 6, 1151--1239. 


%
\bibitem{taka-wata}
S. Takanobu and S. Watanabe:
{\it 
Asymptotic expansion formulas of the Schilder type for a class of conditional
Wiener functional integrations}, in 
\lq \lq
Asymptotic problems in probability theory: Wiener functionals and asymptotics 
\rq \rq
(Sanda/Kyoto, 1990), pp. 194--241, 
Pitman Res. Notes Math. Ser., {\bf 284}, 
Longman Sci. Tech., Harlow, 1993. 

\end{thebibliography}
\end{document}